\documentclass{amsart}

\usepackage{amsthm}
\usepackage{amsfonts}
\usepackage{amsmath}
\usepackage{amssymb}
\usepackage{mathrsfs}
\usepackage{fullpage}
\usepackage{stmaryrd}
\usepackage{hyperref}
\usepackage{verbatim}

\theoremstyle{plain}

\theoremstyle{definition}

\theoremstyle{remark}

\newcommand{\Begin}[2]{\begin{#1}\label{#2}}

\newcommand{\Z}{\mathbb{Z}}

\newcommand{\bSigma}{\mathbf{\Sigma}}
\newcommand{\bDelta}{\mathbf{\Delta}}
\newcommand{\bbC}{\mathbb{C}}

\newcommand{\bbP}{\mathbb{P}}
\newcommand{\bbQ}{\mathbb{Q}}
\newcommand{\bbR}{\mathbb{R}}
\newcommand{\bbS}{\mathbb{S}}

\newcommand{\forces}{\Vdash}
\newcommand{\analytic}{{\bSigma_1^1}}
\newcommand{\lanalytic}{{\Sigma_1^1}}

\newcommand{\lcoanalytic}{{\Pi_1^1}}
\newcommand{\borel}{{\bDelta_1^1}}
\newcommand{\lborel}{{\Delta_1^1}}
\newcommand{\cantorspace}{{{}^\omega 2}}
\newcommand{\bairespace}{{{}^\omega\omega}}
\newcommand{\finBinarySequence}{{{}^{<\omega}2}}

\newcommand{\reals}{\bbR}

\newcommand{\pcantorspace}{{{}^\omega(\cantorspace)}}
\newcommand{\group}{\bigoplus_{n \in \omega} \Z_2}
\newcommand{\wgroup}{{}^\omega(\group)}
\newcommand{\supp}{\mathrm{supp}}
\newcommand{\dom}{\mathrm{dom}}
\newcommand{\grid}{\mathrm{grid}}
\newcommand{\splitt}{\mathrm{split}}
\newcommand{\hatPE}{\widehat{\bbP}_{E_0}}
\newcommand{\PE}{{\bbP_{E_0}}}
\newcommand{\prune}{\mathsf{prune}}
\newcommand{\tprune}{\mathsf{2prune}}
\newcommand{\fprune}{\mathsf{3prune}}
\newcommand{\switch}{\mathsf{switch}}
\newcommand{\Etail}{E_\mathrm{tail}}
\newcommand{\Gal}{\mathrm{Gal}}

\begin{document}

\title{Definable Combinatorics of Some Borel Equivalence Relations}

\author{William Chan}
\address{Department of Mathematics, University of North Texas, Denton, TX 76203}
\email{William.Chan@unt.edu}

\author{Connor Meehan}
\address{Department of Mathematics, California Institute of Technology, Pasadena, CA 91125}
\email{cgmeehan@caltech.edu}

\begin{abstract}
If $X$ is a set, $E$ is an equivalence relation on $X$, and $n \in \omega$, then define 
$$[X]^n_E = \{(x_0, ..., x_{n - 1}) \in {}^nX : (\forall i,j)(i \neq j \Rightarrow \neg(x_i \ E \ x_j))\}.$$ 

For $n \in \omega$, a set $X$ has the $n$-J\'onsson property if and only if for every function $f : [X]^n_= \rightarrow X$, there exists some $Y \subseteq X$ with $X$ and $Y$ in bijection so that $f[[Y]^n_=] \neq X$. A set $X$ has the J\'onsson property if and only for every function $f : (\bigcup_{n \in \omega}[X]^n_=) \rightarrow X$, there exists some $Y \subseteq X$ with $X$ and $Y$ in bijection so that $f[\bigcup_{n \in \omega} [Y]^n_=] \neq X$. 

Let $n \in \omega$, $X$ be a Polish space, and $E$ be an equivalence relation on $X$. $E$ has the $n$-Mycielski property if and only if for all comeager $C \subseteq {}^nX$, there is some $\borel$ $A \subseteq X$ so that $E \leq_\borel E \upharpoonright A$ and $[A]^n_E \subseteq C$. 

The following equivalence relations will be considered: $E_0$ is defined on $\cantorspace$ by $x \ E_0 \ y$ if and only if $(\exists n)(\forall k > n)(x(k) = y(k))$. $E_1$ is defined on $\pcantorspace$ by $x \ E_1 \ y$ if and only if $(\exists n)(\forall k > n)(x(k) = y(k))$. $E_2$ is defined on $\cantorspace$ by $x \ E_2 \ y$ if and only if $\sum\{\frac{1}{n + 1} : n \in x \ \triangle \ y\} < \infty$, where $\triangle$ denotes the symmetric difference. $E_3$ is defined on $\pcantorspace$ by $x \ E_3 \ y$ if and only if $(\forall n)(x(n) \ E_0 \ y(n))$. 

Holshouser and Jackson have shown that $\bbR$ is J\'onsson under $\mathsf{AD}$. It will be shown that $E_0$ does not have the $3$-Mycielski property and that $E_1$, $E_2$, and $E_3$ do not have the $2$-Mycielski property. Under $\mathsf{ZF + AD}$, $\cantorspace \slash E_0$ does not have the $3$-J\'onsson property. 
\end{abstract}

\maketitle\let\thefootnote\relax\footnote{September 7, 2017

The authors were partially supported by NSF grant DMS-1464475. The first author was also supported by NSF grant DMS-1703708.}

%\tableofcontents

\section{Introduction}\label{Introduction}
The J\'onsson property and other combinatorial partition properties of well-ordered sets have been studied by set theorists under the axiom of choice, large cardinal axioms, and the axiom of determinacy. Holshouser and Jackson began the study of the J\'onsson property using definability techniques for sets which generally cannot be well-ordered in a definable manner.

Let $X$ be a set and $E$ an equivalence relation on $X$. For each $n \in \omega$, let $[X]^n_E$ be the collection of tuples $(x_0, ..., x_{n - 1}) \in {}^nX$ so that for all $i \neq j$, $\neg(x_i \  E \ x_j)$. Let $[X]^{<\omega}_E = \bigcup_{n \in \omega} [X]_E^n$. For each $n \in \omega$, $X$ has the \textit{$n$-J\'onsson property} if and only if for every function $f : [X]^n_= \rightarrow X$, there is some $Y \subseteq X$ with $Y$ in bijection with $X$ and $f[[Y]^n_=] \neq X$. $X$ has the \textit{J\'onsson property} if and only if for every function $f : [X]^{<\omega}_= \rightarrow X$, there is some $Y \subseteq X$ with $Y$ in bijection with $X$, and $f[[Y]^{<\omega}_=] \neq X$. 

Holshouser and Jackson showed that $\cantorspace$ has the J\'onsson property under the axiom of determinacy, $\mathsf{AD}$. Let $f : [\cantorspace]^{<\omega}_= \rightarrow \cantorspace$. For each $n \in \omega$, let $f_n : [X]^n_= \rightarrow X$ be $f \upharpoonright [X]^n_=$. Their proof has two notable tasks: 

\noindent(1) Holshouser and Jackson first (assuming all sets have the Baire property) choose comeager sets $C_n \subseteq {}^n(\cantorspace)$ so that $f_n \upharpoonright C_n$ is continuous. Then a single perfect set $P \subseteq \cantorspace$ is found so that for each $n$, $f_n \upharpoonright [P]^n_=$ is continuous. To obtain this perfect set $P$, they use a classical theorem of Mycielski which states: If $C_n$ is a sequence of comeager subsets of ${}^n(\cantorspace)$, then there is some perfect set $P \subseteq \cantorspace$ so that $[P]^n_= \subseteq C_n$ for all $n$. 

\noindent(2) Since each $f_n$ is continuous on $[P]^n_=$, they use a fusion argument to simultaneously prune $P$ to a smaller perfect set $Q \subseteq P$ so that there exists some real that is missed by each $f_n$ on $[Q]^n_=$. 

Holshouser and Jackson ask whether other sets which may not be well-ordered in some choiceless setting like $\mathsf{AD}$ could also have the J\'onsson property. They observed that under $\mathsf{ZF + AD + V = L(\bbR)}$, every set $X \in L_\Theta(\bbR)$ has a surjective function $f : \bbR \rightarrow X$. Define an equivalence relation on $\bbR$ by $x \ E \ y$ if and only if $f(x) = f(y)$. Then $X$ is in bijection with $\bbR \slash E$. The study of the J\'onsson property for sets in $L_\Theta(\bbR)$ is equivalent to studying the J\'onsson property for quotients of $\bbR$ by equivalence relations on $\bbR$. Note that $\bbR$ is in bijection with $\bbR \slash =$. 

Through dichotomy results of Harrington, Hjorth, Kechris, Louveau, and others, the equivalence relations $=$, $E_0$, $E_1$, $E_2$, and $E_3$ occupy special positions in the structure of $\borel$ equivalence relations under $\borel$ reducibilities. $=$ is the identity equivalence relation on $\cantorspace$. $E_0$ is defined on $\cantorspace$ by $x \ E_0 \ y$ if and only if $(\exists n)(\forall k > n)(x(k) = y(k))$. $E_1$ is defined on $\pcantorspace$ by $x \ E_1 \ y$ if and only if $(\exists n)(\forall k > n)(x(k) = y(k))$. $E_2$ is defined on $\cantorspace$ by $x \ E_2 \ y$ if and only if $\sum\{\frac{1}{n + 1} : n \in x \ \triangle \ y\} < \infty$, where $\triangle$ denotes the symmetric difference. $E_3$ is defined on $\pcantorspace$ by $x \ E_3 \ y$ if and only if $(\forall n)(x(n) \ E_0 \ y(n))$. 

Holshouser and Jackson asked whether the methods applicable for showing $\bbR \slash = $ has the J\'onsson property could be used to show the quotients of these other $\borel$ equivalence relations could be J\'onsson. An important aspect of their proof for $\bbR$ was the theorem of Mycielski. They defined the Mycielski property for arbitrary equivalence relations as follows: Let $E$ be an equivalence relation on a Polish space $X$. For each $n \in \omega$, $E$ has the \textit{$n$-Mycielski property} if and only if for every comeager $C \subseteq {}^nX$, there exists some $\borel$ $A \subseteq X$ so that $E \leq_\borel E \upharpoonright A$ and $[A]^n_{E} \subseteq C$. 

They asked whether any of the $\borel$ equivalence relations mentioned above have the $n$-Mycielski property for various $n \in \omega$ and whether the Mycielski property could be used to prove the J\'onsson property for the quotient of any of these equivalence relations. Holshouser and Jackson began this study by showing that $E_0$ has the $2$-Mycielski property and this can be used to show $\cantorspace \slash E_0$ has the $2$-J\'onsson property. This paper will show that the Mycielski property fails in most cases:
\newline
\newline\noindent \textbf{Theorem \ref{E0 does not have Mycielski property}.} The equivalence relation $E_0$ does not have the $3$-Mycielski property.
\newline
\newline\noindent \textbf{Theorem \ref{E1 does not have Mycielski property}.} The equivalence relation $E_1$ does not have the $2$-Mycielski property.
\newline
\newline\noindent\textbf{Theorem \ref{E2 does not have Mycielski property}.} The equivalence relation $E_2$ does not have the $2$-Mycielski property.
\newline
\newline\noindent\textbf{Theorem \ref{E3 does not have mycielski property}.} The equivalence relation $E_3$ does not have the $2$-Mycielski property.
\newline
\newline\indent These results require understanding the structure of $\borel$ subsets of $\cantorspace$ or $\pcantorspace$ so that $E_0 \leq_\borel E_0 \upharpoonright A$ (or $E_1 \leq_\borel E_1 \upharpoonright A$, etc.) that come from the proofs of the dichotomy results. Kanovei, Sabok, and Zapletal in \cite{Borel-Equivalence-Relations}, \cite{Canonical-Ramsey-Theory-on-Polish-Spaces}, \cite{Descriptive-Set-Theory-and-Definable-Forcing}, and \cite{Forcing-Idealized} have studied the forcing of such $\borel$ sets for each of these equivalence relations.

Given that the Mycielski property fails in general, a reflection on Holshouser and Jackson's proof of the J\'onsson property for $\cantorspace$ shows that it is only used to find some perfect set $P$ so that $f_n \upharpoonright [P]^n_=$ is nicely behaved (i.e., continuous). This paper will give a forcing style proof of Holshouser and Jackson results that $\cantorspace$ is J\'onsson and $\cantorspace \slash E_0$ is $2$-J\'onsson assuming all functions satisfy a certain definability condition expressed in Lemma \ref{absoluteness implies R-jonsson property}. This definability condition follows from the Mycielski property for the equivalence relation and the assumption that all sets have the Baire property. All $\borel$ functions have this definability condition and under the axiom of choice and large cardinal assumptions, projective and even more complex sets also satisfy this condition. 

Following part (2) of Holshouser and Jackson's template for $\cantorspace$, suppose one could find some $\borel$ set $B \subseteq \cantorspace$ with $E_0 \leq_\borel E_0 \upharpoonright B$ and $f \upharpoonright [B]^3_{E_0}$ is continuous for some function $f$. Could one then somehow prune $B$ to some $C \subseteq B$ so that $E_0 \leq_\borel E_0 \upharpoonright C$ and $f[[C]^3_{E_0}] \neq \cantorspace$, or even better, miss an $E_0$-class? This paper will have some discussion on how these continuity and surjectivity properties for $E_0$ and $E_2$ can fail.

This shows both part (1) and part (2) of the proof of Holshouser and Jackson establishing $\cantorspace$ is J\'onsson fail for $E_0$ and several other $\borel$ equivalence relations. Moreover, for $E_0$, it will in fact be shown that $\cantorspace \slash E_0$ is not J\'onsson under determinacy:
\newline
\newline\noindent\textbf{Theorem \ref{E0 does not have 3-Jonsson property}.} $(\mathsf{ZF + AD})$ $\cantorspace \slash E_0$ does not have the $3$-J\'onsson property and hence is not J\'onsson.
\newline
\newline\indent Here are some historical remarks about the J\'onsson property: Under the axiom of choice, the J\'onsson property is usually studied on cardinals. Cardinals possessing the J\'onsson property are called \textit{J\'onsson cardinals}. For $n \in \omega$, let $\mathscr{P}^n(X)$ denote the collection of all $n$-element subsets of $X$. Since there is a well-ordering, the J\'onsson property is usually defined using $\mathscr{P}^n(X)$ rather than $[X]^n_=$. When this paper discusses the J\'onsson property using $\mathscr{P}^n(X)$, it will be refered to as the classical J\'onsson property. 

Under the axiom of choice, J\'onsson cardinals also have model-theoretic characterizations. The existence of J\'onsson cardinals imply $V \neq L$. Moreover, it has large cardinal consistency strength: for instance, it implies $0^\sharp$ exists. Erd\H{o}s and Hajnal (\cite{On-a-Problem-of-B-Jonsson} and \cite{Some-Weak-Versions-of-Large-Cardinal-Axioms}) showed that if $2^\kappa = \kappa^+$, then $\kappa^+$ is not a J\'onsson cardinal. Hence under $\mathsf{CH}$, $2^{\aleph_0}$ is not a J\'onsson cardinal. Every real valued measurable cardinal is J\'onsson (see \cite{Some-Weak-Versions-of-Large-Cardinal-Axioms} Corollary 11.1). Solovay showed the consistency of a measurable cardinal implies the consistency of $2^{\aleph_0}$ being real valued measurable. Hence it is consistent relative to a measurable cardinal that $2^{\aleph_0}$ is a J\'onsson cardinal. The sets $\cantorspace$, $\cantorspace \slash E_0$, $\pcantorspace \slash E_1$, $\cantorspace \slash E_2$, $\pcantorspace \slash E_3$ are all in bijection with each other using the axiom of choice. Hence if $\mathsf{CH}$ holds, these quotients do not have the J\'onsson property and if $2^{\aleph_0}$ is real valued measurable, then all these quotients do have the J\'onsson property.

Under $\mathsf{AD}$, the J\'onsson property and other combinatorial partition properties of cardinals were already studied during the 1960s and 1970s. Assuming $\mathsf{AD}$, for each $n \in \omega$, $\aleph_n$ is a J\'onsson cardinal (\cite{Infinitary-Combinatorics-and-the-Axiom-of-Determinateness}). More recently Woodin had shown that under $\mathsf{ZF + AD^+}$, every cardinal $\kappa < \Theta$ has the J\'onsson property. Also \cite{Determinacy-and-Jonsson-Cardinals-in-LR} showed that in $\mathsf{ZF + AD + V = L(\bbR)}$, every cardinal $\kappa < \Theta$ is J\'onsson. \cite{Determinacy-and-Jonsson-Cardinals-in-LR} asked whether $\cantorspace$, which cannot be well-ordered, has the J\'onsson property. In analogy, they asked if every set in $L_\Theta(\reals)$ has the J\'onsson property. Holshouser and Jackson's answer to this question for $\cantorspace$ begins the work that is carried out in this paper.

Throughout, results attributed to Holshouser and Jackson can be found in \cite{Partition-Properties-for-Hyperfinite-Quotients} and \cite{Partition-Properties-for-Non-Ordinal-Sets}.
\newline
\newline\indent This paper is organized as follows:

Section \ref{Basic Information} contains definitions of the main concepts and some basic facts about determinacy. 

Section \ref{reals have the jonsson property} will give a proof of the result of Holshouser and Jackson which shows $\cantorspace$ has the J\'onsson property if all sets have the Baire property. The proof uses forcing arguments and fusion. This section will have some discussions about how absoluteness available under $\mathsf{AD}^+$ can be used to prove this result without using the Mycielski property. However, throughout the paper, a flexible fusion argument is necessary for handling the combinatorics. It is unclear what the relation is between properness, fusion, and the J\'onsson property for the five equivalence relations considered.

Upon considering the J\'onsson property for $\cantorspace$, a natural question is whether there is a function $f : \mathscr{P}^\omega(\cantorspace) \rightarrow \cantorspace$ so that for all $A \subseteq \cantorspace$ with $A \approx \cantorspace$, $f[\mathscr{P}^\omega(A)] = \cantorspace$. Such a function is called an \textit{$\omega$-J\'onsson function} for $\cantorspace$. Under the axiom of choice, \cite{On-a-Problem-of-B-Jonsson} showed that every set has an $\omega$-J\'onsson function. Section \ref{w-jonsson function for reals} gives an example under $\mathsf{ZF + AC^\bbR_\omega}$ (choice for countable sets of nonempty subsets of $\cantorspace$) of a $\borel$ $\omega$-J\'onsson function for $\cantorspace$.

From the effective proof of the $E_0$-dichotomy, every $\analytic$ set $A \subseteq \cantorspace$ so that $E_0 \leq_\borel E_0 \upharpoonright A$ contains the body of a perfect tree with certain symmetry restrictions, known as an $E_0$-tree. Section \ref{the structure of E0} will modify the proof of the $E_0$-dichotomy using Gandy-Harrington methods to prove a structure theorem for $\analytic$ sets with the same $E_0$-saturation on which the restriction of $E_0$ is not smooth: For example, if $A$ and $B$ are two $\analytic$ sets with $E_0 \leq_\borel E_0 \upharpoonright A$, $E_0 \leq_\borel E_0 \upharpoonright B$, and $[A]_{E_0} = [B]_{E_0}$, then there are $E_0$ trees $p$ and $q$ with $[p] \subseteq A$, $[q] \subseteq B$, and $p$ and $q$ are the same except possibly at the stem. This is needed to show the failure of the weak $3$-Mycielski property (see Definition \ref{weak Mycielski property}) for $E_0$. 

Section \ref{R/E0 has 2-Jonsson property} will introduce the forcing $\hatPE^2$. This forcing will be used to prove the result of Holshouser and Jackson stating that $\cantorspace \slash E_0$ has the $2$-J\'onsson property.

Let $X$ and $Y$ be sets. Let $n \in \omega$. Define $X \rightarrow (X)^n_Y$ to mean that for any function $f : \mathscr{P}^n(X) \rightarrow Y$, there is some $Z \subseteq X$ with $Z \approx X$ and $|f[\mathscr{P}^n(Z)]| = 1$. Define $X \mapsto (X)^n_Y$ to mean that for any function $f : [X]^n_= \rightarrow Y$, there is some $Z \subseteq X$ with $Z \approx X$ and $|f[[Z]^n_=]| = 1$. Section \ref{Partition Properties of R/E0 in dimension 2} will show that $\cantorspace \slash E_0 \mapsto (\cantorspace \slash E_0)^2_{n}$ holds for all $n \in \omega$. 

Section \ref{E0 does not have 3-Mycielski property} will show that $E_0$ does not have the $3$-Mycielski property or weak $3$-Mycielski property.

Section \ref{surjectivity and continuity aspects of E0} will produce a continuous function $Q: [\cantorspace]_{E_0}^3 \rightarrow \cantorspace$ so that for every $\analytic$ $A \subseteq \cantorspace$ with $E_0 \leq_\borel E_0 \upharpoonright A$, $Q[[A]^3_{E_0}] = \cantorspace$. A modification of this function yields a $\borel$ function $K : {}^3(\cantorspace) \rightarrow \cantorspace$ so that on any such set $A$, $K \upharpoonright [A]^3_{E_0}$ is not continuous. 

Section \ref{R/E0 does not have 3-Jonsson property} will use the function produced in the previous section to show that $\cantorspace \slash E_0$ does not have the $3$-J\'onsson property or the classical $3$-J\'onsson property under $\mathsf{ZF + AD}$. (In particular, $\cantorspace \slash E_0$ is not J\'onsson under $\mathsf{AD}$.)

Section \ref{Failure of Partition Properties of R/E0 in Dimension Higher than 2} will use the classical $3$-J\'onsson map for $\cantorspace \slash E_0$ to show the failure of $\cantorspace \slash E_0 \rightarrow (\cantorspace \slash E_0)^3_2$. 

The fusion argument related to the proper forcing $\hatPE^2$ was used to establish many of the combinatorial properties of $E_0$ in dimension two. Given the failure of these properties in dimension three, a natural question would be whether the three dimensional analog $\hatPE^3$ is proper and possesses a reasonable fusion. Section \ref{hatPE3 is proper} will show that $\hatPE^3$ is proper by having some type of fusion argument. However, there is far less control of this fusion.

Section \ref{E1 Does Not Have the 2-Mycielski Property} will show that $E_1$ does not have the $2$-Mycielski property. 

Section \ref{The Structure of E2} will modify the proof of the $E_2$-dichotomy result using Gandy-Harrington methods to give structural result about $E_2$-big $\analytic$ sets with the same $E_2$-saturation: For example, if $A$ and $B$ are two $\analytic$ sets with $E_2 \leq_\borel E_2 \upharpoonright A$, $E_2 \leq_\borel E_2 \upharpoonright B$, and $[A]_{E_2} = [B]_{E_2}$, then there are two $E_2$-trees (perfect trees with certain properties) $p$ and $q$ so that $[p] \subseteq A$, $[q] \subseteq B$, $[[p]]_{E_2} = [[q]]_{E_2}$, and $p$ and $q$ resemble each other in specific ways.

Section \ref{E2 Does Not Have the 2-Mycielski Property} will use results of the previous section to show $E_2$ does not have the $2$-Mycielski property and the weak $2$-Mycielski property.

Section \ref{Surjectivity and Continuity Aspects of E2} will produce a continuous function $Q : [\cantorspace]^3_{E_2} \rightarrow \cantorspace$ so that on any $\analytic$ set $A$ with $E_2 \leq_\borel E_2 \upharpoonright A$, $Q[[A]_{E_2}^3] = \cantorspace$. There is also a $\borel$ function $P' : {}^3(\cantorspace) \rightarrow \cantorspace$ so that for any such set $A$, $P' \upharpoonright A$ is not continuous.

Section \ref{The Structure of E3} contains no new results but just gives the rather lengthy characterization of $\analytic$ sets $A \subseteq \pcantorspace$ so that $E_3 \leq_\borel E_3 \upharpoonright A$ which comes from the $E_3$-dichotomy result. This structure result is applied in Section \ref{E3 Does Not Have the 2-Mycielski Property} to show that $E_3$ does not have the $2$-Mycielski property.

Section \ref{Completeness of Ultrafilters on Quotients} will study the completeness of non-principal ultrafilters on quotients of Polish spaces by equivalence relations. 
\\*
\\*\indent The authors would like to thank Jared Holshouser and Stephen Jackson whose talks and subsequent discussions motivated the work that appears here. The authors would also like to thank Alexander Kechris for comments and discussions about this paper. 

\section{Basic Information}\label{Basic Information}
\Begin{definition}{string notation}
Let $\sigma \in \finBinarySequence$. Suppose $|\sigma| = k$. Then $\tilde \sigma \in \cantorspace$ is defined by $\tilde\sigma(n) = \sigma(j)$ where $0 \leq j < k$ and $j \equiv n \ \mathrm{mod} \ k$.

For example, $\tilde 0$, $\tilde 1$, $\widetilde{01}$, etc. will appear frequently.
\end{definition}

\Begin{definition}{basic open neighborhood cantorspace}
Let $\sigma \in \finBinarySequence$. Let $N_\sigma = \{x \in \cantorspace : x \supseteq \sigma\}$. 

$\{N_\sigma : \sigma \in \finBinarySequence\}$ is a basis for the topology on $\cantorspace$. 

Let $\sigma,\tau \in \finBinarySequence$. Let $N_{\sigma,\tau} = \{(x,y) \in {}^2(\cantorspace) : x \in N_\sigma \wedge y \in N_\tau\}$. 

$\{N_{\sigma,\tau} : \sigma,\tau \in \finBinarySequence\}$ is a basis for the topology on ${}^2(\cantorspace)$. 

$N_{\sigma,\tau,\rho}$ is defined similarly for ${}^3(\cantorspace)$.
\end{definition}

\Begin{definition}{basic open neighborhood pcantorspace}
Let $n \in \omega$ and $\sigma : n \rightarrow \finBinarySequence$. Let $N_\sigma = \{x \in \pcantorspace : (\forall k < n)(\sigma(k) \subseteq x(k))\}$. 

$\{N_\sigma : \sigma \in {}^{<\omega}(\finBinarySequence)\}$ is a basis for the topology on $\pcantorspace$.

Let $\sigma, \tau : n \rightarrow \finBinarySequence$. Let $N_{\sigma,\tau} = \{(x,y) \in {}^2(\pcantorspace) : x \in N_\sigma \wedge y \in N_\tau\}$. 

$\{N_{\sigma,\tau} : \sigma,\tau \in {}^{<\omega}(\finBinarySequence) \wedge |\sigma| = |\tau|\}$ is a basis for the topology on ${}^2(\pcantorspace)$. 
\end{definition}

\Begin{definition}{bijection notation}
Let $A$ and $B$ be two sets. $A \approx B$ denotes that there is a bijection between $A$ and $B$.
\end{definition}

Often this paper will consider settings where the full axiom of choice may fail. In such contexts, not all sets have a cardinal, i.e. is in bijection with an ordinal. Similarity of size is more appropriately given by the existence of bijections. Recall the following method of producing bijections between sets which is provable in $\mathsf{ZF}$:

\Begin{fact}{cantor-schroder-bernstein}
(Cantor-Schr\"oder-Bernstein) $(\mathsf{ZF})$ Let $X$ and $Y$ be two sets. Suppose there are injections $\Phi: X \rightarrow Y$ and $\Psi : Y \rightarrow X$. Then there is a bijection $\Lambda : X \rightarrow Y$.
\end{fact}

\Begin{definition}{tuples notation}
Let $X$ and $Y$ be sets. ${}^X Y$ is the set of functions from $X$ to $Y$. 

$\mathscr{P}(X)$ is the power set of $X$.

Let $n \in \omega$. Define
$$\mathscr{P}^n(X) = \{F \in \mathscr{P}(X) : F \approx n\}$$
$$\mathscr{P}^{<\omega}(X) = \bigcup_{n \in \omega} \mathscr{P}^n(X)$$

Let $E$ be an equivalence relation on a set $X$. Let $n \in \omega$. Define
$$[X]^{n}_E = \{(x_0, ..., x_{n - 1}) \in {}^nX : (\forall i,j < n)(i \neq j \Rightarrow \neg(x_i \ E \ x_j))\}$$
$$[X]^{<\omega}_E = \bigcup_{n \in \omega} [X]_E^n$$
\end{definition}

\Begin{definition}{jonsson property}
Let $X$ be a set and $n \in \omega$. A set $X$ has the $n$-J\'onsson property if and only for all functions $f : [X]^n_= \rightarrow X$, there is some $Y \subseteq X$ so that $Y \approx X$ and $f[[Y]^n_=] \neq X$. $X$ has the J\'onsson property if and only if for all $f : [X]^{<\omega}_= \rightarrow X$, there is some $Y \subseteq X$ so that $Y \approx X$ and $f[[Y]^{<\omega}_=] \neq X$. 

A set $X$ has the classical $n$-J\'onsson property (or classical J\'onsson property) if and only the above holds with $[X]_n^=$ (or $[X]^{<\omega}_=$) replaced with $\mathscr{P}^n(X)$ (or $\mathscr{P}^{<\omega}(X)$, respectively).
\end{definition}

If $X$ is a wellordered set, one can identify a finite set $F \subseteq X$ with the increasing enumeration of its elements. Such a presentation is helpful for defining useful functions on $\mathscr{P}^n(X)$. In the absence of choice, it is easier to define functions when one considers order tuples from $[X]^n_=$. For this reason, the paper will be mostly concerned about the J\'onsson property as defined above rather than the classical J\'onsson property, although the classical version will be discussed in Section \ref{R/E0 does not have 3-Jonsson property}.

\Begin{definition}{jonsson functions}
Let $X$ be a set. $[X]^\omega_=$ and $\mathscr{P}^\omega(X)$ are defined as above (with $\omega$ in place of $n \in \omega$).

Let $N \in \omega \cup \{\omega\}$. A $N$-J\'onsson function for $X$ is a function $\Phi : [X]^N_= \rightarrow X$ so that for any $Y \subseteq X$ with $Y \approx X$, $\Phi[[Y]^N_=] = X$.  

A classical $N$-J\'onsson function for $X$ is defined in the same way as the above with $\mathscr{P}^N(X)$ instead of $[X]^N_=$. 
\end{definition}

With the axiom of choice, \cite{On-a-Problem-of-B-Jonsson} showed that every set has an $\omega$-J\'onsson map. The existence of $\omega$-J\'onsson maps for certain cardinals is where Kunen's original proof of the Kunen inconsistency used the axiom of choice. Note that for $N \in \omega \cup \{\omega\}$, a counterexample to the $N$-J\'onsson property for some set is equivalent to the existence of an $n$-J\'onsson function for that set.

\Begin{definition}{borel reducibility}
Let $X$ and $Y$ be Polish spaces. Let $E$ and $F$ be equivalence relations on $X$ and $Y$, respectively. A $\borel$ reduction between $X$ and $Y$ is a $\borel$ function $\Phi : X \rightarrow Y$ such that for all $a,b \in X$, $a \ E \ b$ if and only if $\Phi(a) \ F \ \Phi(b)$. 

This situtation is denoted by $E \leq_\borel F$. Define $E \equiv_\borel F$ if and only if $E \leq_\borel F$ and $F \leq_\borel E$. 
\end{definition}

\Begin{definition}{equivalence relation saturation definition}
Let $E$ be an equivalence relation on a set $X$. If $x \in X$, then $[x]_E = \{y \in X : y \ E \ x\}$ is the $E$-class of $x$. Let $A \subseteq X$. $[A]_E = \{y \in X : (\exists x \in A)(x \ E \ y)\}$ is the $E$-saturation of $A$.
\end{definition}

\Begin{definition}{mycielski property}
Let $X$ be a Polish space and $E$ be an equivalence relation on $X$. Let $n \in \omega$. $X$ has the $n$-Mycielski property if and only if for every $C \subseteq {}^n X$ which is comeager in ${}^nX$, there is a $\borel$ set $A \subseteq X$ so that $E \equiv_\borel E \upharpoonright A$ and $[A]^n_E \subseteq C$. 

$E$ has the Mycielski property if and only if for all sequences $(C_n : n \in \omega)$ such that for all $n \in \omega$, $C_n \subseteq {}^nX$ is comeager in ${}^nX$, there is a some set $A \subseteq X$ so that $E \equiv_\borel E \upharpoonright A$ and for all $n \in \omega$, $[A]^n_E \subseteq C_n$. 
\end{definition}

The Mycielski property of equivalence relations comes from the following eponymous result:

\Begin{fact}{Mycielski theorem}
(Mycielski) Let $(C_n : n \in \omega)$ be a sequence such that for each $n \in \omega$, $C_n \subseteq {}^n(\cantorspace)$ is a comeager subset of ${}^n(\cantorspace)$. Then there is a perfect set $P \subseteq \cantorspace$ so that for all $n \in \omega$, $[P]^n_= \subseteq C_n$. 
\end{fact}

\Begin{definition}{cotinuity property}
Let $E$ be an equivalence relation on a Polish space $X$. Let $n \in \omega$. $E$ has the $n$-continuity property if and only if for every function $f: {}^n X \rightarrow X$, there is some $\borel$ $A \subseteq X$ so that $E \equiv_\borel E \upharpoonright A$ and $f \upharpoonright [A]^n_E$ is continuous. 
\end{definition}

\Begin{fact}{mycielski implies continuity}
Let $E$ be an equivalence relation on a Polish $X$ which has the $n$-Mycielski property. Then for every function $f : {}^n X \rightarrow X$ with the property of Baire (i.e. $f^{-1}[U]$ has the Baire property for every open set $U$), there is some $\borel$ $A \subseteq X$ with $E \equiv_\borel E \upharpoonright A$ so that $f \upharpoonright [A]^n_E$ is continuous. Hence if every set has the Baire property, then $E$ has the $n$-continuity property.
\end{fact}

\begin{proof}
Let $f : {}^n X \rightarrow X$. Since $f$ is Baire measurable, there is some $C \subseteq {}^nX$ so that $f \upharpoonright C$ is continuous. By the $n$-Mycielski property, there is some $A \subseteq X$ with $E \equiv_\borel E \upharpoonright A$ so that $[A]^n_E \subseteq C$. $f \upharpoonright [A]^n_E$ is continuous. 
\end{proof}

In place of the axiom of choice, the paper will often use the axiom of determinacy. The following is a quick description of determinacy:

\Begin{definition}{determinacy}
Let $X$ be a set. Let $A \subseteq {}^\omega X$. The game $G_A$ is defined as follows: Player 1 plays $a_i \in X$, and player 2 plays $b_i \in X$ for each $i \in \omega$. At turn $2i$, player 1 plays $a_i$, and at turn $2i + 1$, player 2 plays $b_i$. Let $f \in {}^\omega X$ be defined by $f(2i) = a_i$ and $f(2i + 1) = b_i$. Player 1 wins this play of $G_A$ if and only if $f \in A$. Player 2 wins otherwise. 

A winning strategy for player 1 is a function $\tau : {}^{<\omega} X \rightarrow X$ so that for any $(b_i : i \in \omega)$ if $(a_i : i \in \omega)$ is defined recursive by $a_0 = \tau(\emptyset)$ and $a_{n + 1} = \tau(a_0 ...a_nb_n)$, then player 1 wins the resulting play of $G_A$. A winning strategy for player 2 is defined similarly. 

The Axiom of Determinacy for $X$, denoted $\mathsf{AD}_X$, is the statement that for all $A \subseteq {}^\omega X$, $G_A$ has a winning strategy for some player.

$\mathsf{AD}$ refers to $\mathsf{AD}_2$ or equivalently $\mathsf{AD}_\omega$. $\mathsf{AD}_\bbR$ will also be used. Note that $\mathsf{AD}_\reals$ often will refer to $\mathsf{AD}_\cantorspace$ or $\mathsf{AD}_\bairespace$. 
\end{definition}

$\mathsf{AD}$ implies classical regularity properties for sets of reals: Every set of reals has the Baire property and is Lebesgue measurable. Every uncountable set of reals has a perfect subset. Every function on the reals is continuous on a comeager set. 

Uniformization however is more subtle:

\Begin{definition}{uniformization definition}
Let $R \subseteq \cantorspace \times \cantorspace$. Let $R^x = \{y : (x,y) \in R\}$. Suppose for all $x \in \cantorspace$, $R^x \neq \emptyset$. $R$ can be uniformized if and only if there is a some function $f : \cantorspace \rightarrow \cantorspace$ so that for all $x \in \cantorspace$, $(x,f(x)) \in R$. Such a function $f$ is called a uniformization of $R$. 
\end{definition}

\Begin{fact}{ADR uniformization}
$(\mathsf{ZF + AD_\reals})$ Every relation can be uniformized.
\end{fact}

\begin{proof}
Suppose $R \subseteq \cantorspace \times \cantorspace$ with the property that for all $x$, $R^x \neq \emptyset$. Consider the two step game where player 1 plays $a \in \cantorspace$ and player 2 responds with $b \in \cantorspace$. Player 2 wins if and only if $(a,b) \in R$. Clearly player 1 can not have a winning strategy. Any winning strategy for player 2 yields a uniformization of $R$. 
\end{proof}

Woodin has shown that if there is a measurable cardinal with infinitely many Woodin cardinals below it, then $L(\bbR) \models \mathsf{AD}$. Solovay showed in \cite{The-Independence-of-DC-from-AD} Lemma 2.2 and Corollary 2.4 that the relation $R(x,y)$ if and only $y$ is not ordinal definable from $x$ can not be uniformized in $L(\reals)$. Hence $\mathsf{AD}_\bbR$ is stronger than than $\mathsf{AD}$. $\mathsf{AD}$ is not capable of proving full uniformization.

\Begin{definition}{lifting functions on quotients}
Let $E$ be an equivalence relation on $\cantorspace$ and $n \in \omega$. Let $f : (\cantorspace \slash E)^n \rightarrow \cantorspace \slash E$. A lift of $f$ is a function $F : {}^n(\cantorspace) \rightarrow \cantorspace$ with the property that for all $(x_0, ...,x_n) \in {}^n(\cantorspace)$, $[F(x_0, ..., x_{n - 1})]_E = f([x_0]_E, ..., [x_{n - 1}]_E)$.
\end{definition}

\Begin{fact}{uniformization implies lifting}
Let $E$ be an equivalence relation on $\cantorspace$ and $n \in \omega$. Let $f : {}^n(\cantorspace \slash E) \rightarrow \cantorspace \slash E$. Define $R_f(x_0, ..., x_{n - 1}, y) \Leftrightarrow y \in f([x_0]_E, ..., [x_{n - 1}]_E)$. If $F$ is a uniformization of $R_f$ (with respect to the last variable), then $F$ is a lift of $f$. 

Under $\mathsf{AD}_\reals$, every such function has a lift.
\end{fact}

Many natural models of $\mathsf{AD}$ such as $L(\bbR)$ are not models of $\mathsf{AD}_\bbR$. However, functions on quotients of equivalence relations with all classes countable can still be uniformized: $\mathsf{AD}^+$ is a strengthening of $\mathsf{AD}$ which holds in all known models of $\mathsf{AD}$ (in particular  $L(\reals)$). See \cite{Axiom-of-Determinacy-Forcing-Axioms} Definition 9.6 for the definition of $\mathsf{AD}^+$. It is open whether $\mathsf{AD}$ and $\mathsf{AD}^+$ are equivalent. Also $\mathsf{AD}_\reals + \mathsf{DC}$ implies $\mathsf{AD}^+$, and it is open whether this holds without $\mathsf{DC}$.

\Begin{fact}{woodin countable section uniformization}
(Countable Section Uniformization) (Woodin) $(\mathsf{AD}^+)$ Let $R \subseteq \cantorspace \times \cantorspace$ have the property that for all $x \in \cantorspace$, $R^x$ is countable. Then $R$ can be uniformized. 
\end{fact}

\begin{proof}
See \cite{Ramsey-Ultrafiler-and-Countable-to-One-Uniformation} Theorem 3.2 for a proof. 
\end{proof}

\Begin{fact}{countable equivalence relation lifting}
$(\mathsf{AD^+})$ Let $E$ be an equivalence relation on $\cantorspace$ with all classes countable. Let $n \in \omega$ and $f : {}^n(\cantorspace \slash E) \rightarrow \cantorspace \slash E$. Then $f$ has a lift. 
\end{fact}

For the results of this paper, all results that require lifts can be replace by lift on some comeager set. The benefit is that such lift follows from comeager uniformization which is provable in just $\mathsf{ZF + AD}$.

\Begin{fact}{comeager uniformization}
(Comeager Unformization) $(\mathsf{ZF + AD})$ Let $R \subseteq \cantorspace \times \cantorspace$ be a result such that $(\forall x)(\exists y)R(x,y)$, then there is a comeager $C \subseteq \cantorspace$ and some function $f : C \rightarrow \cantorspace$ so that $(\forall x \in C)R(x,f(x))$. 
\end{fact}

By shrinking to an appropriate comeager set, one can assume that the uniformizing function is also continuous.

Often to use the techniques of forcing over countable elementary structures, the axiom of determinacy will need to be augmented by dependent choice ($\mathsf{DC}$). Kechris \cite{The-Axiom-of-Determinacy-Implies-Dependent-Choice} proved that $\mathsf{AD}$ and $\mathsf{AD + DC}$ have the same consistency strength by showing if $L(\reals) \models \mathsf{AD}$, then $L(\reals) \models \mathsf{DC}$. However, Solovay \cite {The-Independence-of-DC-from-AD} showed that $\mathsf{AD}_\reals + \mathsf{DC}$ has strictly stronger consistency strength than $\mathsf{AD}_\reals$.

If one is ultimately interested in functions $F : {}^n(\cantorspace) \rightarrow \cantorspace$ which are lifts of some function $f: {}^n(\cantorspace \slash E) \rightarrow \cantorspace \slash E$ only in order to infer information about $f$, then the demand in the Mycielski property that one considers tuples coming from a single set $A \subseteq \cantorspace$ such that $E \equiv_\borel E \upharpoonright A$ seems restrictive. If one ultimately will collapse back to the quotient, two sets $A$ and $B$ with the same $E$-saturation should work equally well. This motivates the following concepts:

\Begin{definition}{E product}
Let $n \in \omega$ and $E$ be an equivalence relation on some Polish space $X$. Let $(A_i : i < n)$ be a sequence of subsets of $X$. Define 
$$\prod_{i < n}^E A_i = \{(x_0, ..., x_{n - 1}) : (\forall i)(x_i \in A_i) \wedge (\forall i \neq j)(\neg(x_i \ E \ x_j))\}.$$ 
This set will sometimes be denoted $A_0 \times_E ... \times_E A_{n - 1}$.
\end{definition}

\Begin{definition}{weak Mycielski property}
Let $E$ be an equivalence relation on a Polish space $X$. Let $n \in \omega$. $E$ has the $n$-weak-Mycielski property if and only if for any $C \subseteq {}^nX$ which is comeager in ${}^nX$, there are $\borel$ sets $(A_i : i < n)$ with the property that for each $i < n$, $E \equiv_\borel E \upharpoonright A_i$ and $\prod_{i < n}^E A_i \subseteq C$. 
\end{definition}

\section{$\cantorspace$ Has the J\'onsson Property}\label{reals have the jonsson property}
This section will give a forcing style proof of Holshouser and Jackson's result that $\cantorspace$ has the J\'onsson property under some determinacy assumptions. The J\'onsson property for $\cantorspace$ will follow from a flexible fusion argument for Sacks forcing and the fact that under determinacy assumptions, every function is definable (on some perfect set) with certain absoluteness properties between countable structures and the real universe. Continuous functions will satisfy this property, and so the Baire property and the Mycielski property for $=$ can be used to show every function has such a definition on some perfect set. This definability can also be achieved by absoluteness phenomena that occur under $\mathsf{AD}^+$. Later, it will be shown that the Mycielski property fails for all the other simple equivalence relations considered; the hope is that such a definability and absoluteness approach could establish J\'onsson type properties without the Mycielski property. In the following, the fusion argument is essential for the combinatorics of the forcing argument. It is unclear what the relation is between fusion (or properness), the Mycielski property, and the J\'onsson property.

\Begin{definition}{Sacks forcing}
A tree $p$ on $2$ is a subset of $\finBinarySequence$ so that if $s \in p$ and $t \subseteq s$, then $t \in p$. $p$ is a perfect tree if and only if for all $s \in p$, there is a $t \supseteq s$ so that $t\hat{\ }0, t\hat{\ }1 \in p$. 

Let $\bbS$ denote the collection of all perfect trees on $2$, $\leq_\bbS = \subseteq$, and $1_\bbS = \finBinarySequence$. $(\bbS, \leq_\bbS, 1_\bbS)$ is Sacks forcing, denoted by just $\bbS$.

Let $p \in \bbS$. $s \in p$ is a split node if and only if $s\hat{\ }0, s\hat{\ }1 \in p$. $s \in p$ is a split of $p$ if and only if $s \upharpoonright (|s| - 1)$ is a split node of $p$. For $n \in \omega$, $s$ is a $n$-split of $p$ if and only if $s$ is a $\subseteq$-minimal element of $p$ with exactly $n$-many proper initial segments which are split nodes of $p$. 

Let $\splitt^n(p)$ denote the set of $n$-splits of $p$. Note that $|\splitt^n(p)| = 2^n$ and $\splitt^0(p) = \{\emptyset\}$. 

If $p,q \in \bbS$, define $p \leq_\bbS^n q$ if and only if $p \leq_\bbS q$ and $\splitt^n(p) = \splitt^n(q)$. 

If $p \in \bbS$ and $s \in p$, then define $p_s = \{t \in p : t \subseteq s \vee s \subseteq t\}$. 

Let $p \in \bbS$. Let $\Lambda$ be defined as follows:

\noindent (i) $\Lambda(p,\emptyset) = \emptyset$. 

\noindent (ii) Suppose $\Lambda(p,s)$ has been defined for all $s \in {}^n2$. Fix an $s \in {}^n2$ and $i \in 2$. Let $t \supseteq \Lambda(p,s)$ be the minimal split node of $p$ extending $\Lambda(p,s)$. Let $\Lambda(p, s\hat{\ }i) = t\hat{\ }i$. 

Let $\Xi(p,s) = p_{\Lambda(p,s)}$. 

For $n \in \omega$, let $\bbS^n$ denote the $n$-fold product of $\bbS$. If $p \in \bbS$, then let $p^n \in \bbS^n$ be defined so that for all $i < n$, $p^n(n) = p$. 

Let $n \in \omega$ and $m < n$. There is an $\bbS^n$-name $x^{n,m}_\mathrm{gen}$ which names the $m^\text{th}$ Sacks-generic real coming from an $\bbS^n$-generic filter. 
\end{definition}

\Begin{fact}{fusion lemma}
A fusion sequence is a sequence $\langle p_n : n \in \omega\rangle$ in $\bbS$ so that for all $n \in \omega$, $p_{n + 1} \leq_\bbS^n p_n$. The fusion of this sequence is $p_\omega = \bigcap_{n \in \omega} p_n$.

$p_\omega$ is a condition in $\bbS$.
\end{fact}

\Begin{lemma}{absoluteness implies R-jonsson property}
Let $f: [\cantorspace]^{<\omega}_= \rightarrow \cantorspace$. Let $f_n = f \upharpoonright [\cantorspace]^n_=$. Suppose there is a countable model $M$ of some sufficiently large fragment of $\mathsf{ZF}$, $p \in \bbS \cap M$, and a $\bbS^n$-name $\tau_n \in M$ so that $p^n \forces \tau_n \in \cantorspace$ and whenever $G^n \subseteq \bbS^n$ is $\bbS^n$-generic over $M$ with $p^n \in G^n$, $\tau_n[G^n] = f_n(x_\mathrm{gen}^{n,0}[G^n], ..., x_\mathrm{gen}^{n,n - 1}[G^n])$. Then there exists a $q \in \bbS$ so that $f[[[q]]^{<\omega}_=] \neq \cantorspace$. 
\end{lemma}

\begin{proof}
For each $n \in \omega$, let $(D_m^n : m \in \omega)$ be a sequence of dense open subsets of $\bbS^n$ in $M$ so that for all $m$, $D_{m + 1}^n \subseteq D_m^n$ and if $D$ is a dense open subset of $\bbS^n$ in $M$, then there is some $m$ so that $D_m^n \subseteq D$. 

Let $z \in \cantorspace \setminus (\cantorspace)^M$. 

A fusion sequence $\langle p_n : n \in \omega\rangle$ with $p_0 = p$ will be constructed with the following properties:

For all $n > 0$, $m \leq n$, and $(\sigma_0, ..., \sigma_{m - 1}) \in {}^m({}^n2)$ so that $\sigma_i \neq \sigma_j$ if $i \neq j$:

\noindent (i) $(\Xi(p_n,\sigma_0), ..., \Xi(p_n, \sigma_{m - 1})) \in D_n^m$.

\noindent (ii) There are some $k \in \omega$ and $i \in 2$ so that $z(k) \neq i$ and $(\Xi(p_n,\sigma_0), ..., \Xi(p_n, \sigma_{m - 1})) \forces_{\bbS^m}^M \tau_m(\check k) = \check i$. 

Suppose this fusion sequence $\langle p_n : n \in \omega \rangle$ could be constructed. Let $q$ be its fusion. Fix $m > 0$. Suppose $(x_0, ..., x_{m - 1}) \in [[q]]_=^m$. Let $G^m_{(x_0, ..., x_{m - 1})} = \{(p_0, ..., p_{m - 1}) \in \bbS^m \cap M : (\forall i < m)(x_i \in [p_i])\}$. Note that $G^m_{(x_0, ..., x_{m - 1})}$ is a $\bbS^m$ generic filter over $M$: There is some $L$ so that for all $k \geq L$, there are $\sigma_i^k \in {}^k2$ with the property that for all $i < m$, $x_i \in \Xi(p_k, \sigma^k_i)$ and for all $i \neq j$, $\sigma_i^k \neq \sigma_j^k$. Then for all $k \geq L$, $(\Xi(p_k,\sigma_0^k), ..., \Xi(p_k, \sigma_{m - 1}^k)) \in G^m_{(x_0, ..., x_{m - 1})}$. (i) asserts that this element belongs to $D^m_k$. Hence $G^m_{(x_0, ..., x_{m - 1})}$ is $\bbS^m$-generic over $M$. 

By (ii), $\tau_m[G_{(x_0, ..., x_{m - 1})}^m] \neq z$. Also 
$$\tau_m[G_{(x_0, ..., x_{m - 1})}^m] = f_m(x_\mathrm{gen}^{m,0}[G_{(x_0, ..., x_{m - 1})}^m], ..., x_\mathrm{gen}^{m,m -1}[G_{(x_0, ..., x_{m - 1})}^m]) = f_m(x_0, ..., x_{m - 1}).$$
Hence $z \notin f_m[[[q]]^m_=]$. Thus $f[[[q]]^{<\omega}_=] \neq \cantorspace$. 

The construction of the fusion sequence remains: Let $p_0 = p$. 

Suppose $p_n$ has been constructed with the above properties. For some $J \in \omega$, let $(\bar{\sigma}_k : k < J)$ enumerate all tuples of strings $(\sigma_0, ..., \sigma_{m - 1})$ where $m \leq n + 1$, $\sigma_i \in {}^{n + 1}2$, and if $i \neq j$, $\sigma_i \neq \sigma_j$. 

Next, one construct a sequence $r_{-1}, ..., r_{J - 1}$ as follows: Let $r_{-1} = p_n$. Suppose $r_k$ for $k < J - 1$ has been constructed. Suppose $\bar\sigma_{k + 1} = (\sigma_0, ..., \sigma_{m - 1})$. 

(Case I) There is some $(u_0, ..., u_{m - 1}) \leq_{\bbS^m} (\Xi(r_k, \sigma_0), ..., \Xi(r_k, \sigma_{m - 1}))$ and $c \in (\cantorspace)^M$ so that 
$$(u_0, ..., u_{m - 1}) \forces_{\bbS^m}^M \tau_m = \check c.$$ 
Also as $D^m_{n + 1}$ is dense open in $\bbS^m$, one may choose $(u_0, ..., u_{m - 1})$ satisfying the above and $(u_0, ..., u_{m - 1}) \in D^m_{n + 1}$. Note that since $z \notin M$ and $c \in M$, there must be some $j \in \omega$ and $i \in 2$ so that $c(j) \neq i$ and $z(j) = i$. Now let $r_{k + 1} \in \bbS$ be so that for all $\sigma \in {}^{n + 1} 2$
$$\Xi(r_{k + 1}, \sigma) = \begin{cases}
u_i & \quad (\exists i)(0 \leq i \leq m - 1)(\sigma = \sigma_i) \\
\Xi(r_k, \sigma) & \quad \text{otherwise}
\end{cases}.$$

(Case II) For all $(u_0, ..., u_{m - 1}) \leq_{\bbS^m} (\Xi(r_k, \sigma_0), ..., \Xi(r_k, \sigma_{m - 1}))$, 
$$(u_0, ..., u_{m - 1}) \forces_{\bbS^m}^M \tau_m \notin M.$$ 
Hence there are $(u_0, ..., u_{m - 1}) \leq_{\bbS^m} (\Xi(r_k, \sigma_0), ..., \Xi(r_k, \sigma_{m - 1}))$, $(v_0, ..., v_{m - 1}) \leq_{\bbS^m} (\Xi(r_k, \sigma_0), ..., \Xi(r_k, \sigma_{m - 1}))$, and $j \in \omega$ so that
$$(u_0, ..., u_{m - 1}) \forces_{\bbS^m}^M \tau_m(\check j) = \check 0$$
$$(v_0, ..., v_{m - 1}) \forces_{\bbS^m}^m \tau_m(\check j) = \check 1$$
Without loss of generality, suppose that $z(j) = 1$. Moreover since $D_{n + 1}^m$ is dense open, one may assume that $(u_0, ..., u_{m - 1}) \in D^m_{n + 1}$. Define $r_{k + 1}$ in the same way as in Case I.

Finally, let $p_{n + 1} = r_{J - 1}$. $p_{n + 1} \leq_{\bbS}^n p_n$ and condition (i) and (ii) are satisfied. This completes the construction.
\end{proof}

\Begin{fact}{continuity satisfy sacks name condition}
($\mathsf{ZF + DC}$) Let $p \in \bbS$ and $n \in \omega$. If $f_n : [[p]]^n_= \rightarrow \cantorspace$ is continuous, then there is a countable elementary $M \prec V_\Xi$ (for $\Xi$ some sufficiently large cardinal) and a name $\tau_n$ so that $M$, $\tau_n$, and $p$ satisfy the conditions of Lemma \ref{absoluteness implies R-jonsson property}. 
\end{fact}

\begin{proof}
If $f_n$ is continuous, then $f_n$ has $\Sigma_1^1$ and $\Pi_1^1$ formulas with parameters from $\cantorspace$ defining it. Let $M \prec V_\Xi$ be a countable elementary substructure containing $p$ and all the parameters used to define $f_n$. (This requires $\mathsf{DC}$.) Using Mostowski's absoluteness, $f_n$ (as defined by this formula) continues to define a function in the forcing extension $M[G]$, where $G \subseteq \bbS^n$ is $\bbS^n$-generic over $M$. So there is some $\bbS^n$-name $\tau_n \in M$ so that $p^n \forces_{\bbS^n}^M \tau_n = f_n(x_\mathrm{gen}^{n,0}, ..., x_\mathrm{gen}^{n, n - 1})$. Suppose $G^n \subseteq \bbS^n$ is $\bbS^n$-generic over $M$ and contains $p^n$. Then $M[G^n] \models \tau_n[G^n] = f_n(x_\mathrm{gen}^{n,0}[G^n], ..., x_\mathrm{gen}^{n, n - 1}[G^n])$. Let $\pi : M[G^n] \rightarrow N$ be the Mostowski collapse of $M[G^n]$. Since reals are not moved by the Mostowki collapse map $\pi$, $\pi(f_n)$ is still defined by the same formula. So 
$$N \models \pi(\tau_n[G^n]) = \tau_n[G^n] = \pi(f_n(x_\mathrm{gen}^{n,0}[G^n], ..., x_\mathrm{gen}^{n, n - 1}[G^n])) = f_n(x_\mathrm{gen}^{n,0}[G^n], ..., x_\mathrm{gen}^{n, n - 1}[G^n]).$$
Then applying Mostowski absoluteness, $\tau_n[G^n] = f_n(x_\mathrm{gen}^{n,0}[G^n], ..., x_\mathrm{gen}^{n, n - 1}[G^n])$. 
\end{proof}

\Begin{theorem}{R has Jonsson property BP}
(Holshouser-Jackson) Assume $\mathsf{ZF + DC}$ and all sets of reals have the Baire property. Then $\cantorspace$ has the J\'onsson property.
\end{theorem}

\begin{proof}
Let $f : [\cantorspace]^{<\omega}_= \rightarrow \cantorspace$. Let $f_n : [\bbR]^{n}_= \rightarrow \bbR$ be defined by $f_n = f \upharpoonright [\cantorspace]^n_=$. Since all sets of reals have the Baire property, there are comeager subsets $C_n \subseteq {}^n(\cantorspace)$ so that $f_n \upharpoonright C_n$ is continuous. By the theorem of Mycielski (i.e. $=$ has the Mycielski property), there is a perfect tree $p$ so that $[[p]]_=^n \subseteq C_n$ for all $n \in \omega$. Hence for all $n \in \omega$, $f_n \upharpoonright [[p]]_=^n$ is a continuous function. Then Fact \ref{continuity satisfy sacks name condition} and Lemma \ref{absoluteness implies R-jonsson property} imply that $\cantorspace$ has the J\'onsson property.
\end{proof}

\Begin{remark}{R jonsson and DC}
As a consequence of phrasing this argument using forcing, one needed to introduce countable elementary substructures. $\mathsf{DC}$ is needed in general to obtain useful countable elementary substructures. A more direct topological argument can be used to avoid DC.
\end{remark}

\section{$\omega$-J\'onsson Function for $\cantorspace$}\label{w-jonsson function for reals}
Let $\mathsf{AC}_\omega^\bbR$ be the axiom of countable choice for $\cantorspace$: If $\mathcal{E}$ is a countable set of nonempty subsets of $\cantorspace$, then $\mathcal{E}$ has a choice function. 

Note that $\mathsf{ZF + AD}$ implies $\mathsf{AC_\omega^\reals}$. 

Using the axiom of choice, every set has an $\omega$-J\'onsson function. However, just $\mathsf{ZF + AC_\omega^\reals}$ implies there is a $\lborel$ classical $\omega$-J\'onsson function for $\cantorspace$. In fact, a slightly stronger statement holds:

\Begin{theorem}{omega Jonsson function for cantorspace}
$(\mathsf{ZF + AC_\omega^\reals})$ There is a $\lborel$ function $\Phi: \mathscr{P}^\omega(\cantorspace) \rightarrow \cantorspace$ so that if $B \subseteq \cantorspace$ is uncountable, then $\Phi[\mathscr{P}^\omega(B)] = \cantorspace$. 

There is a $\lborel$ classical $\omega$-J\'onsson function for $\cantorspace$. 
\end{theorem}

\begin{proof}
Let $A$ be a countable subset of $\cantorspace$.

Let $a_\emptyset^A$ be the longest element of $\finBinarySequence$ which is an initial segment of every element of $A$. 

If $a^A_\sigma$ is not defined, then $a^A_{\sigma \hat{\ }i}$ is not defined for $i \in 2$. If  $a^A_\sigma$ is defined, let $a^A_{\sigma\hat{\ }i}$ be the longest element of $\finBinarySequence$ which is an initial segment of every element of $A \cap N_{a^A_{\sigma}\hat{\ }i}$, if it exists. Otherwise $a^A_{\sigma\hat{\ }i}$ is undefined. (Note this happens if and only if $A \cap N_{a_\sigma\hat{\ } i}$ is a singleton.)

For $A \in \mathscr{P}^\omega(\cantorspace)$, let $\Psi(A)$ be the collection of $\sigma \in \finBinarySequence$ so that $a_\sigma^A$ is defined. $\Psi(A)$ is an infinite tree on $2$ with possibly dead nodes and is not perfect. $\Psi$ is a $\lborel$ function.

Using some recursive coding, let $X$ be the collection of reals coding infinite binary trees (which may have dead branches). $X$ is an uncountable $\Pi_1^0$ set. 

Let $T \in X$. Let $\hat T = \{\sigma \hat{\ }\tilde 0, \sigma \hat{\ } 1 \hat{\ }\tilde 0 : \sigma \in T\}$.  $\hat T \in \mathscr{P}^\omega(\cantorspace)$. One seeks to show that $\Psi(\hat{T}) = T$.  To see this, the following claim is helpful: If $\sigma \in T$, then $a^{\hat T}_\sigma = \sigma$ and if $\sigma \notin T$, then $a^{\hat T}_\sigma$ is undefined.

This claim is proved by induction: $\emptyset \in T$ so $\tilde 0, 1\hat{\ }\tilde 0 \in \hat T$. $a_\emptyset^{\hat T} = \emptyset$. Suppose this holds for $\sigma$. Suppose $\sigma\hat{\ }i \in T$. Then $\sigma\hat{\ }i\hat{\ }0\hat{\ }\tilde 0$ and $\sigma\hat{\ }i\hat{\ }1\hat{\ }\tilde 0$ are both in $\hat T$. By induction, $a_{\sigma}^{\hat T} = \sigma$. The longest string which is an initial segment of every element of $\hat T \cap N_{a_{\sigma}^{\hat T}\hat{\ }i} = \hat T \cap N_{\sigma\hat{\ }i}$ is $\sigma\hat{\ }i$. This shows $a_{\sigma\hat{\ }i}^{\hat T} = \sigma\hat{\ }i$. Suppose $\sigma\hat{\ }i \notin T$. Either $a_\sigma^{\hat T}$ is undefined or $a_\sigma^{\hat T}$ is defined. If $a_\sigma^{\hat T}$ is undefined, then $a_{\sigma\hat{\ }i}^{\hat T}$ is undefined. Suppose $a_{\sigma}^{\hat T}$ is defined. By induction, $a_\sigma^{\hat T} = \sigma$. $\sigma \in T$ implies that $\sigma\hat{\ }0\hat{\ }\tilde 0$ and $\sigma\hat{\ }1\hat{\ }\tilde 0$ are both in $\hat T$. Since $\sigma\hat{\ }i \notin T$, $N_{\sigma\hat{\ }i} \cap \hat T = N_{a_\sigma^{\hat T}\hat{\ }i} \cap \hat T = \{\sigma\hat{\ }i\hat{\ }\tilde 0\}$. $a_{\sigma\hat{\ }i}^{\hat T}$ is undefined. This completes the proof of the claim.  

This shows $\Phi(\hat T) = T$. Hence $\Psi[\mathscr{P}^\omega(\cantorspace)] = X$. Let $\Gamma : X \rightarrow \cantorspace$ be a $\lborel$ bijection. Let $\Phi = \Gamma \circ \Psi$. 

Let $B$ be an uncountable subset of $\cantorspace$. Then there is an uncountable $C \subseteq B$ which has no isolated points. One way to see this is to note that using a countable basis, the Cantor-Bendixson process must stop at a countable ordinal. The fixed point starting from $B$ would be an uncountable set with no isolated points. More directly: a condensation point of $B$ is a point so that every open set containing that point contains uncountable many points of $B$. Let $C$ be the set of condensation points of $B$ which are in $B$. 

Fix such a set $C$. Let $\mathcal{E} = \{N_\sigma \cap C : \sigma \in \finBinarySequence \wedge N_\sigma \cap C \neq \emptyset\}$. $\mathcal{E}$ is a countable set. Using $\mathsf{AC}_\omega^\bbR$, let $\Lambda$ be a choice function for $\mathcal{E}$. Let $T$ be any infinite binary tree on $2$. 

The following objects will be constructed: 

\noindent (I) $c_s \in C$ for each $s \in \finBinarySequence$.

\noindent (II) A strictly increasing sequence $(k_i : i \in \{-1\} \cup \omega)$ of integers.

For each $n \in \omega$, let $A_n = \{c_s : s \in {}^n 2\}$. The objects above will satisfy the following properties:

\noindent (i) If $s \in T$, then $a_s^{A_{|s| + 1}}$ is defined and has length less than $k_{|s|}$. If $s \notin T$, then $A_{|s| + 1} \cap N_{c_s \upharpoonright k_{|s|}} = \{c_s\}$.  

\noindent (ii) If $s \in T$, then $c_{s\hat{\ }i} \supseteq a_s^{A_{|s| + 1}} \hat{\ } i$ for each $i \in 2$. 

\noindent (iii) For all $m$, if $n > m$, then $\{x \upharpoonright k_m : x \in A_m\} = \{x \upharpoonright k_m : x \in A_n\}$.  

Let $c_\emptyset$ be any element of $C$. Let $k_{-1} = 0$. Let $A_0 = \{c_\emptyset\}$. 

Suppose for $m \in \omega$, $c_s \in C$ for all $s \in {}^m2$ and $k_{m - 1}$ have been defined. Suppose properties (i) to (iii) hold for $t \in \finBinarySequence$ with $|t| < m$. Let $s \in T \cap {}^m2$. Since $c_s \in C$ and $C$ has no isolated points, there is some $m_s > k_{m - 1}$ so that $N_{(c_s \upharpoonright m_s) \hat{\ } (1 - c_s(m_s))} \cap C \neq \emptyset$. Let $c_{s\hat{\ }c_s(m_s)} = c_s$ and let $c_{s\hat{\ }(1 - c_s(m_s))} = \Lambda(N_{(c_s \upharpoonright m_s) \hat{\ }(1 - c_s(m_s))})$. If $s \notin T$, then let $c_{s\hat{\ } i} = c_s$ for each $i \in 2$. Let $k_m = \sup\{m_s + 1 : s \in {}^m 2 \cap T\}$. 

Since for each $s \in T \cap {}^m2$, $m_s > k_{m - 1}$, (i) to (iii) still hold for $t$ with $|t| < m$. Let $s \in T$. Using the induction hypothsis for (ii) on $s \upharpoonright m - 1$, one has that $c_s \supseteq a_{s \upharpoonright m - 1}^{A_m} \hat{\ }s(m - 1)$. $c_{s\hat{\ }0}$ and $c_{s\hat{\ } 1}$ extend $a_{s \upharpoonright m - 1}^{A_m}\hat{\ } s(m - 1)$. This shows that $a_{s}^{A_{m + 1}}$ is defined. In fact, $a_s^{A_{m + 1}} = c_s \upharpoonright m_s$. If $s \notin T$, (i) is clear from the construction. Properties (i) to (iii) hold for $s \in {}^{m}2$. 

Let $A = \bigcup_{n \in \omega} A_n$. Note that $A$ is countably infinite and $A \subseteq C \subseteq B$. From the above properties, if $s \in T$, then $a_s^A$ is defined and in fact equal to $a_s^{A_{|s| + 1}}$. Suppose $s \notin T$. Let $t \subseteq s$ be maximal with $t \in T$. The above properties imply that $A \cap N_{a_t^A \hat{\ }s(|t|)} = \{c_{t\hat{\ }s(|t|)}\} = \{c_s\}$. Hence $a_{t\hat{\ } s(|t|)}^A$ is not defined and hence $a_s^A$ is not defined. Thus $T = \Psi(A)$. This shows that $\Phi[\mathscr{P}^\omega(B)] = \cantorspace$. $\Phi$ is an $\omega$-J\'onsson function for $\cantorspace$.
\end{proof}

\Begin{question}{ACR and omega jonsson map}
Under  $\mathsf{ZF} + \neg \mathsf{AC}^\bbR_\omega$, can there be a classical $\omega$-J\'onsson function for $\cantorspace$?

The first statement of Theorem \ref{omega Jonsson function for cantorspace} may not be true without $\neg\mathrm{AC}^\bbR_\omega$: Let $\bbC_\omega$ denote the finite support product of Cohen forcing $\bbC$. Let $G \subseteq \bbC_\omega$ be $\bbC_\omega$-generic over $L$. For each $n \in \omega$, let $c_n$ be the $n^\text{th}$-Cohen generic real naturally added by $G$. Let $A = \{c_n : n \in \omega\}$. Let $H = (\mathrm{HOD}(A \cup \{A\}))^{L[G]}$. $H$ is called the Cohen-Halpern-L\'evy model. In $H$, $A$ has no countably infinite subsets. Hence the first statement of Theorem \ref{omega Jonsson function for cantorspace} cannot hold. However $A$ is not in bijection with $\cantorspace$. This suggest the following natural question:

In $H$, is there a classical $\omega$-J\'onsson function for $\cantorspace$?
\end{question}

\section{The Structure of $E_0$}\label{the structure of E0}
\Begin{definition}{E0 equivalence relation}
$E_0$ is the equivalence relation defined on $\cantorspace$ by $x \ E_0 \ y$ if and only if $(\exists n)(\forall k > n)(x(k) = y(k))$. 
\end{definition}

\Begin{definition}{E0 tree}
Let $\{s, v_n^i : i \in 2 \wedge n \in \omega\} \subseteq \finBinarySequence$ have the property that for all $n \in \omega$ and $i \in 2$, $\langle i \rangle \subseteq v_n^i$ and $|v_n^0| = |v_n^1|$.  

Let $\varphi(\emptyset) = s$. If $\sigma \in \finBinarySequence$ and $|\sigma| > 0$, then let $\varphi(\sigma) = s \hat{\ }v_0^{\sigma(0)}\hat{\ } ...\hat{\ }v_{|\sigma| - 1}^{\sigma(|\sigma| - 1)}$. 

A perfect tree $p$ is an $E_0$-tree if and only if there is a sequence $\{s, v_n^i : i \in 2 \wedge n \in \omega\}$ with the above properties so that $p$ is the $\subseteq$-downward closure of $\{\varphi(\sigma) : \sigma \in \finBinarySequence\}$. 

Let $\PE$ be the collection of all perfect $E_0$ trees. If $p,q \in \PE$, then $p \leq_\PE q$ if and only if $p \subseteq q$. Let $1_{\PE} = \finBinarySequence$. $(\PE,\leq_\PE,1_\PE)$ is forcing with perfect $E_0$-trees.

If $p \in \PE$, then the notation $s^p$ and $v_n^{i,p}$ will be used to denote the strings witnessing $p$ is a perfect $E_0$-tree.

Let $\Phi : \cantorspace \rightarrow [p]$ be defined by $\Phi(x) = \bigcup_{n \in \omega} \varphi(x \upharpoonright n)$, where $\varphi$ is associated with the $E_0$-tree $p$ as above. $\Phi$ is the canonical homeomorphism of $\cantorspace$ onto $[p]$, and $\Phi$ is a reduction witnessing $E_0 \leq_\borel E_0 \upharpoonright [p]$. 
\end{definition}

\Begin{fact}{E0 characterization}
Suppose $B$ is a $\analytic$ set so that $E_0 \leq_\borel E_0 \upharpoonright B$. Then there is an $E_0$-tree $p$ so that $[p] \subseteq B$. 
\end{fact}

\begin{proof}
This is implicit in \cite{Glimm-Effros-Dichotomy-for-Borel-Equivalence-Relations}. See \cite{Descriptive-Set-Theory-and-Definable-Forcing} Lemma 2.3.29 and \cite{Borel-Equivalence-Relations} Theorem 10.8.3.
\end{proof}

The weak Mycielski property for $E_0$ considers $E_0$-products of $\borel$ sets $A_0, ..., A_{n - 1}$ so that $E_0 \equiv_\borel E_0 \upharpoonright A_i$ and $[A_i]_{E_0} = [A_j]_{E_0}$. Showing the failure of the weak Mycielski property requires finding some structure shared by all of the sets $A_0, ..., A_{n - 1}$. For instance, are there perfect $E_0$-trees $p_i$ so that $[p_i] \subseteq A_i$ and $[[p_i]]_{E_0} = [[p_j]]_{E_0}$? How similar can $p_i$ and $p_j$ be chosen to be? 

A simpler solution using the $\sigma$-additivity of the $E_0$-ideal, which follows Fact \ref{E0 characterization}, will be given first. A stronger result giving more information using effective methods will follow.

\Begin{fact}{E0 sigma additivity}
For each $n \in \omega$, suppose $A_n \subseteq \cantorspace$ is $\analytic$ and there is no $E_0$-tree $p$ so that $[p] \subseteq A_n$. Then there is no $E_0$ tree $p$ so that $[p] \subseteq \bigcup_{n \in \omega} A_n$. 
\end{fact}

\begin{proof}
Suppose there is some $E_0$-tree $p$ so that $[p] \subseteq \bigcup_{n \in \omega} A_n$. Let $\Phi : \cantorspace \rightarrow [p]$ be the canonical injective reduction witnessing $E_0 \leq_\borel E_0 \upharpoonright [p]$. For each $n \in \omega$, $\Phi^{-1}[A_n]$ is a $\analytic$ set, and $\cantorspace = \bigcup_{n \in \omega} \Phi^{-1}[A_n]$. There is some $m \in \omega$ so that $\Phi^{-1}[A_m]$ is nonmeager. Therefore, there is some continuous injective function $\Psi : \cantorspace \rightarrow \Phi^{-1}[A_m]$ which witnesses $E_0 \leq_\borel E_0 \upharpoonright \Phi^{-1}[A_m]$. $\Phi \circ \Psi$ witnesses $E_0 \leq_{\borel} E_0 \upharpoonright A_m$. This implies there is some $E_0$-tree $q$ so that $[q] \subseteq A_m$. Contradiction. 
\end{proof}

\Begin{definition}{E0 tail notations}
If $x \in \cantorspace$ and $n \in \omega$, let $x_{\geq n} \in \cantorspace$ be defined by $x_{\geq n} (k) = x(n + k)$. 

If $A \subseteq \cantorspace$, then let $(A)_{\geq n} = \{z : (\exists x \in A)(z = x_{\geq n})\}$. 
\end{definition}

\Begin{definition}{switch functions}
Let $s \in \finBinarySequence$. Define $\switch_s : \cantorspace \rightarrow \cantorspace$ by
$$\switch_s(x)(n) = \begin{cases}
s(n) & \quad n < |s| \\
x(n) & \quad \text{otherwise}
\end{cases}.$$

Also if $\sigma \in \finBinarySequence$, $\switch_s(\sigma) \in {}^{|\sigma|}2$ is defined as above just for $n < |\sigma|$. 
\end{definition}

\Begin{theorem}{E0 tree inside set same saturation}
Let $n \in \omega$. For $k < n$, let $A_k \subseteq \cantorspace$ be $\analytic$ so that $E_0 \leq_\borel E_0 \upharpoonright A_k$ and for all $k < n - 1$, $[A_k]_{E_0} \subseteq [A_{k + 1}]_{E_0}$. Then there exists $E_0$-trees $p_k$ so that $[p_k] \subseteq A_k$ and for all $a,b < n$, $|s^{p_a}| = |s^{p_b}|$, and $v_m^{i,p_a} = v_m^{i,p_b}$ for all $m \in \omega$ and $i \in 2$. 
\end{theorem}

\begin{proof}
For each $c \in \omega$, let 
$$E_c = \left\{(x_0, ..., x_{n - 1}) \in \prod_{k < n} A_k : (\forall i, j < n)((x_i)_{\geq c}  = (x_j)_{\geq c})\right\}.$$ 
For each $c \in \omega$, $E_c$ is a $\analytic$ set. For $k < n$, let $\pi_k: {}^n(\cantorspace) \rightarrow \cantorspace$ be the projection map onto the $k^\text{th}$ coordinate. $\pi_0[E_c]$ is $\analytic$ for each $c \in \omega$. Since $[A_i]_{E_{0}} \subseteq [A_{i + 1}]_{E_0}$, $\bigcup_{c \in \omega} \pi_0[E_c] = A_0$. By Fact \ref{E0 sigma additivity}, there is some $m \in \omega$ so that $\pi_0[E_m]$ contains the body of an $E_0$-tree $q_0$. By choosing an appropriate subtree, one may assume that $|s^{q_0}| > m$. Let $s_0 = s^{q_0} \upharpoonright m$. 

Fix $k < n - 1$. Suppose the $E_0$-tree $q_k$ and $s_{k} \in {}^m2$ have been constructed so that $s_k \subseteq s^{q_k}$ and $[q_k] \subseteq \pi_k[E_m]$. Then 
$$\left[\bigcup_{s \in {}^{m}2} \switch_s[[q_k]] \cap \pi_{k + 1}[E_m]\right]_{E_0} = [[q_k]]_{E_0}.$$
By Fact \ref{E0 sigma additivity}, there is some $s_{k + 1} \in {}^m2$ so that $\switch_{s_{k + 1}}[[q_{k}]] \cap \pi_{k + 1}[E_m]$ contains an $E_0$-tree. Let $q_{k + 1}$ be such an $E_0$-tree. Note that $\switch_{s_k}[q_{k + 1}]$ is an $E_0$-subtree of $q_k$.

For $i < n$, let $p_i = \switch_{s_i}[q_{n - 1}]$. Note that $[p_i] \subseteq \pi_i[E_m] \subseteq A_i$.
\end{proof}

The rest of this section will prove a result that implies Theorem \ref{E0 tree inside set same saturation} using an effective definability condition. The methods from \cite{Borel-Equivalence-Relations} Theorem 10.8.3 will be used to simultaneously produce $E_0$-trees, which are very similar to each other, through several sets.

\Begin{definition}{gandy-harrington forcing}
Let $z \in \cantorspace$. Let $\bbP_z$ be the forcing of nonempty $\Sigma_1^1(z)$ sets ordered by inclusion with largest element $\bbP_z = \cantorspace$. $\bbP_z$ is $z$-Gandy-Harrington forcing.
\end{definition}

\Begin{fact}{gandy-harrington dense open nonempty}
There is a $z$-recursive (in a suitable sense) collection $\mathcal{D} = \{D_n : n \in \omega\}$ of dense open subsets of $\bbP_z$ so that if $G \subseteq \bbP_z$ is generic for $\mathcal{D}$, then $\bigcap G \neq \emptyset$. 
\end{fact}

\begin{proof}
See \cite{Borel-Equivalence-Relations} Theorem 2.10.4.
\end{proof}

\Begin{fact}{coding of borel sets}
Let $z \in \cantorspace$. There is a $\Pi_1^1(z)$ set $D \subseteq \omega$, $\Sigma_1^1(z)$ set $P \subseteq \omega \times \cantorspace$, and $\Pi_1^1(z)$ set $Q \subseteq \omega \times \cantorspace$ with the following properties: 

(i) For all $e \in D$, $P^e = Q^e$, where if $X \subseteq \omega \times \cantorspace$, then $X^e = \{x \in \cantorspace : (e,x) \in X\}$. 

(ii) If $X \subseteq \cantorspace$ is $\lborel(z)$, then there is some $e \in D$ so that $X = P^e = Q^e$. 
\end{fact}

\Begin{definition}{E0 S and H sets}
Let $z \in \cantorspace$. Let $S_z$ be the union of all $\lborel(z)$ sets $C$ so that for all $x,y \in C$, $\neg(x \ E_0 \ y)$.

Let $H_z = \cantorspace \setminus S_z$. 
\end{definition}

\Begin{fact}{property of E0 S and H sets}
Let $z \in \cantorspace$. $S_z$ is $\lcoanalytic(z)$. $H_z$ is $\lanalytic(z)$. If $X \cap H_z \neq \emptyset$ and $X$ is $\lanalytic(z)$, then there exists $x,y \in X$ with $x \neq y$ and $x \ E_0 \ y$. $H_z$ is $E_0$-saturated.
\end{fact}

\begin{proof}
Let $D$, $P$, and $Q$ be the sets from Fact \ref{coding of borel sets}. Note that
$$x \in S_z \Leftrightarrow (\exists e)(e \in D \wedge x \in Q^e \wedge (\forall f,g)((f \neq g \wedge f,g \in P^e) \Rightarrow \neg(f \ E_0 \ g))).$$
$S_z$ is $\lcoanalytic(z)$. Hence $H_z$ is $\lanalytic(z)$. 

Let $\mathcal{A}$ be the collection of all $\lanalytic(z)$ subsets of $\cantorspace$ whose elements are pairwise $E_0$-inequivalent. Let $U \subseteq \omega \times \cantorspace$ be a universal $\lanalytic(z)$ set. 
$$\{e : U^e \in \mathcal{A}\} = \{e : (\forall f,g)((f,g \in U^e \wedge f \neq g) \Rightarrow \neg(f \ E_0 \ g))\}$$
The above is a $\lcoanalytic(z)$ set. So $\mathcal{A}$ is a collection of $\lanalytic(z)$ sets which is $\lcoanalytic(z)$ in the codes. By $\lanalytic(z)$-reflection (see \cite{Borel-Equivalence-Relations}, Theorem 2.7.1), every $\lanalytic(z)$ set $X$ whose elements are $E_0$-inequivalent has a $\lborel(z)$ set $C$ whose elements are $E_0$-inequivalent and $X \subseteq C$. 

Suppose $X$ is $\lanalytic(z)$, $X \cap H_z \neq \emptyset$, and the elements of $X$ are pairwise $E_0$-inequivalent. By the previous paragraph, there is some $\lborel(z)$ set $C$ which also $E_0$-inequivalent and $X \subseteq C$. Then $X \subseteq S_z$. Contradiction.

Suppose $x \in H_z$, $y \notin H_z$, and $x \ E_0 \ y$. Let $n \in \omega$ be so that $x_{\geq n} = y_{\geq n}$. $y \in S_z$ implies that there is some $\lborel(z)$ $E_0$-inequivalent set $X$ so that $y \in X$. $\mathsf{switch}_{x\upharpoonright n}[X]$ is a $\lborel(z)$ $E_0$-inequivalent set containing $x$. This contradicts $x \in H_z$. This shows $H_z$ is $E_0$-saturated.
\end{proof}

\Begin{lemma}{same tail meet dense same tail}
Let $z \in \cantorspace$ and $n,\ell \in \omega$. Suppose $(\hat{B}_i : i < n)$ is a collection of nonempty $\lanalytic(z)$ sets. Suppose for all $i,j < n$, $(\hat{B}_i)_{\geq \ell} = (\hat{B}_j)_{\geq \ell}$. Let $D \subseteq \bbP_z$ be a dense open subset of the forcing $\bbP_z$. Then there is a collection $(B_i : i < n)$ of nonempty $\Sigma_1^1(z)$ sets so that for all $i,j < n$, $B_i \in D$, $B_i \subseteq \hat{B}_i$, and $(B_i)_{\geq \ell} = (B_j)_{\geq \ell}$. 
\end{lemma}

\begin{proof}
For $k < n$, $\lanalytic(z)$ sets $\{B_i^k : -1 \leq k < n \wedge 0 \leq i < n\}$ will be constructed with the properties that 

\noindent (i) For all $i < n$, $B_i^i \in D$. 

\noindent (ii) If $-1 \leq k < n - 1$ and $0 \leq i < n$, then $B_i^{k + 1} \subseteq B_i^{k}$. 

\noindent (iii) For all $-1 \leq k < n$ and $0 \leq i,j \leq n$, $(B_i^k)_{\geq \ell} = (B_j^k)_{\geq \ell}$.

Note that this implies that if $k \geq i$, then $B_i^k \in D$. 

Let $B^{-1}_i = \hat{B}_i$. (iii) is satisfied. 

Suppose for $-1 \leq k < n - 1$ and $0 \leq i < n$, $B_i^k$ has been constructed with the desired properties. Since $D$ is dense open, there is some nonempty $\Sigma_1^1(z)$ set, denoted $B_{k + 1}^{k + 1}$, so that $B_{k + 1}^{k + 1} \subseteq B_{k + 1}^k$ and $B_{k + 1}^{k + 1} \in D$. For $0 \leq i < n$, let $B_{i}^{k + 1} = \{x \in B_i^k : (\exists z)(z \in B_{k + 1}^{k + 1} \wedge x_{\geq \ell} = z_{\geq \ell})\}$. All the conditions are satisfied.

Finally, let $B_i = B_i^{n - 1}$.

\end{proof}

\Begin{theorem}{effective structure of E0 sets}
Let $z \in \cantorspace$ and $n \in \omega$. Let $(A_a : a < n)$ be a collection of $\lanalytic(z)$ sets so that $\bigcap_{a < n} [A_a \cap H_z]_{E_0} \neq \emptyset$. Then there are $E_0$-trees $(p_a : a \in n)$ so that for all $a,b < n$, $k \in \omega$ and $i < 2$, 

(i) $|s^{p_a}| = |s^{p_b}|$ and $v_k^{i,p_a} = v_k^{i,p_b}$

(ii) $[p_a] \subseteq A_a$. 
\end{theorem}

\begin{proof}
The following objects will be constructed: For each $a < n$, $k \in \omega$, $i \in 2$, and $t \in \finBinarySequence$, 

\noindent (a) $w_t^a, s^a, v_k^i \in \finBinarySequence$ 

\noindent (b) $\ell_k \in \omega$ and for all $k \in \omega$, $\ell_{k} \leq m_k < \ell_{k +1}$

\noindent (c) $\lanalytic(z)$ nonempty sets $X_t^a$ 

\noindent with the following properties

\noindent (i) For all $a < n$ and $t \in \finBinarySequence$, $|w^a_t| = \ell_{|t|}$. For all $a < n$ and $k \in \omega$, $|s^a| = |s^b|$ and $|v_k^0| = |v_k^1|$. For all $k \in \omega$ and $i \in 2$, $\langle i \rangle \subseteq v_i^k$.  For all $a < n$ and $t \in \finBinarySequence$, if $|t| = 1$, then $s^a \subseteq w_t^a$ and if $t > 1$, then  $s^a\hat{\ }v_0^{t(0)}\hat{\ } ...\hat{\ } v_{|t| - 2}^{t(|t| - 2)} \subseteq w_t^a$. 

\noindent (ii) If $a \in n$, $t \in \finBinarySequence$, $u \in \finBinarySequence$, and $t \subseteq u$, then $X_u^a \subseteq X_t^a$, $X_t^a \subseteq N_{w_t^a}$, and $X_\emptyset^a \subseteq A_a \cap H_z$. 

\noindent (iii) Let $\mathcal{D} = (D_n : n \in \omega)$ be the collection of dense open subset of $\bbP_z$ from Fact \ref{gandy-harrington dense open nonempty}. For all $t \in \finBinarySequence$, $X_t^a \in D_{|t|}$. 

\noindent (iv) For all $k < \omega$, $\ell_k < \ell_{k + 1}$. For all $k < \omega$, $t,u \in {}^k2$, and $a,b \in n$, $(X_t^a)_{\geq \ell_k} = (X_u^b)_{\geq \ell_k}$.

Suppose objects with these properties can be constructed. For $a < n$, let $p^a$ be the $E_0$-tree given by $s^{p_a} = s^a$ and $v^{i,p_a}_k = v_k^i$. Let $\Phi^a : \cantorspace \rightarrow [p_a]$ be the canonical map associated with the $E_0$-tree $p_a$. For each $x \in \cantorspace$, let $G_x^a$ be the $\bbP_z$-filter generated by the upward closure of $\{X^a_{x \upharpoonright k} : k \in \omega\}$. $G_x^a \cap D_k \neq \emptyset$ since $X_{a \upharpoonright k}^a \in D_k$. $G_x^a$ is a filter generic for $\mathcal{D}$. By Fact \ref{gandy-harrington dense open nonempty}, $\bigcap G_x^a \neq \emptyset$. By (i) and (ii), $\bigcap G_x^a = \{\Phi^a(x)\}$. Thus $\Phi^a(x) \in X_\emptyset^a \subseteq A_a \cap H_z$. Therefore, $[p_a] \subseteq A_a$. 

Next, the construction will be described. Since $\bigcap_{a < n} [A_a \cap H_z]_{E_0} \neq \emptyset$, let $(x_a : a < n)$ be elements of $\cantorspace$ so that for all $a,b < n$, $x_a \in A_a \cap H_z$ and $x_a \ E_0 \ x_b$. Choose $\ell_0 \in \omega$ so that for all $a,b < n$, $(x_a)_{\geq \ell_0} = (x_b)_{\geq \ell_0}$. Let $w_\emptyset^a = x_a \upharpoonright \ell_0$. Let $Z = \{r \in \cantorspace : (\exists y_0, ..., y_{n - 1})(\bigwedge_{a < n} y_a \in A_a \cap H_z \wedge w_\emptyset^a\hat{\ }r = y_a)\}$. $Z$ is a nonempty $\lanalytic(z)$ set. Let $B_a = \{x \in A_a \cap H_z : (\exists r)(r \in Z \wedge x \supseteq w_\emptyset^a \wedge x_{\geq \ell_0} = r)\}$. For each $a < n$, $B_a$ is a nonempty $\lanalytic(z)$ set. Note that for all $a,b < n$, $(B_a)_{\geq \ell_0} = (B_b)_{\geq \ell_0}$. Applying Lemma \ref{same tail meet dense same tail}, find sets $(X_\emptyset^a : a \in \omega)$ so that $X_\emptyset^a \subseteq B_a$ and $X_{\emptyset}^a \in D_0$. 

Suppose $X_t^a$ has been constructed for $a < n$ and $t \in {}^k2$. $X_{0^k}^0 \subseteq H_z$ is nonempty and $\lanalytic(z)$. By Fact \ref{property of E0 S and H sets}, there are $x,y \in X_{0^k}^0$ so that $x \neq y$ and $x \ E_0 \ y$. By (i) and (ii), $x \upharpoonright \ell_k = y \upharpoonright \ell_k$. Therefore, there is some $\ell_{k + 1} > \ell_k$ so that $x_{\geq \ell_{k + 1}} = y_{\geq \ell_{k + 1}}$. Let $m_k$ be smallest $m$ so that $\ell_k \leq m < \ell_{k + 1}$ and $x(m) \neq y(m)$. Without loss of generality, suppose that $x(m_k) = 0$ and $y(m_k) = 1$. For $i \in 2$, let $w_{0^k\hat{\ }0}^0 = x \upharpoonright \ell_{k + 1}$ and $w^0_{0^k\hat{\ }1} = y \upharpoonright \ell_{k + 1}$. For $a \in n$ and $t \in \finBinarySequence$, let $w_{t\hat{\ }i}^a = \mathrm{switch}_{w_{t}^a} w_{0^k\hat{\ }i}^0$. 

If $k = 0$, then let $s^a = w_{\langle 0 \rangle}^a \upharpoonright m_0$. If $k > 0$, then let $L_k = |s^0| + \sum_{0 \leq j \leq k - 1} |v_j^0|$. Let $v_k^i$ be the string of length $m_k - L_k$ defined by $v_k^i(j) = w_{0^{k}\hat{\ }i}^0(L_k + j)$. 

Let $Z = \{z \in \cantorspace : (\exists x,y)(x,y \in X_{0^n}^0 \wedge w_{0^k\hat{\ }0} \subseteq x \wedge w_{0^k\hat{\ }1} \subseteq y \wedge z = x_{\geq \ell_{k + 1}} = y_{\geq \ell_{k + 1}})\}$. $Z$ is a nonempty $\lanalytic(z)$. For $t \in {}^{k + 1}2$ and $a < n$, let $B_t^a = \{x \in X_{t \upharpoonright k} : (\exists z)(z \in Z \wedge x = w_{t}^a \hat{\ } z)\}$. $B_t^a$ is a nonempty $\lanalytic(z)$ set and $(B_t^a)_{\geq \ell_{k + 1}} = (B_u^b)_{\geq \ell_{k + 1}}$ for all $a,b \in n$ and $t,u \in {}^{k + 1}2$. Apply Lemma \ref{same tail meet dense same tail} to get $X_t^a \subseteq B_t^a$ so that $X_t^a \in D_{k + 1}$ and $(X_t^a)_{\geq \ell_{k + 1}} = (X_u^b)_{\geq \ell_{k + 1}}$. This completes the construction.
\end{proof}

\section{$\cantorspace \slash E_0$ Has the $2$-J\'onsson Property}\label{R/E0 has 2-Jonsson property}
\Begin{theorem}{2-Mycielski property for E0}
(Holshouser-Jackson) $E_0$ has the $2$-Mycielski property.
\end{theorem}

This can be proved by producing an $E_0$-tree $p$ using an $E_0$-tree fusion argument to ensure that $[[p]]_{E_0}^2$ meets a fixed countable sequence of dense (topologically) open subsets of ${}^2(\cantorspace)$. The fusion argument is quite similar to the fusion argument used in the forcing style proof of the $2$-J\'onsson property for $\cantorspace \slash E_0$ in this section.

By Theorem \ref{E0 tree inside set same saturation}, if $(A_a : a < n)$ is a sequence of $\analytic$ sets so that $[A_a]_{E_0} = [A_b]_{E_0}$ for all $a,b < n$, then there is a sequence of $E_0$-trees $(p_a : a \in n)$ so that for all $a,b < n$ and $i <2$, $|s^{p_a}| = |s^{p_b}|$ and $v_i^{p_a} = v_i^{p_b}$. This motivates the definition of the following forcings:

\Begin{definition}{forcing with E0 tree with same saturation}
Let $n > 0$, let $\widehat\bbP_{E_0}^n$ be the collection of $n$-tuples of $E_0$-trees $(p_0, ..., p_{n - 1})$ so that for all $a,b < n$ and $i < 2$, $|s^{p_a}| = |s^{p_b}|$ and $v_i^{p_a} = v_i^{p_b}$. Let $\leq_{\hatPE^n}$ be coordinatewise $\leq_{\bbP_{E_0}}$. Let $1_{\hatPE^n} = (1_{\bbP_{E_0}}, ..., 1_{\bbP_{E_0}})$. $(\hatPE^n, \leq_{\hatPE^n}, 1_{\hatPE^n})$ is forcing with $n$ $E_0$-trees with the same $E_0$-saturation.

Let $x_\mathrm{gen}^0$ and $x_\mathrm{gen}^1$ be the $\hatPE^2$ names for the left and right generic real added by a generic filter for $\hatPE^2$.
\end{definition}

\Begin{definition}{hat PE0 pruning notation}
If $p, q \in \PE$, then let $p \leq_{\PE}^n q$ be defined in the same way as $p \leq_{\bbS}^n q$, when $p$ and $q$ are considered as conditions in Sacks forcing $\bbS$. 

Let $\leq_{\hatPE^2}^n$ be the coordinate-wise ordering using $\leq_{\PE}^n$.

A sequence $\langle p_n : n \in \omega\rangle$ of conditions of $\PE$ is a fusion sequence if and only if $p_{n + 1} \leq_{\PE}^n p_n$ for all $n \in \omega$. Similarly, a sequence $\langle (p_n,q_n) : n \in \omega\rangle$ of conditions in $\hatPE^2$ is a fusion sequence if and only if $(p_{n + 1},q_{n + 1}) \leq_{\hatPE^2}^n (p_n,q_n)$ for all $n \in \omega$.

Suppose $p \in \PE$. Let $n \in \omega$ with $n > 0$. Let $u,v \in {}^n 2$ with $u(n - 1) \neq v(n - 1)$. Suppose $(p',q') \in \hatPE^2$ with the property that $p' \leq_{\PE} \Xi(p,u)$ and $q' \leq_{\PE} \Xi(p,v)$. Let $A = \{s \in {}^n 2 : s(n - 1) = u(n - 1)\}$ and $B = \{s \in {}^n 2 : s(n - 1) = v(n - 1)\}$. Then define $\prune^{(u,v)}_{(p',q')}(p) \in \PE$ by
$$\prune^{(u,v)}_{(p',q')}(p) = \bigcup_{t \in A} \switch_{s^{\Xi(p,t)}}(p') \cup \bigcup_{t \in B} \switch_{s^{\Xi(p,t)}}(q')$$

Suppose $(p,q) \in \hatPE^2$. Let $n \in \omega$, $n > 0$, and $u,v \in {}^n 2$ so that $u(n - 1) \neq v(n - 1)$. Let $(p',q') \leq_{\hatPE^n} (p,q)$ so that $p' \leq_{\PE} \Xi(p,u)$ and $q' \leq_{\PE} \Xi(q,v)$. Define $\tprune^{(u,v)}_{(p',q')}(p,q)$ by
$$\tprune_{(p',q')}^{(u,v)}(p,q) = \Big(\prune^{(u,v)}_{(p', \switch_{s^{p}}(q'))}(p), \prune^{(u,v)}_{(\switch_{s^q}(p'),q')}(q)\Big)$$
Perhaps more concretely: Let $A = \{s \in {}^n2 : s(n - 1) = u(n - 1)\}$ and $B = \{s \in {}^n2 : s(n - 1) = v(n - 1)\}$. 
$$\tprune_{(p',q')}^{(u,v)}(p,q) = \Big(\bigcup_{t \in A} \switch_{s^{\Xi(p,t)}}(p') \cup \bigcup_{t \in B} \switch_{s^{\Xi(p,t)}}(q'), \bigcup_{t \in A} \switch_{s^{\Xi(q,t)}}(p') \cup \bigcup_{t \in B} \switch_{s^{\Xi(q,t)}}(q')\Big)$$
\end{definition}

\Begin{fact}{fusion and prune of hat PE0}
Suppose $(p,q)$, $(p',q')$, $n$, $u$, and $v$ are as in Definition \ref{hat PE0 pruning notation}. Then $\tprune_{(p',q')}^{(u,v)}(p,q) \in \hatPE^2$ and if $\tprune_{(p',q')}^{(u,v)}(p,q) = (x,y)$, then $s^x = s^p$ and $s^y = s^q$. 

If $|u| = n$, then $\tprune_{(p',q')}^{(u,v)}(p,q) \leq_{\hatPE^2}^n (p,q)$. 

If $\langle p_n : n \in \omega\rangle$ is a fusion sequence of conditions in $\PE$, then $p_\omega = \bigcap_{n \in \omega} p_n$ is a condition in $\PE$ and is called the fusion of the fusion sequence.

Similarly, if $\langle (p_n,q_n) : n \in \omega\rangle$ is a fusion sequence of conditions in $\hatPE^2$, then $(p_\omega,q_\omega) = (\bigcap_{n \in \omega} p_n, \bigcap_{n \in \omega} q_n)$ is a condition in $\hatPE^2$ and is called the fusion of the fusion sequence.
\end{fact}

\Begin{fact}{2-hat E0 forcing is proper}
(\cite{Canonical-Ramsey-Theory-on-Polish-Spaces} Proposition 7.6) $\hatPE^2$ is a proper forcing. In fact, for any countable model $M $ and $(p,q) \in (\hatPE^2)^M$, there is a $(p',q') \leq_{\hatPE^2} (p,q)$ which is a $(M,\hatPE^2)$-master condition and $[p'] \times_{E_0} [q']$ consists of pairs of reals which are $\hatPE^2$-generic over $M$. 

Moreover, if $\tau$ is a $\hatPE^2$-name in $M$ such that $(p,q) \forces_{\hatPE^2}^M \tau \in \cantorspace$, then one can even find $(p',q')$ with the above properties and so that either

\noindent (i) there is some $z \in \cantorspace \cap M$ with $z \ E_0 \ \tilde 0$ so that $(p',q') \forces_{\hatPE^2} \tau = \check z$ 

\noindent or

\noindent (ii) $(p',q') \forces_{\hatPE^2} \neg(\tau \ E_0 \ \tilde 0)$. 
\end{fact}

\begin{proof}
Let $(D_n : n \in \omega)$ be a decreasing sequence of dense open subsets of $\hatPE^2$ in $M$ with the property that if $D$ is a dense open subset of $\hatPE^2$ in $M$, then there is some $k \in \omega$ so that $D_k \subseteq D$. 

There are two cases: (Note that $\tilde 0$ can be replaced by any real in $M$ in the following argument and hence as well in the statement of the fact.)

(Case I) There is some $(p',q') \leq_{\hatPE^2} (p,q)$ so that $(p',q') \forces_{\hatPE^2}^M \tau \ E_0 \ \tilde 0$. Then there is some $z \in (\cantorspace)^M$ and $(p'',q'')$ so that $(p'',q'') \forces_{\hatPE^2}^M \tau = \check z$. Since $D_0$ is dense open in $\hatPE^2$, one may assume $(p'',q'') \in D_0$. Let $(p_0,q_0) = (p'',q'')$. 

(Case II) $(p,q) \forces_{\hatPE^2}^M \neg(\tau \ E_0 \ \tilde 0)$. Let $(p',q')$ be any condition below $(p,q)$ so that $(p',q') \in D_0$. Let $(p_0, q_0) = (p',q')$. 

In either case, a fusion sequence $((p_n, q_n) : n \in \omega)$ of conditions in $(\hatPE^2)^M$ will be constructed with the following property:

\noindent (i) For all $n \in \omega$, $(u,v) \in {}^n 2 \times {}^n 2$ so that $u(n - 1) \neq v(n - 1)$, $(\Xi(p_n, u), \Xi(q_n, v)) \in D_{n}$. 

Suppose this can be done to produce a fusion sequence $\langle (p_n, q_n) : n \in \omega\rangle$. Let $(p',q')$ be the fusion of this fusion sequence. First, it will be shown that (i) implies that $(p',q')$ is a $(M,\hatPE^2)$-master condition with the property that $[p'] \times_{E_0} [q']$ consists entirely of pairs of reals which are $\hatPE^2$-generic over $M$. 

Let $D$ be a dense open subset of $\hatPE^2$ with $D \in M$. By the choice of $(D_n : n \in \omega)$, there is some $k \in \omega$ so that $D_k \subseteq D$. Let $G \subseteq \hatPE^2$ be $\hatPE^2$-generic over $M$ containing $(p',q')$. Let $K = |\Lambda(p_k,0^k)|$. Let $E_K$ be the collection of $(p,q) \in \hatPE^2$ so that $|s^p| > K$ and there is some $j$ with $K < j < |s^p|$ so that $s^p(j) \neq s^q(j)$. $E_K$ is dense open. Since $G$ is generic, $G \cap E_K \neq \emptyset$. As $G$ is generic and $(p',q') \in G$, one may assume that there is some $(p'',q'') \leq_{\hatPE^2} (p',q')$ with $(p'',q'') \in G \cap E_k$. So there is some $J > k$ and $u,v \in {}^J2$ with $u(J - 1) \neq v(J - 1)$ so that $(p'',q'') \leq_{\hatPE^2} (\Xi(p',u), \Xi(q',v)) \leq_{\hatPE^2} (\Xi(p_J,u), \Xi(q_J,v))$. Since $(p'',q'') \in G$ and $G$ is a filter, $(\Xi(p_J,u), \Xi(q_J,v)) \in G$. However $(\Xi(p_J, u),\Xi(q_J,v)) \in D_J \subseteq D_k$ by (i). Note that $(\Xi(p_J,u), \Xi(q_J,v)) \in M$. Since $G$ was an arbitrary generic containing $(p',q')$, it has been shown that $(p',q') \forces_{\hatPE^2}^V \dot G \cap \check M \cap \check D \neq \emptyset$. $(p',q')$ is a $(M,\hatPE^2)$-master condition. $\hatPE^2$ is proper.

Now suppose $(a,b) \in [p'] \times_{E_0} [q']$. Let $G_{(a,b)} = \{(p,q) \in \hatPE^2 \cap M : (a,b) \in [p]\times[q]\}$. Let $D$ be a dense open set. There is some $k$ so that $D_k \subseteq D$. Since $\neg(a \ E_0 \ b)$, there is some $j > k$ and some $u,v \in {}^j2$ with $u(j - 1) \neq v(j - 1)$ so that $\Lambda(p',u) \subseteq a$ and $\Lambda(q',v) \subseteq b$. Then $(a,b) \in [\Xi(p_j,u)] \times_{E_0} [\Xi(q_j,v)]$. Therefore $(\Xi(p_j,u), \Xi(q_j,v)) \in G_{(a,b)} \cap D_k \subseteq G_{(a,b)} \cap D$. $G_{(a,b)}$ is $\hatPE^2$-generic over $M$. $(a,b)$ is a $\hatPE^2$-generic pair of reals over $M$. 

It is clear using the forcing theorems that if Case I holds, then statement (i) of the fact holds and if Case II holds, then statement (ii) of the fact holds. 

Now it remains to construct the fusion sequence:

$(p_0, q_0)$ is already given depending on the case. The rest of the construction is the same for both cases.

Suppose $(p_n,q_n)$ have been constructed with the desired properties. For some $J \in \omega$, let $(u_0,v_0)$, ..., $(u_{J-1}, v_{J-1})$ enumerate all $(u,v) \in {}^{n}2 \times {}^n2$ with $u(n -1) \neq v(n - 1)$.

A sequence of conditions in $\hatPE^2$, $(x_{-1},y_{-1})$, ..., $(x_{J - 1}, y_{J - 1})$ will be constructed:

Let $(x_{-1},y_{-1}) = (p_n,q_n)$. 

Suppose $(x_k,y_k)$ has been constructed for $k < J -1$.

Since $D_{n + 1}$ is dense open, find some $(p',q') \leq_{\hatPE^2} (\Xi(x_k, u_k), \Xi(y_k,v_k))$ so that $(p',q') \in D_{n + 1}$. Let $(x_{k + 1}, y_{k +1}) = \tprune^{(u_k,v_k)}_{(p',q')}(x_k,y_k)$.
\end{proof}

\Begin{lemma}{2-Jonsson R mod E0 lemma}
Let $f : [\cantorspace]^{2}_{E_0} \rightarrow \cantorspace$. Suppose there is a countable model $M$ of some sufficiently large fragment of $\mathsf{ZF}$, $(p,q) \in \hatPE^2 \cap M$, and $\tau \in M^{\hatPE^2}$ so that $(p,q) \forces_{\hatPE^2}^M \tau \in \cantorspace$ and whenever $G \subseteq \hatPE^2$ is $\hatPE^2$-generic over $M$ with $(p,q) \in G$, then $\tau[G] = f(x_\mathrm{gen}^0[G], x_\mathrm{gen}^1[G])$. Then there is a $(p',q') \leq_{\hatPE^2} (p,q)$ so that $[f[[p'] \times_{E_0} [q']]]_{E_0} \neq \cantorspace$. 
\end{lemma}

\begin{proof}
Let $(p',q') \in \hatPE^2$ be given by Fact \ref{2-hat E0 forcing is proper}. Then exactly one of the two happens:

(i) For all $G \subseteq \hatPE^2$ which are generic over $M$,  $M[G] \models \tau[G] \ E_0 \ \tilde 0$. By absoluteness, $\tau[G] \ E_0 \ \tilde 0$.

(ii) For all $G \subseteq \hatPE^2$ which are generic over $M$,  $M[G] \models \neg(\tau[G] \ E_0 \ \tilde 0)$. By absoluteness, $\neg(\tau[G] \ E_0 \ \tilde 0)$.

Let $(a,b) \in [p'] \times_{E_0} [q']$. By Fact \ref{2-hat E0 forcing is proper}, there is some $\hatPE^2$-generic filter $G_{(a,b)}$ so that $x_\mathrm{gen}^0[G_{(a,b)}] = a$ and $x_\mathrm{gen}^1[G_{(a,b)}] = b$.

If (i) holds, then $f(a,b) = f(x_\mathrm{gen}^0[G_{(a,b)}], x_\mathrm{gen}^1[G_{(a,b)}]) = \tau[G_{(a,b)}]$ which is $E_0$ related to $\tilde 0$. So $[f[[p'] \times_{E_0} [q']]]_{E_0} = [\tilde 0]_{E_0} \neq \cantorspace$.

If (ii) holds, then $f(a,b) = f(x_\mathrm{gen}^0[G_{(a,b)}], x_\mathrm{gen}^1[G_{(a,b)}]) = \tau[G_{(a,b)}]$, but $\neg(\tau[G_{(a,b)}] \ E_0 \ \tilde 0)$. So $\tilde 0 \notin [f[[p'] \times_{E_0} [q']]]_{E_0}$. 
\end{proof}

\Begin{theorem}{BP 2-Jonsson R mod E0}
(Holshouser-Jackson) $(\mathsf{ZF + DC + AD})$ $\cantorspace \slash E_0$ has the $2$-J\'onsson property.
\end{theorem}

\begin{proof}
Let $F : [\cantorspace \slash E_0]^2_= \rightarrow \cantorspace \slash E_0$. Define the relation $R \subseteq \cantorspace \times \cantorspace \times \cantorspace$ by 
$$R(x,y,z) \Leftrightarrow z \in F([x]_{E_0}, [y]_{E_0})$$
$\mathsf{AD}$ can prove comeager uniformization (Fact \ref{comeager uniformization}): There is a comeager set $C \subseteq [\cantorspace]^2_{E_0}$ and a continuous function $f: C \rightarrow \cantorspace$ which uniformizes $R$ on $C$. Since $E_0$ has the $2$-Mycielski property, let $p$ be an $E_0$-tree so that $[[p]]^2_{E_0} \subseteq C$. So $f \upharpoonright [[p]]^2_{E_0}$ is a continuous function. By a similar argument as in Fact \ref{continuity satisfy sacks name condition}, one can find a name $\tau$ satisfying Lemma \ref{2-Jonsson R mod E0 lemma} using the condition $(p,p) \in \hatPE^2$. Then Lemma \ref{2-Jonsson R mod E0 lemma} gives some $(p',q') \leq_{\hatPE^2} (p,p)$ so that $[f[[p']\times_{E_0}[q']]]_{E_0}$ is either $[\tilde 0]_{E_0}$ or does not contain $\tilde 0$. Since $(p',q') \in \hatPE^2$, $[p']$ and $[q']$ have the same $E_0$-saturation. Let $A = [[p']]_{E_0} = [[q']]_{E_0}$. Note that $E_0 \equiv_\borel E_0 \upharpoonright A$. Let $ B = A \slash E_0$. There is a bijection between $B$ and $\cantorspace \slash E_0$. Moreover, $F[[B]_=^2]$ is either the singleton $\{[\tilde 0]_{E_0}\}$ or is missing the element $[\tilde 0]_{E_0}$. In either case, $F[[B]^2_=] \neq \cantorspace \slash E_0$. 
\end{proof}

\section{Partition Properties of $\cantorspace \slash E_0$ in Dimension 2}\label{Partition Properties of R/E0 in dimension 2}

\Begin{definition}{partition relations}
Let $X$ and $Y$ be sets. Let $n \in \omega$. Denote $X \rightarrow (X)^n_Y$ to mean that for any function $f : \mathscr{P}^n(X) \rightarrow Y$, there is some $Z \subseteq X$ with $Z \approx X$ and $|f[\mathscr{P}^n(Z)]| = 1$. 

Denote $X \mapsto (X)^n_Y$ to mean that for any function $f : [X]^n_= \rightarrow Y$, there is some $Z \subseteq X$ with $Z \approx X$ and $|f[[Z]^n_=]| = 1$. 
\end{definition}

\Begin{fact}{partition relation for R}
(Galvin) Assuming $\mathsf{ZF}$ and all sets of reals have the Baire property, $\cantorspace \rightarrow (\cantorspace)^2_n$ for all $n \in \omega$. 
\end{fact}

Note that $\cantorspace \mapsto (\cantorspace)^2_2$ is not true: If $x,y \in \cantorspace$ and $x \neq y$, then define $d(x,y) = \min\{n : x(n) \neq y(n)\}$. Define $f: [\cantorspace]^2_= \rightarrow 2$ by $f(x,y) = x(d(x,y))$. Note that $f(x,y) \neq f(y,x)$. It is impossible to find a homogeneous set for this coloring of $[\cantorspace]^2_=$.  

However, under $\mathsf{AD}$, $\cantorspace \slash E_0 \mapsto (\cantorspace \slash E_0)^2_n$ does hold:

\Begin{lemma}{definability implies use one color in partition}
Let $n > 1$. Let $F : [\cantorspace]^2_{E_0} \rightarrow n$ be a function. Suppose there is a countable model $M$ of some sufficiently large fragment of $\mathsf{ZF}$, $(p,q) \in \hatPE^2 \cap M$, and $\tau$ is a $\hatPE^2$-name in $M$ so that $(p,q) \forces_{\hatPE^2}^M \tau \in \check n$ and whenever $G \subseteq \hatPE^2$ is $\hatPE^2$-generic over $M$ with $(p,q) \in G$, $\tau[G] = F(x_\mathrm{gen}^0[G],x_\mathrm{gen}^1[G])$. Then there is a $(p',q') \leq_{\hatPE^2} (p,q)$ so that $|F[ [p'] \times_{E_0} [q']]| = 1$.
\end{lemma}

\begin{proof}
Since $(p,q) \forces_{\hatPE^2}^M \tau \in \hat n$, there is some $(r,s) \leq_{\hatPE^2} (p,q)$ and some $k \in n$ so that $(r,s) \forces_{\hatPE^2}^M \tau = \check k$. Using Fact \ref{2-hat E0 forcing is proper}, let $(p',q') \leq_{\hatPE^2} (r,s)$ be a $(M,\hatPE^2)$-master condition so that $[p'] \times_{E_0} [q']$ consists of pairs of reals which are $\hatPE^2$-generic over $M$. Using the forcing theorem and the assumptions, $(p',q')$ works. 
\end{proof}

\Begin{theorem}{E0 dimension 2 partition property}
$(\mathsf{ZF + DC + AD})$ $\cantorspace \slash E_0 \mapsto (\cantorspace \slash E_0)^2_n$ for all $n \in \omega$.
\end{theorem}

\begin{proof}
Let $f : [\cantorspace \slash E_0]_=^2 \rightarrow n$. Define a relation $R \subseteq \cantorspace \times \cantorspace \times n$ by
$$R(x,y,k) \Leftrightarrow k = f([x]_{E_0},[y]_{E_0}).$$
$\mathsf{AD}$ can prove comeager uniformization (Fact \ref{comeager uniformization}) so there is some comeager $C \subseteq [\cantorspace]^2_{E_0}$ and a continuous function $F : C \rightarrow n$ which uniformizes $R$ on $C$. Since $E_0$ has the $2$-Mycielski property, let $p$ be an $E_0$-tree so that $[[p]]^2_{E_0} \subseteq C$. $F \upharpoonright [[p]]^2_{E_0}$ is a continuous function.  As before, one can find a $\hatPE^2$-name $\tau$ satisfying Lemma \ref{definability implies use one color in partition} using the the condition $(p,p)$. Then using Lemma \ref{definability implies use one color in partition}, there is a $(r,s) \leq_{\hatPE^2} (p,p)$ so that $|F[[r] \times_{E_0} [s]]| = 1$. Let $k$ be the unique element in this set.

Let $A = [r]_{E_0} = [s]_{E_0}$. $A \slash E_0 \approx \cantorspace \slash E_0$. Suppose $x,y \in A \slash E_0$ and $x \neq y$. There is some $a,b \in A$ with $a \in [r]$, $b \in [s]$, $a \in x$, and $b \in y$. $F(a,b) = k$. Therefore, $R(a,b,k)$. By definition, $f(x,y) = k$. So $|f[[A]^2_=]| = 1$. 
\end{proof}

\Begin{remark}{E0 dimension 2 partition property remark}
Before, it was mentioned that $\cantorspace \mapsto (\cantorspace)^2_2$ is not true. Using the notation from the above proof: Observe that the function $F$ produced in the above proof is $E_0$-invariant in the sense that if $x \ E_0 \ x'$ and $y \ E_0 \ y'$ then $F(x,y) = F(x',y')$. Also since $(r,s) \in \hatPE^2$, if $a \in A$, then there is some $a' \in [r]$ and $a'' \in [s]$ with $a \ E_0 \ a' \ E_0 \ a''$. These two facts are essential in proving $\cantorspace \slash E_0 \mapsto (\cantorspace \slash E_0)^2_n$. 

Later it will be shown that the partition relation for $\cantorspace \slash E_0$ will fail in higher dimension. The counterexample is closely connected to the failure of the $3$-J\'onsson property.
\end{remark}

\section{$E_0$ Does Not Have the 3-Mycielski Property}\label{E0 does not have 3-Mycielski property}

An earlier section mentioned that Holshouser and Jackson proved $E_0$ has the $2$-Mycielski property and $\cantorspace \slash E_0$ has the $2$-J\'onsson property. The next few sections will show that dimension $2$ is the best possible for these combinatorial properties. This section will show the failure of the $3$-Mycielski property and the weak $3$-Mycielski property for $E_0$.

\Begin{theorem}{E0 does not have Mycielski property}
Let $D \subseteq {}^3(\cantorspace)$ be defined by 
$$D = \{(x,y,z) \in {}^3(\cantorspace) : (\exists n)(x(n) \neq y(n) \wedge x(n) \neq z(n) \wedge y(n + 1) \neq z(n + 1))\}.$$
$D$ is dense open in ${}^3(\cantorspace)$.

If $p$ is an $E_0$-tree with associated objects $\{s, v_n^i : i \in 2 \wedge n \in \omega\}$, $\varphi$, and $\Phi$ as in Definition \ref{E0 tree}, then 
$$(\Phi(\widetilde{010}), \Phi(\widetilde{110}), \Phi(\tilde{1})) \notin D.$$

$E_0$ does not have the 3-Mycielski property.
\end{theorem}

\begin{proof}
Suppose $(x,y,z) \in D$. Then there is an $n \in \omega$ so that $x(n) \neq y(n)$, $x(n) \neq z(n)$, and $y(n + 1) \neq z(n + 1)$. Let $\sigma = x \upharpoonright (n + 2)$, $\tau = y \upharpoonright (n + 2)$, and $\rho = z \upharpoonright (n + 2)$. Then $N_{\sigma,\tau,\rho} \subseteq D$. $D$ is open. 

Suppose $\sigma,\tau,\rho \in \finBinarySequence$ and $|\sigma| = |\tau| = |\rho| = k$. Let $\sigma' = \sigma\hat{\ }00$, $\tau' = \tau\hat{\ }10$, and $\rho' = \rho\hat{\ }11$. Then $N_{\sigma',\tau',\rho'} \subseteq N_{\sigma,\tau,\rho}$ and $N_{\sigma',\tau',\rho'} \subseteq D$. $D$ is dense open.

Let $L_{-1} = |s|$. For $n \in \omega$, let $L_n = |s| + \sum_{k \leq n} |v_k^0|$. Note that if $x(n) = y(n)$, then for all $L_{n - 1} \leq k < L_n$, $\Phi(x)(k) = \Phi(y)(k)$. If $x(n) \neq y(n)$, then $\Phi(x)(L_{n - 1}) \neq \Phi(y)(L_{n - 1})$, and there may be other $L_{n - 1} \leq k < L_n$ so that $\Phi(x)(k) \neq \Phi(y)(k)$. 

Suppose $\Phi(\widetilde{010})(n) \neq \Phi(\widetilde{110})(n)$ and $\Phi(\widetilde{010})(n) \neq \Phi(\tilde 1)(n)$. Then there exists some $k$ so that $L_{3k - 1} \leq n < L_{3k}$. If $n \neq L_{3k} - 1$, then $n + 1 < L_{3k}$. Since $\widetilde{110}(3k) = 1 = \tilde{1}(3k)$, $\Phi(\widetilde{110})(n + 1) = \Phi(\tilde 1)(n + 1)$. Hence if $n \ne L_{3k} - 1$, then $n$ cannot be used to witness that $(\Phi(\widetilde{101}), \Phi(\widetilde{110}), \Phi(\tilde{1})) \in D$. Suppose $n = L_{3k} - 1$. Since $\widetilde{110}(3k + 1) = 1 = \tilde{1}(3k + 1)$, one has that $\Phi(\widetilde{110})(n + 1) = \Phi(\widetilde{110})(L_{3k}) = \Phi(\tilde 1)(L_{3k}) = \Phi(\tilde 1)(n + 1)$. If $n = L_{3k} - 1$, then $n$ does not witness membership in $D$. Hence $(\Phi(\widetilde{010}), \Phi(\widetilde{110}), \Phi(\tilde 1)) \notin D$. 

Since $\neg(\widetilde{010} \ E_0 \ \widetilde{110})$, $\neg(\widetilde{010} \ E_0 \ \tilde 1)$, and $\neg(\widetilde{110} \ E_0 \ \tilde 1)$, $(\Phi(\widetilde{010}), \Phi(\widetilde{110}), \Phi(\tilde 1)) \in [[p]]^3_{E_0}$. Hence $[[p]]^3_{E_0} \not\subseteq D$. 

Suppose $B$ is $\borel$ so that $E_0 \upharpoonright B \equiv_\borel E_0$. By Fact \ref{E0 characterization}, there is some $E_0$-tree $p$ so that $[p] \subseteq B$. By the above, $[B]^3_{E_0} \not\subseteq D$. $E_0$ does not have the 3-Mycielski property.
\end{proof}

The $3$-Mycielski property asks for a single $\borel$ set $A$ with $E_0 \leq_\borel E_0 \upharpoonright A$ so that $[A]^3_{E_0} = A \times_{E_0} A \times_{E_0} A$ is contained inside a comeager set. If one is interested in combinatorial properties of the quotient $\cantorspace \slash E_0$, such as the J\'onsson property, then one is only concerned with three sets $A$, $B$, and $C$ with the same $E_0$-saturation. With this consideration, the $3$-Mycielski property seems unnecessarily restrictive. The weak $3$-Mycielski property was defined to remove this demand.

One other curiosity of the $3$-Mycielski property is that Theorem \ref{E0 does not have Mycielski property} allows a (topologically) dense open subset of ${}^3(\cantorspace)$ to be a counterexample to the $3$-Mycielski property. Let $D \subseteq {}^3(\cantorspace)$ be any dense open set. There are three strings $\sigma$, $\tau$, and $\gamma$ of the same length so that $N_{\sigma,\gamma,\tau} \subseteq D$. Let $p,q,r$ be any three perfect $E_0$-trees so that $s^p = \sigma$, $s^q = \tau$, $s^r = \gamma$, and for all $n \in \omega$ and $i \in 2$, $v_n^{i,p} = v_n^{i,q} = v_n^{i,r}$. Then $[p] \times_{E_0} [q] \times_{E_0} [r] \subseteq D$. Also $[[p]]_{E_0} = [[q]]_{E_0} = [[r]]_{E_0}$. So no dense open set can be a counterexample to the weak $3$-Mycielski property. 

Using the more informative structure theorem for $E_0$ proved above and the argument in Theorem \ref{E0 does not have Mycielski property}, a comeager subset of ${}^3(\cantorspace)$ is used to show $E_0$ does not have the weak $3$-Mycielski property.

\Begin{theorem}{E0 does not have weak 3-mycielski property}
For each $k \in \omega$, let 
$$D_k = \{(x,y,z) \in {}^3(\cantorspace) : (\exists n \geq k)(x(n) \neq y(n) \wedge x(n) \neq z(n) \wedge y(n + 1) \neq z(n + 1))\}.$$
Each $D_k$ is dense open.

Let $C = \bigcap_{n \in \omega} D_n$. $C$ is comeager. 

For any $\borel$ sets $A_0$, $A_1$, and $A_2$ such that 

\noindent (I) $E_0 \leq_{\borel} E_0 \upharpoonright A_0$, $E_0 \leq_{\borel} E_0 \upharpoonright A_1$, $E_0 \leq_{\borel} E_0 \upharpoonright A_2$

\noindent (II) $[A_0]_{E_0} = [A_1]_{E_0} = [A_2]_{E_0}$, 

$A_0 \times_{E_0} A_1 \times_{E_0} A_2 \not\subseteq C$. 

$C$ does not have the weak $3$-Mycielski property.
\end{theorem}

\begin{proof}
Let $A_0, A_1, A_2$ be any three $\borel$ sets so that $E_0 \leq_\borel E_0 \upharpoonright A_i$ and have the same $E_0$-saturation. By Theorem \ref{E0 tree inside set same saturation}, there are $E_0$-trees, $p$, $q$, and $r$ so that 

(i) $|s^p| = |s^q| = |s^r| = k$ 

(ii) $v_n^{i,p} = v_n^{i,q} = v_n^{i,r}$ for all $n \in \omega$ and $i \in 2$

(iii) $[p] \subseteq A_0$, $[q] \subseteq A_1$, and $[r] \subseteq A_2$. 

Note that the only differences among the three $E_0$-trees occurs in the stems. Hence by the same argument as in Theorem \ref{E0 does not have Mycielski property}, $[p] \times_{E_0} [q] \times_{E_0} [r] \not\subseteq D_k$. Hence $A_0 \times_{E_0} A_1 \times_{E_0} A_2 \not\subseteq D_k$. So $A_0 \times_{E_0} A_1 \times_{E_0} A_2 \not\subseteq C$.
\end{proof}

\section{Surjectivity and Continuity Aspects of $E_0$}\label{surjectivity and continuity aspects of E0}

From Holshouser and Jackson's proof of the J\'onsson property for $\cantorspace$, the Mycielski property was used primarily to show an arbitrary function $f$ on ${}^n(\cantorspace)$ could be continuous on $[[p]]_=^n$ for some perfect tree $p$. Using the continuity of $f \upharpoonright [[p]]_=^n$, they show that there is some perfect subtree $q \subseteq p$ so that $f[[[q]]_=^n] \neq \cantorspace$.

As the previous section shows that $E_0$ does not have the $3$-Mycielski property, it is natural to ask if by some other means it is possible to find for any $f : {}^3(\cantorspace) \rightarrow \cantorspace$, some $\borel$ set $A$ so that $E_0 \leq_\borel E_0 \upharpoonright A$ and $f \upharpoonright [A]^3_{E_0}$ is continuous. Also if the function $f \upharpoonright [A]^3_{E_0}$ is continuous, is it possible to find some $\borel$ $B \subseteq A$ with $E_0 \leq_\borel E_0 \upharpoonright B$ so that $f[[B]^3_{E_0}]$ does not meet all $E_0$ equivalence classes? This section will provide an example to show both of these properties can fail. This example will also be modified in the next section to show the failure of the $3$-J\'onsson property for $E_0$. 

\Begin{fact}{function to 3w surjective E0 tree}
Let $A = \{x \in {}^\omega 3 : (\forall n)(x(n) \neq x(n + 1))\}$. There is a continuous function $P : [\cantorspace]^3_{E_0} \rightarrow A$ so that for any $E_0$-tree $p$, $P[[[p]]_{E_0}^3] = A$. Moreover, if $p$, $q$, and $r$ are $E_0$-trees so that $|s^p| = |s^q| = |s^r|$ and for all $i \in 2$ and $n \in \omega$, $v_n^{i,p} = v_n^{i,q} = v_n^{i,r}$, then $P[[[p]] \times_{E_0} [[q]] \times_{E_0} [[r]]]$ meets all $E_0$-classes of $A$, where the latter $E_0$ is defined on ${}^\omega 3$. 
\end{fact}

\begin{proof}
Let $(x,y,z) \in [\cantorspace]^3_{E_0}$. Let $L_0$ be the largest $N \in \omega$ so that $x\upharpoonright N = y\upharpoonright N = z \upharpoonright N$. Define 
$$a_0 = \begin{cases}
0 & \quad x(L_0) = y(L_0) \\
1 & \quad x(L_0) = z(L_0) \\
2 & \quad y(L_0) = z(L_0)
\end{cases}.$$

Suppose $L_n$ and $a_n$ have been defined. Let $L_{n + 1}$ be the smallest $N > L_n$ so that $x(N) \neq y(N)$ if $a_n = 0$, $x(N) \neq z(N)$ if $a_n = 1$, and $y(N) \neq z(N)$ if $a_n = 2$. Define
$$a_{n + 1} = \begin{cases}
0 & \quad x(L_{n + 1}) = y(L_{n + 1}) \\
1 & \quad x(L_{n + 1}) = z(L_{n + 1}) \\
2 & \quad y(L_{n + 1}) = z(L_{n +1 })
\end{cases}.$$

Define $P(x,y,z) \in A$ by $P(x,y,z)(n) = a_n$. $P$ is continuous.

Now let $p$ be an $E_0$ tree. Let $s$ and $v_n^i$, for $n \in \omega$ and $i \in 2$, be associated with the $E_0$-tree $p$. Let $\Phi : \cantorspace \rightarrow [p]$ be the canonical homeomorphism.

Let $v \in A$. Let
$$(a_i,b_i) = \begin{cases}
(0,1) & \quad v(i) = 0 \\
(1,0) & \quad v(i) = 1 \\
(1,1) & \quad v(i) = 2
\end{cases}.$$
Let $a,b \in \cantorspace$ be defined by $a(n) = a_n$ and $b(n) = b_n$. Then $(\Phi(\tilde 0), \Phi(a), \Phi(b)) \in [[p]]^3_{E_0}$ and $P((\Phi(\tilde 0), \Phi(a), \Phi(b))) \allowbreak = v$. Hence $P[[[p]]^3_{E_0}] = A$. The second statement is proved similarly after noting the three $E_0$-trees are the same after their stems.
\end{proof}

\Begin{theorem}{cont function surjective E0 tree}
There is a continuous $Q : [\cantorspace]^3_{E_0} \rightarrow \cantorspace$ so that for any $E_0$-tree $p$, $Q[[[p]]^3_{E_0}] = \cantorspace$. 
\end{theorem}

\begin{proof}
Let $t_0 = 00$, $t_1 = 01$, and $t_2 = 10$. 

Let $Q' : {}^\omega 3 \rightarrow \cantorspace$ be defined by $Q'(x) = t_{x(0)}\hat{\ } t_{x(1)} \hat{\ } ...$. $Q'$ is a continuous injection.

$Q'[A]$ is a perfect subset of $\cantorspace$. $Q'[A] = [T]$ for some perfect tree $T$. Let $Q'' : Q'[A] \rightarrow \cantorspace$ be the continuous bijection naturally induced by $T$. 

Let $Q = Q'' \circ Q' \circ P$. $Q$ has the desired property.
\end{proof}

\Begin{corollary}{borel function surjective on any power E0}
There is a $\borel$ function $K : {}^3(\cantorspace) \rightarrow \cantorspace$ so that on any $\analytic$ set $A$ so that $E_0 \leq_\borel E_0 \upharpoonright A$, $K [[A]^3_{E_0}] = \cantorspace$. Moreover, for any $\analytic$ sets $A_0$, $A_1$, and $A_2$ so that for all $i < 3$, $E_0 \leq_\borel E_0 \upharpoonright A_i$ and $[A_0]_{E_0} = [A_1]_{E_0} = [A_2]_{E_0}$, $[K[\prod_{i < 3}^{E_0} A_i]]_{E_0} = \cantorspace$. 
\end{corollary}

\Begin{fact}{not continuous on any 3 product E0 product into A}
There is a $\borel$ function $P' : {}^3(\cantorspace) \rightarrow A$ so that for all $E_0$-tree $p$, $P' \upharpoonright [[p]]^3_{E_0}$ is not continuous.
\end{fact}

\begin{proof}
Let $P$ be the function from Fact \ref{function to 3w surjective E0 tree}. Define $P'$ by
$$P'(x,y,z) = \begin{cases}
\widetilde{01} & \quad (x,y,z) \notin [\cantorspace]^3_{E_0}\\
P(x,y,z) & \quad (\forall k)(\exists n > k)(P(x,y,z)(n) = 2) \\
\widetilde{01} & \quad (\exists k)(\forall n > k)(P(x,y,z)(n) < 2)
\end{cases}.$$

Suppose $P'$ is continuous on some $[[p]]_{E_0}^3$. Let $s \in \finBinarySequence$ be so that for all $n < |s| - 1$, $s(n) \neq s(n + 1)$ and there exists some $n < |s|$ so that $s(n) = 2$. By continuity, $(P')^{-1}[N_s] \cap [[p]]^3_{E_0}$ is open in $[[p]]_{E_0}^3$. There is some $u,v,w \in \finBinarySequence$ so that $|u| = |v| = |w|$ and $N_{\varphi(u), \varphi(v), \varphi(w)} \cap [[p]]_{E_0}^3 \subseteq (P')^{-1}[N_s] \cap [[p]]_{E_0}^3$. 

Let $x = u\hat{\ }\tilde{0}$, $y = v\hat{\ }\widetilde{01}$, and $z = w\hat{\ }\widetilde{10}$. Then $(\Phi(x),\Phi(y),\Phi(z)) \in (P')^{-1}[N_s] \cap [[p]]^3_{E_0}$. However, there is a $k$ so that for all $n > k$, $P(\Phi(x),\Phi(y),\Phi(z))(n) < 2$. Therefore, $P'(\Phi(x),\Phi(y),\Phi(z)) = \widetilde{01}$. However, $\widetilde{01} \notin N_s$ since there is some $n$ so that $s(n) = 2$. 
\end{proof}

\Begin{theorem}{not continuous on any 3 product}
There is a $\borel$ function $K : {}^3(\cantorspace) \rightarrow \cantorspace$ so that for all $\analytic$ sets $A$ so that $E_0 \leq_\borel E_0 \upharpoonright A$, $K\upharpoonright [A]^3_{E_0}$ is not continuous. 
\end{theorem}

\begin{proof}
Let $Q'$ be the function from the proof of Theorem \ref{cont function surjective E0 tree}. Let $P'$ be the function from Fact \ref{not continuous on any 3 product E0 product into A}. $K = Q' \circ P'$ works.
\end{proof}

As a consequence, one has another proof of the failure of the $3$-Mycielski property for $E_0$. 

\Begin{corollary}{failure of continuity for E0 failure 3-Mycielski for E0}
$E_0$ does not have the $3$-Mycielski property.
\end{corollary}

\begin{proof}
Let $C \subseteq {}^3(\cantorspace)$ be any comeager set so that $K \upharpoonright C$ is a continuous function, where $K$ is from Fact \ref{not continuous on any 3 product}. Then $C$ witnesses the failure of the 3-Mycielski property for $E_0$. 
\end{proof}

\section{$\cantorspace \slash E_0$ Does Not Have the $3$-J\'onsson Property}\label{R/E0 does not have 3-Jonsson property}

\Begin{definition}{tail equivalence relation}
For $n \in \omega$, let $\Etail^n$ be the equivalence relation defined on ${}^\omega n$ by $x \ \Etail^n \ y$ if and only if $(\exists r)(\exists s)(\forall a)(x(r + a) = y(s + a))$. 
\end{definition}

\Begin{fact}{function P into A is E0 to tail invariant}
The function $P:[\cantorspace]^3_= \rightarrow A$ from Fact \ref{function to 3w surjective E0 tree} is $E_0$ to $\Etail^3$ invariant, which means for all $(x,y,z), (x',y',z') \in [\cantorspace]^3_{E_0}$ such that $x \ E_0 \ x'$, $y \ E_0 \ y'$, and $z \ E_0 \ z'$, $P(x,y,z) \ \Etail^3 \ P(x',y',z')$. 
\end{fact}

\begin{proof}
Using the notation from Fact \ref{function to 3w surjective E0 tree}, let $(L_k : k \in \omega)$ and $(a_k : k \in \omega)$ be the $L$ and $a$ sequences for $(x,y,z)$ and let $(J_k : k \in \omega)$ and $(b_k : k \in \omega)$ be the $L$ and $a$ sequences for $(x',y',z')$. 

Let $M \in \omega$ be so that $x_{\geq M} = x'_{\geq M}$, $y_{\geq M} = y'_{\geq M}$, and $z_{\geq M} = z'_{\geq M}$. Let $r \in \omega$ be largest so that $L_r < M$, and let $s \in \omega$ be largest so that $J_s < M$. 

(Case I) Suppose $L_{r + 1} = J_{s + 1}$. Then $a_{r + 1} = b_{s + 1}$. Since for all $n \geq 1$, $L_{r + n}, J_{s + n} \geq M$, $L_{r + n} = J_{s + n}$ and $a_{r + n} = b_{s + n}$. Hence $P(x,y,z) \ \Etail^3 \ P(x',y',z')$. 

(Case II) Suppose $a_r = b_s$. Then one must have for all $n \in \omega$, $L_{r + n} = J_{s + n}$ and $a_{r + n} = b_{s + n}$. Hence $P(x,y,z) \ \Etail^3 \ P(x',y',z')$. 

(Case III) Suppose $a_r \neq b_s$ and $L_{r + 1} \neq J_{s + 1}$. Without loss of generality, $L_{r + 1} < J_{s + 1}$, $a_r = 0$, and $b_s = 2$. This implies that for all $M \leq k < L_{r + 1}$, $x(k) = x'(k) = y(k) = y'(k) = z(k) = z'(k)$. However, $x(L_{r + 1}) \neq y(L_{r + 1})$ and $y(L_{r + 1}) = y'(L_{r + 1}) = z'(L_{r + 1}) = z(L_{r + 1})$ because $L_{r + 1} < J_{s + 1}$. Hence $a_{r + 1} = 2$. For any $k$ so that $L_{r + 1} \leq k < J_{s + 1}$, $y(k) = y'(k) = z(k) = z'(k)$. However, $y'(J_{s + 1}) \neq z'(J_{s + 1})$ hence $y(J_{s + 1}) \neq z(J_{s + 1})$ and $J_{s + 1}$ is the smallest $N > L_{r + 1}$ for which this happens. Hence $L_{r + 2} = J_{s + 1}$. Also $a_{r + 2} = b_{s + 1}$. Hence for all $n \in \omega$, $a_{(r + 2) + n} = b_{(s + 1) + n}$. This implies $P(x,y,z) \ \Etail^3 \ P(x',y',z')$. 
\end{proof}

\Begin{fact}{Etail3 reduce into Etail on A}
$\Etail^2 \leq_\borel \Etail^3 \upharpoonright A$. Hence $E_0 \equiv_\borel \Etail^3 \upharpoonright A$.
\end{fact}

\begin{proof}
Let $\Phi : \cantorspace \rightarrow A$ be defined by $\Phi(x) = x \oplus \tilde 2$, where 
$$(x \oplus y)(n) = \begin{cases}
x(k) & \quad n = 2k \\
y(k) & \quad n = 2k + 1
\end{cases}$$

Suppose $x \ \Etail^2 \ y$. Then there are some $a,b \in \omega$ so that for all $n$, $x(a + n) = y(b + n)$. For all $n \in \omega$, $\Phi(x)(2a + n) = \Phi(y)(2b + n)$. 

Suppose $\neg(x \ \Etail^2 \  y)$. Let $a,b \in \omega$. Suppose $a$ is even and $b$ is odd. Then $\Phi(x)(a + 0) \in 2$ but $\Phi(y)(b + 0) = 2$. The same argument works if $a$ is odd and $b$ is even. Suppose $a$ and $b$ are both even. Let $a = 2a'$ and $b = 2b'$. Since $\neg (x \ \Etail^2 \ y)$, there is some $k$ so that $x(a' + k) \neq y(b' + k)$. Then $\Phi(x)(a + 2k) = x(a' + k) \neq y(b' + k) = \Phi(y)(b + 2k)$. Suppose $a$ and $b$ are both odd. $a = 2a' + 1$ and $b = 2b' + 1$. Since $\neg(x \ \Etail^2 \ y)$, there is some $k$ so that $x((a' + 1) + k) \neq y((b' + 1) + k)$. Hence $\Phi(x)(a + (1 + 2k)) \neq \Phi(y)(b + (1 + 2k))$. This shows $\neg(\Phi(x) \ \Etail^3 \ \Phi(y))$. 

Since $\Etail^2 \equiv_\borel E_0$ and $\Etail^3 \upharpoonright A \leq_\borel \Etail^3 \equiv_\borel E_0$, one has $E_0 \equiv_\borel \Etail^3 \upharpoonright A$. 
\end{proof}

\Begin{theorem}{E0 does not have 3-Jonsson property}
$(\mathsf{ZF + AD})$ $\cantorspace \slash E_0$ does not have the 3-J\'onsson property. 
\end{theorem}

\begin{proof}
Let $P$ be the function from Fact \ref{function P into A is E0 to tail invariant}. Let $\bar{P} : [\cantorspace \slash E_0]^3_= \rightarrow A \slash \Etail^3$ be defined by $\bar{P}(a,b,c) = d$ if and only if 
$$(\forall x,y,z)((x \in a \wedge y \in b \wedge z \in c) \Rightarrow P(x,y,z) \in d).$$ 
By Fact \ref{function P into A is E0 to tail invariant}, $\bar{P}$ is a well defined function. Since $E_0 \equiv_\borel \Etail^3 \upharpoonright A$ by Fact \ref{Etail3 reduce into Etail on A}, let $U : A \slash \Etail^3 \rightarrow \cantorspace \slash E_0$ be a bijection (given by Fact \ref{cantor-schroder-bernstein}).

Let $F : [\cantorspace \slash E_0]^3_= \rightarrow \cantorspace \slash E_0$ be defined by $U \circ \bar{P}$. 

Let $X \subseteq \cantorspace \slash E_0$ be such that there is a bijection $B : \cantorspace \slash E_0 \rightarrow X$. By Fact \ref{comeager uniformization}, $\mathsf{AD}$ implies $B$ has a lift $B': D \rightarrow \bigcup X$, where $D \subseteq \cantorspace$ is some comeager set. Since $B$ was a bijection, $B'$ is a reduction of $E_0 \upharpoonright D$ to $E_0 \upharpoonright \bigcup X$. Using $\mathsf{AD}$, there is a comeager set $C \subseteq D$ so that $B' \upharpoonright C : C \rightarrow \bigcup X$ is a continuous reduction of $E_0 \upharpoonright C$ to $E_0 \upharpoonright \bigcup X$. There is a continuous function witnessing $E_0 \leq_{\borel} E_0 \upharpoonright C$. By composition, there is a continuous reduction witnessing $E_0 \leq_\borel E_0 \upharpoonright B'[C]$. There is an $E_0$-tree $p$ so that $[p] \subseteq B'[C] \subseteq \bigcup X$. By Fact \ref{function to 3w surjective E0 tree}, $P[[[p]]_{E_0}^3] = A$. This implies that $\bar{P}[[X]_=^3] = A \slash \Etail^3$. Since $U$ is a bijection, $U[\bar{P}[[X]^3_=]] = F[[X]^3_=] = \cantorspace \slash E_0$. $F$ witnesses $\cantorspace \slash E_0$ does not have the $3$-J\'onsson property.
\end{proof}

As mentioned earlier, since this paper is often concerned with sets without well-orderings, the J\'onsson property is defined using sets of tuples $[A]^n_=$. The usual definition of the J\'onsson property (of cardinals) involve the $n$-elements subsets of $A$, $\mathscr{P}^n(A)$. This paper calls this the classical $n$-J\'onsson property. With a slight modification, one can also obtain the failure of the classical $3$-J\'onsson property for $\cantorspace \slash E_0$. 

\Begin{definition}{equivalence relation F}
Let $S_3$ be the permutation group on $3 = \{0,1,2\}$. $S_3$ acts on ${}^\omega 3$ in the natural way: if $p \in S_3$ and $x \in {}^\omega 3$, then $(p\cdot x)(n) = p(x(n))$.

Let $F$ be an equivalence relation on ${}^\omega 3$ defined by $x \ F \ y$ if and only if $(\exists p \in S_3)(p \cdot x \ \Etail^3 \ y)$. 
\end{definition}

\Begin{fact}{F is borel bireducible to E0}
Let $A = \{x \in {}^\omega 3 : (\forall n)(x(n) \neq x(n + 1))\}$. $F \upharpoonright A \equiv_\borel E_0$. 
\end{fact}

\begin{proof}
Note that $\Etail^3 \upharpoonright A$ is hyperfinite by Fact \ref{Etail3 reduce into Etail on A}. Note that $\Etail^3 \upharpoonright A \subseteq F \upharpoonright A$ and each $F \upharpoonright A$ equivalence class is a union of at most six $\Etail^3 \upharpoonright A$ equivalence classes. By a result of Jackson, $F \upharpoonright A$ is hyperfinite. Hence, $F \upharpoonright A \leq_\borel E_0$. 

Next a reduction  $\Phi: \cantorspace \rightarrow A$ will be produced witnessing $\Etail^2 \leq_\borel F \upharpoonright A$: 
$$\Phi(x) = x(0)\hat{\ }2012102\hat{\ }x(1)\hat{\ }2012102\hat{\ }x(2)...$$

If $x \ \Etail^2 \ y$, then $\Phi(x) \ F \ \Phi(y)$.

Suppose $\Phi(x) \ F \ \Phi(y)$. This means there is some $g \in S_3$ so that $g \cdot \Phi(x) \ \Etail^3 \  \Phi(y)$. Consider what happens for each $g \in S_3$: $g$ will be presented in cycle notation.

$g = \mathrm{id}$: It is clear that $g \cdot \Phi(x) \ \Etail^3 \ \Phi(y)$ implies that $x \ \Etail^2 \ y$. 

$g = (0,1)$: Then a portion of $g \cdot \Phi(x)$ looks like
$$...\hat{\ }g(x(i))\hat{\ }2102012\hat{\ }g(x(i + 1))\hat{\ }2102012\hat{\ }g(x(i + 2))\hat{\ }...$$

$g = (0,2)$: 
$$...\hat{\ }g(x(i))\hat{\ }0210120\hat{\ }g(x(i + 1))\hat{\ }0210120\hat{\ }g(x(i + 2))\hat{\ }...$$

$g = (0,1,2)$:
$$...\hat{\ }g(x(i))\hat{\ }0120210\hat{\ }g(x(i + 1))\hat{\ }0120210\hat{\ }g(x(i + 2))\hat{\ }...$$

$g = (0,2,1)$:
$$...\hat{\ }g(x(i))\hat{\ }1201021\hat{\ }g(x(i + 1))\hat{\ }1201021\hat{\ }g(x(i + 2))\hat{\ }...$$

In all these cases, $\Phi(y)$ will contain a block of $2012102$, but $g \cdot \Phi(x)$ can not possibly contain such a block. So it is impossible that $g \cdot \Phi(x) \ \Etail^3 \ \Phi(y)$. 

$g = (1,2)$:
$$...\hat{\ }g(x(i))\hat{\ }1021201\hat{\ }g(x(i + 1))\hat{\ }1021201\hat{\ }g(x(i + 2))\hat{\ }...$$
The only way that some tail of $g \cdot \Phi(x)$ contains blocks of $2012102$ is if $x \ \Etail^2 \ \tilde 1$. This however forces $g \cdot \Phi(x) \ \Etail^3 \ \Phi(\tilde 1)$. This implies that both $x \ \Etail^2 \ \tilde 1$ and $y \ \Etail^2 \ \tilde 1$. $x \ \Etail^2 \ y$. 

This shows that $\Phi$ is a reduction of $\Etail^2$ into $F \upharpoonright A$. 

This completes the proof that $E_0 \equiv_\borel F \upharpoonright A$. 
\end{proof}

\Begin{theorem}{E0 does not have classical 3-jonsson}
$(\mathsf{ZF + AD})$ $\cantorspace \slash E_0$ does not have the classical $3$-J\'onsson property.
\end{theorem}

\begin{proof}
Let $P : [\cantorspace]^3_{E_0} \rightarrow A$ be the function from Fact \ref{function to 3w surjective E0 tree}. Fact \ref{function P into A is E0 to tail invariant} shows that $P$ is $E_0$ to $\Etail^3$ invariant. Note that if $(x_0, x_1, x_2) \in [\cantorspace]^3_{E_0}$ and $g \in S_3$, then there is some other $h \in S_3$ so that $P(x_{g(0)}, x_{g(1)}, x_{g(2)}) = h \cdot P(x_0,x_1,x_2)$. 

Define a function $\Psi : \mathscr{P}^3(\cantorspace \slash E_0) \rightarrow A \slash F$ as follows: Let $D \in \mathscr{P}^3(\cantorspace \slash E_0)$. Choose any $(x_0, x_1, x_2) \in [\cantorspace]^3_{E_0}$ so that $D = \{[x_0]_{E_0}, [x_1]_{E_0}, [x_2]_{E_0}\}$. Let $\Psi(A) = [P(x_0,x_1,x_2)]_F$.

By the above observations, $\Psi$ is a well-defined surjection onto $A \slash F$. By Fact \ref{F is borel bireducible to E0}, there is a bijection $\Gamma : A \slash F \rightarrow \cantorspace \slash E_0$. Let $\Phi = \Gamma \circ \Psi$. By an argument similar to Theorem \ref{E0 does not have 3-Jonsson property}, $\Phi$ witnesses $\cantorspace \slash E_0$ does not have the classical $3$-J\'onsson property.
\end{proof}

\section{Failure of Partition Properties of $\cantorspace \slash E_0$ in Dimension Higher Than 2}\label{Failure of Partition Properties of R/E0 in Dimension Higher than 2}

This section will use the failure of the classical $3$-J\'onsson property to show that the classical partition property in dimension three fails for $\bbR \slash E_0$. Note that for any $Y$, the failure of $\cantorspace \slash E_0 \rightarrow (\cantorspace \slash E_0)^3_Y$ implies the failure of $\cantorspace \slash E_0 \mapsto (\cantorspace)^3_Y$. 

\Begin{theorem}{failure partition R/E0 in higher dimensions}
$(\mathsf{ZF + AD})$ For any set $Y$ with at least two elements, $\cantorspace \slash E_0 \rightarrow (\cantorspace \slash E_0)^3_{Y}$ fails.

In fact, if $Y$ is a set so that there is a partition of $\cantorspace \slash E_0$ by nonempty sets indexed by elements of $Y$, then there is map $f : \mathscr{P}^3(\cantorspace \ \slash E_0) \rightarrow Y$ with the property that for all $C \subseteq \cantorspace \slash E_0$ with $C \approx \cantorspace \slash E_0$, $f[\mathscr{P}^3(C)] \approx Y$.
\end{theorem}

\begin{proof}
Let $a,b \in Y$. Partition $\cantorspace \slash E_0$ into two nonempty disjoint sets $A$ and $B$. Let $\Lambda : \cantorspace \slash E_0 \rightarrow Y$ be defined by
$$\Lambda(x) = \begin{cases}
a & \quad x \in A \\
b & \quad x \in B
\end{cases}$$
Let $\Phi$ be the classical $3$-J\'onsson function from the proof of Theorem \ref{E0 does not have classical 3-jonsson}.

Define $f : \mathscr{P}^3(\cantorspace \slash E_0) \rightarrow Y$ by $f = \Lambda \circ \Phi$.

Suppose $C \subseteq \cantorspace \slash E_0$ and $C \approx \cantorspace \slash E_0$. Suppose $a_0 \in A$ and $b_0 \in B$. Since $\Phi$ is a classical $3$-J\'onsson map, there are some $R,S \in \mathscr{P}^3(C)$ so that $\Phi(R) = a_0$ and $\Phi(S) = b_0$. Since $f(R) = a$ and $f(S) = b$. $|f[\mathscr{P}^3(C)]| = 2$. 

For the second statement, suppose $(A_y : y \in Y)$ is a partition of $\cantorspace \slash E_0$ into nonempty sets. Define $\Lambda : \cantorspace \slash E_0 \rightarrow Y$ by $\Lambda(x) = y$ if and only if $x \in A_y$. Let $f = \Lambda \circ \Phi$. This maps works by an argument like above. 
\end{proof}

Given a set $X$ and $n \in \omega$, one can define $d_X(n)$ to be the smallest element of $\omega$, if it exists, such that for every $k$ and every function $f : \mathscr{P}^n(X) \rightarrow k$, there is some $S \subseteq k$ with $|S| \leq d_X(n)$ and $A \subseteq X$ with $A \approx X$ so that $f[\mathscr{P}^n(A)] \subseteq S$. Say that $d_X(n)$ is infinite if no such integer can be found.

\cite{A-Partition-Theorem-For-Perfect-Sets} showed that assuming the appropriate sets have the Baire property, for every $n, k \in \omega$ and function $f : \mathscr{P}^n(\cantorspace) \rightarrow k$, there is an $S \subseteq k$ with $|S| \leq (n - 1)!$ and a set $A \subseteq \cantorspace$ with $A \approx \cantorspace$ so that $f[\mathscr{P}^n(A)] \subseteq S$. Hence for $n > 0$, $d_\cantorspace(n) \leq (n - 1)!$ assuming $\mathsf{AD}$. 

Under $\mathsf{AD^+}$, $d_{\cantorspace \slash E_0}(2)$ is finite and equal to $1$, but for $n \geq 3$, $d_{\cantorspace \slash E_0}(n)$ is infinite.

\section{$\hatPE^3$ Is Proper}\label{hatPE3 is proper}
Fact \ref{2-hat E0 forcing is proper} shows that $\hatPE^2$ is proper by having a very flexible fusion argument. Moreover, below any condition $(p,q) \in \hatPE^2$ and countable elementary submodels $M$, one can find a $(M,\hatPE)$-master condition $(p',q')$ so that every element of $[p'] \times_{E_0} [q']$ is a $\hatPE^2$-generic real over $M$. This fusion argument for $\hatPE^2$ is also used to prove numerous combinatorial properties in dimension $2$. The analog of most of these properties in dimension $3$ fails. No fusion with the type of property that $\hatPE^2$ has can exist for $\hatPE^3$. The natural question to ask would be whether $\hatPE^3$ is proper at all. 

This section will show that $\hatPE^3$ is proper via a fusion argument. However, one loses control of when exactly dense sets are met.

\Begin{definition}{extension hatPE3 condition by copying}
Suppose $(p,q,r) \in \hatPE^3$. Let $(u,v,z)$ be a triple of strings in $\finBinarySequence$ of the same length $n + 1$ so that $\{0,1\} = \{u(n),v(n),z(n)\}$. Suppose $(p',q',r') \leq_{\hatPE^3} (\Xi(p,u),\Xi(q,v),\Xi(r,z))$. Let $\fprune^{(u,v,z)}_{(p',q',r')}(p,q,r)$ be the unique condition $(a,b,c) \leq_{\hatPE^3}^{n + 1} (p,q,r)$ so that $\Xi(a,u) = p'$, $\Xi(b,v) = q'$ and $\Xi(c,z) = r'$. 

In the above, the relation $\leq^n_{\hatPE^3}$ is defined as coordinate-wise $\leq_\PE^n$. $\fprune^{(u,v,z)}_{(p',q',r')}(p,q,r)$ has an explicit definition that is obtained by copying $p'$, $q'$ and $r'$ below the appriopriate part of $(p,q,r)$ like in Definition \ref{hat PE0 pruning notation}.
\end{definition}

\Begin{fact}{forcing hatPE3 is proper}
(With Zapletal) $\hatPE^3$ is a proper forcing.
\end{fact}

\begin{proof}
Let $(p,q,r) \in \hatPE^3$. Let $\Xi$ be some large regular cardinal. Let $M \prec V_\Xi$ be a countable elementary substructure containing $(p,q,r)$. Let $(D_n : n \in \omega)$ enumerate all the dense open subsets of $\hatPE^3$ that belong to $M$. One may assume that $D_{n + 1} \subseteq D_n$ for all $n \in \omega$. 

If there exists some condition $(p',q',r') \leq^0_{\hatPE^3} (p,q,r)$ with $(p',q',r') \in D_0$, then by elementarity there is such a condition in $M$. Let $(p_0,q_0,r_0)$ be such a condition in $M$. Otherwise, let $(p_0,q_0,r_0) = (p,q,r)$. 

Suppose $(p_n,q_n,r_n)$ have been defined. Let $\{(u_i,v_i,z_i) : i < K\}$ enumerate all the strings in $\finBinarySequence$ with length $n + 1$ so that $\{u(n),v(n),z(n)\} = \{0,1\}$.  

Let $(a_{-1},b_{-1},c_{-1}) = (p_n,q_n,r_n)$. For $i$ with $-1 \leq i < K - 1$, suppose $(a_i,b_i,c_i)$ has been defined. 

Let $(a_{i + 1}^{-1}, b_{i + 1}^{-1}, c_{i + 1}^{-1}) = (a_i,b_i,c_i)$. Suppose for $j$ with $-1 \leq j < n + 1$, $(a_{i + 1}^{j}, b_{i + 1}^j, c_{i + 1}^j)$ has been defined. If there exists some condition below $(\Xi(a_{i + 1}^j, u_{i + 1}), \Xi(b_{i + 1}^j, v_{i + 1}), \Xi(c_{i + 1}^j, z_{i + 1}))$ that belongs to $D_{j + 1}$, then choose, by elementarity, such a condition $(a',b',c') \in M \cap D_{j + 1}$. Let $(a_{i + 1}^{j + 1}, b_{i + 1}^{j + 1}, c_{i + 1}^{j + 1}) = \fprune^{(u_{i + 1},v_{i + 1},z_{i + 1})}_{(a',b',c')}(a_{i + 1}^j, b_{i + 1}^j, c_{i + 1}^j)$. If no such condition exists, then let $(a_{i + 1}^{j + 1}, b_{i + 1}^{j + 1}, c_{i + 1}^{j + 1}) = (a_{i + 1}^j, b_{i + 1}^j, c_{i + 1}^j)$. 

Let $(a_{i + 1}, b_{i + 1}, c_{i + 1}) = (a_{i + 1}^{n + 1}, b_{i + 1}^{n + 1}, c_{i + 1}^{n + 1})$. Let $(p_{n + 1}, q_{n + 1}, r_{n + 1}) = (a_{K - 1}, b_{K - 1}, c_{K - 1})$. Note that $(p_{n + 1}, q_{n + 1}, r_{n + 1}) \leq^{n + 1}_{\hatPE^3} (p_n,q_n,r_n)$. 

$\langle (p_n,q_n,r_n) : n \in \omega \rangle$ forms a fusion sequence in $\hatPE^3$. Let $(p_\omega,q_\omega,r_\omega)$ be the fusion of this fusion sequence. The claim is that this is a $(M,\hatPE^3)$-master condition below $(p,q,r)$.

It needs to be shown for each $n$ that $(p_\omega,q_\omega,r_\omega) \forces_{\hatPE^3} \check M \cap \check D_n \cap \dot G \neq \emptyset$. Let $G$ be any $\hatPE^3$-generic over $M$ containing $(p_\omega,q_\omega,r_\omega)$. There is some $(p',q',r') \in G \cap D_n$. Since $G$ is a filter, there is some $(p'',q'',r'') \leq_{\hatPE^3} (p_\omega,q_\omega,r_\omega)$ so that $(p'',q'',r'') \in G \cap D_n$. By genericity, one may assume there is some $m > n$ and some $(u,v,z) \in {}^m 2$ so that $(p'',q'',r'') \leq_{\hatPE^3} (\Xi(p_\omega, u),\Xi(q_\omega,v), \Xi(r_\omega,z))$. During the construction while producing $(p_m,q_m,r_m)$, the strings $(u,v,z) = (u_i,v_i,z_i)$ for some $i$ in the chosen enumeration of strings. Note that $(p'',q'',r'') \leq_{\hatPE^3} (\Xi(a_i^{n - 1}, u_i), \Xi(b_i^{n - 1}, v_i), \Xi(c_{i}^{n - 1}, z_i))$ and $(p'',q'',r'') \in D_n$. At this stage, one would have chosen some $(a',b',c') \in M \cap D_n$ below $(\Xi(a_i^{n - 1}, u_i), \Xi(b_i^{n - 1}, v_i), \Xi(c_{i}^{n - 1}, z_i))$ and set $(a_i^n,b_i^n,c_i^n) = \fprune_{(a',b',c')}^{(u_i,v_i,z_i)}(a_i^{n - 1},b_i^{n - 1},c_i^{n - 1})$. Note that $(p'',q'',r'') \leq_{\hatPE^3} (\Xi(p_\omega,u), \Xi(q_\omega,v), \Xi(r_\omega, z)) \leq_{\hatPE^3} (a',b',c')$. Since $G$ is a filter and $(p'',q'',r'') \in G$, $(a',b',c') \in G \cap M \cap D_n$. This shows that $\hatPE^3$ is a proper forcing. 
\end{proof}

In the proof, one extends a portion of the three trees to get into a dense set $D$ only if it was possible and otherwise ignored $D$. Because of this, one cannot prove that $[p_\omega] \times_{E_0} \times [q_\omega] \times_{E_0} [r_\omega]$ consists entirely of reals which are $\hatPE^3$-generic over $M$. 

\section{$E_1$ Does Not Have the 2-Mycielski Property}\label{E1 Does Not Have the 2-Mycielski Property}

This section will give an example to show $E_1$ does not have the $2$-Mycielski property. The notation of Definition \ref{basic open neighborhood pcantorspace} will be used in the following.

As in earlier sections, an understanding of the structure theorem of $E_1$-big $\analytic$ sets is essential:

\Begin{definition}{E1 definition}
$E_1$ is the equivalence relation on $\pcantorspace$ defined by $x \ E_1 \ y$ if and only if there exists a $k$ so that for all $n \geq k$, $x(n) = y(n)$. 
\end{definition}

\Begin{definition}{keeping homomorphism}
\cite{Canonical-Ramsey-Theory-on-Polish-Spaces} Let $s$ be an infinite subset of $\omega$. Let $\pi_s : \omega \rightarrow s$ be the unique increasing enumeration of $s$. A homeomorphism $\Phi : \pcantorspace \rightarrow \pcantorspace$ is an $s$-keeping homeomorphism if and only if the following hold:

1. For all $n \in \omega$, if $x(n) \neq y(n)$, then $\Phi(x)(\pi_s(m)) \neq \Phi(y)(\pi_s(m))$ for all $m \leq n$. 

2. For all $n \in \omega$, if for all $m > n$, $x(m) = y(m)$, then for all $k > \pi_s(n)$, $\Phi(x)(k) = \Phi(y)(k)$. 
\end{definition}

\Begin{fact}{E1 bireduciblity characterization}
Let $B \subseteq \pcantorspace$ be $\analytic$. $E_1 \upharpoonright B \equiv_\borel E_1$ if and only if there is some infinite $s \subseteq \omega$ and $s$-keeping homeomorphism $\Phi$ so that $\Phi[\pcantorspace] \subseteq B$. 
\end{fact}

\begin{proof}
This result is implicit in \cite{The-Classification-of-Hypersmooth-Borel-Equivalence-Relations}. See \cite{Canonical-Ramsey-Theory-on-Polish-Spaces}, Section 7.2.1. 
\end{proof}

\Begin{theorem}{E1 does not have Mycielski property}
Let $D = \{(x,y) \in {}^2(\pcantorspace) : (\exists n)(x(n)(0) \neq y(n)(0))\}$. $D$ is dense open and for all $\borel$ $B$ such that $E_1 \upharpoonright B \equiv_\borel E_1$, $[B]^2_{E_1} \not\subseteq D$. 

$E_1$ does not have the 2-Mycielski property.
\end{theorem}

\begin{proof}
Suppose $(x,y) \in D$. There is some $n \in \omega$ so that $x(n)(0) \neq y(n)(0)$. Let $\sigma,\tau : (n + 1) \rightarrow {}^12$ be defined by $\sigma(k) = x(k) \upharpoonright 1$ and $\tau(k) = y(k) \upharpoonright 1$. Then $(x,y) \in N_{\sigma,\tau} \subseteq D$. $D$ is open.

Let $\sigma, \tau : m \rightarrow \finBinarySequence$. Define $\sigma', \tau' : (m + 1) \rightarrow \finBinarySequence$ by 
$$\sigma'(k) = \begin{cases}
\sigma(k) & \quad k < m \\
\langle 0 \rangle & \quad k = m 
\end{cases} 
\ \ \ \text{ and } \ \ \ 
\tau'(k) = \begin{cases}
\tau(k) & \quad k < m \\
\langle 1 \rangle & \quad k = m
\end{cases}.$$ 
$N_{\sigma',\tau'} \subseteq N_{\sigma,\tau}$ and $N_{\sigma',\tau'} \subseteq D$. $D$ is dense open.

Let $s \subseteq \omega$ be infinite. Let $\pi_s : \omega \rightarrow s$ be the unique increasing enumeration of $s$. Let $\Phi : \pcantorspace \rightarrow \pcantorspace$ be an $s$-keeping homeomorphism. 

Let $\delta_{n} : (\pi_s(n) + 1) \rightarrow {}^12$ be defined by $\delta_{n}(k) = \Phi(\bar{0})(k) \upharpoonright 1$. A strictly increasing sequence $\langle m_n : n \in \omega\rangle$ of natural numbers and functions $\sigma_n : m_n \rightarrow {}^{m_n}2$ satisfying the following for all $n \in \omega$ will be defined:

\noindent 1. $\sigma_n(k) \subseteq \sigma_{n + 1}(k)$ for each $k < m_n$. 

\noindent 2. $\Phi[N_{\sigma_n}] \subseteq N_{\delta_n}$.

\noindent 3. There exists a $j < m_n$ such that $\sigma_n(n)(j) = 1$ and for all $k > n$ and $i < m_n$, $\sigma_n(k)(i) = 0$.

Let $m_{-1} = 0$ and $\sigma_{-1} = \delta_{-1} = \emptyset$. 

Suppose $m_n$ and $\sigma_n$ have been defined and satisfy conditions 2 and 3 if $n \geq 0$. Define $y \in N_{\sigma_n}$ by $y(i)(j) = 0$ if $(i,j) \notin m_n \times m_n$. Then $\Phi(y) \in N_{\delta_n}$. Since $y(k) = \bar{0}(k)$ for all $k > n$ and $\Phi$ is an $s$-keeping homeomorphism, $\Phi(y)(k) = \Phi(\bar{0})(k)$ for all $k > \pi_s(n)$. Thus $\Phi(y) \in N_{\delta_{n + 1}}$. By continuity of $\Phi$, there is some $M \geq m_n$ so that if $\tau : M \rightarrow {}^M2$ is defined by $\tau(i) = y(i) \upharpoonright M$, then $\Phi(N_\tau) \subseteq N_{\delta_{n + 1}}$. Let $m_{n + 1} = M + 1$ and define
$$\sigma_{n + 1}(i)(j) = \begin{cases}
1 & \quad i= n + 1 \wedge j = M \\
y(i)(j) & \quad \text{otherwise}
\end{cases}.$$

$m_{n + 1}$ and $\sigma_{n + 1}$ satisfy conditions 1, 2, and 3.

Let $x \in \pcantorspace$ be so that $\{x\} = \bigcap_{n \in \omega} N_{\sigma_n}$. $\neg(\bar{0} \ E_1 \ x)$ since for all $n$, there exists a  $j$ so that $x(n)(j) = 1$ by condition 3. However since $\Phi(x) \in N_{\delta_n}$ for all $n$, $(\Phi(\bar{0}), \Phi(x)) \notin D$. From Definition \ref{keeping homomorphism}, $\Phi$ is a $E_1$ reduction so $\neg(\Phi(\bar{0}) \ E_1 \ \Phi(x))$.  Hence $[\Phi(\pcantorspace)]^2_{E_1} \not\subseteq D$. 

So it has been shown that for all infinite $s \subseteq \omega$ and $s$-keeping homeomorphisms $\Phi$, $[\Phi[\pcantorspace]]^2_{E_1} \not\subseteq D$. By Fact \ref{E1 bireduciblity characterization}, every $\borel$ set $B$ so that $E_1 \upharpoonright B \equiv_\borel E_1$ contains $\Phi[\pcantorspace]$ for some $s$ and some $s$-keeping homeomorphism $\Phi$. Therefore, $[B]^2_{E_1} \not\subseteq D$ for all such $B$. $E_1$ does not have the 2-Mycielski property.
\end{proof}

\section{The Structure of $E_2$}\label{The Structure of E2}

This section will give a proof of a result about the structure of $E_2$-big sets necessary for analyzing the weak-Mycielski property for $E_2$. The proof is similar to but a bit a more technical than the argument of \cite{Borel-Equivalence-Relations} Theorem 15.4.1. Some of the notation and terminology come from \cite{Borel-Equivalence-Relations}.

\Begin{definition}{summation distance notation}
Suppose $x,y \in \cantorspace$. Define
$$\delta(x,y) = \sum_{k \in x \triangle y} \frac{1}{k + 1}.$$
Suppose $m < n \leq \omega$. Suppose $x,y \in {}^N2$ where $n \leq N \leq \omega$. Define
$$\delta_m^n(x,y) = \sum \left\{\frac{1}{k + 1} : (m \leq k < n) \wedge (k \in x \triangle y)\right\}.$$
Let $A, B \subseteq \cantorspace$. Let $m < n \leq \omega$. Let $\epsilon > 0$. Define $\delta_m^n(A,B) < \epsilon$ if and only if for all $x \in A$, there is some $y \in B$ so that $\delta_m^n(x,y) < \epsilon$ and for all $y \in B$, there exists some $x \in A$ so that $\delta_m^n(x,y) < \epsilon$. 
\end{definition}

\Begin{definition}{E2 equivalence relation}
$E_2$ is the equivalence relation on $\cantorspace$ defined $x \ E_2 \ y$ if and only if $\delta(x,y) < \infty$. 
\end{definition}

\Begin{lemma}{delta is a metric}
Let $m < n \leq \omega$. Fix $N$ so that $n \leq N \leq \omega$. If $n < \omega$, then $\delta_m^n$ is a pseudo-metric on ${}^N2$. If $n = \omega$, then $\delta_m^n$ is a pseudo-metric on any $E_2$ equivalence class.
\end{lemma}

\Begin{lemma}{E2 splitting lemma 1}
Let $m,p \in \omega$. Let $q \in \bbQ^+$. Let $(\hat X_i : i < p)$ be a sequence of $\lanalytic(z)$ subsets of $\cantorspace$. Let $(x_i : i < p)$ be a sequence in $\cantorspace$ with the property that $x_i \in \hat X_i$ and $\delta_m^\omega(x_0, x_i) < q$. Then there exists a sequence $(X_i : i < p)$ of $\lanalytic(z)$ sets with $x_i \in X_i$ and $\delta_m^\omega(X_0, X_i) < q$. 
\end{lemma}

\begin{proof}
Let 
$$X_0 = \left\{x \in \hat X_0 : (\exists z_1, ..., z_{p - 1})\left(\bigwedge_{1 \leq i < p} z_i \in \hat X_i \wedge \delta_m^\omega(x,z_i) < q\right)\right\}.$$
For $1 \leq i < p$, define
$$X_i = \left\{x \in \hat X_i : (\exists z)(x \in X_0 \wedge \delta_m^\omega(x,z) < q)\right\}.$$
\end{proof}

\Begin{lemma}{E2 splitting lemma 2}
Let $m,p \in \omega$. Let $q \in \bbQ^+$. Let $(\hat X_i : i < p)$ be a sequence of $\lanalytic(z)$ sets with $\delta_m^\omega(\hat X_0, \hat X_i) < q$ for all $i < p$. Let $j < p$. Suppose $A \subseteq \hat X_j$ is a $\lanalytic(z)$ set. Then there exists a sequence $(X_i : i < p)$ of $\lanalytic(z)$ sets with the property that for all $i < p$, $X_i \subseteq \hat{X}_i$, $X_j = A$, and $\delta_m^\omega(X_0, X_i) < q$.
\end{lemma}

\begin{proof}
Let
$$X_0 = \{x \in \hat X_0 : (\exists z)(z \in A \wedge \delta_m^\omega(x,z) < q)\}.$$
For all $i < p$ and $i \neq j$, let
$$X_i = \{x \in \hat X_i : (\exists z)(z \in X_0 \wedge \delta_m^\omega(x,z) < q)\}$$
Let $X_j = A$. 
\end{proof}

\Begin{lemma}{E2 splitting lemma 3}
Let $m,p \in \omega$. Let $q \in \bbQ^+$. Let $(\hat X_i : i < p)$ be a sequence of $\lanalytic(z)$ sets with $\delta_m^\omega(\hat X_0, \hat X_i) < q$ for all $i < p$. Let $D$ be a dense open subset of $\bbP_z$. Then there exists a sequence $(X_i : i < p)$ of $\lanalytic(z)$ sets with $X_i \subseteq \hat X_i$, $X_i \in D$, and $\delta_m^\omega(X_0,X_i) < q$  for all $i < p$.
\end{lemma}

\begin{proof}
Let $Y_i^{-1} = \hat X_i$. 

One seeks to define $\lanalytic(z)$ sets $Y_i^{j}$ for all $-1 \leq j < p$ with the property that if $-1 \leq j < p -1$, then $Y_i^{j + 1} \subseteq Y_i^{j}$, and for any $-1 \leq j < p$ and $0 \leq i < p$, $Y_j^j \in D$ and $\delta_m^\omega(Y_0^j, Y_i^j) < q$. 

Suppose $Y_i^j$ has been defined with the desired properties for $j < p - 1$ and all $i < p$. Since $D$ is dense open in $\bbP_z$, pick some $A \subseteq Y_{j + 1}^{j}$ so that $A \in D$. Use Lemma \ref{E2 splitting lemma 2} with $\{Y_i^{j} : i < p\}$ and $A \subseteq Y_i^j$ to obtain $\{Y_i^{j + 1} : i < p\}$ with the desired properties. 

Let $X_i = Y_i^{p - 1}$. 
\end{proof}

In the previous three lemmas, the first set was distinguished. In the following argument, sets may be indexed by strings and so in applications of the three lemmas, one will need to indicate what this distinguished set is. 

\Begin{lemma}{E2 splitting lemma 4}
Let $z \in \cantorspace$. Let $k,m,p \in \omega$. Let $r,v \in \bbQ^+$. Let $(B_s^i : s \in {}^k2 \wedge i < p)$ be a sequence of $\lanalytic(z)$ sets. Let $(b_s^i : s \in {}^k 2 \wedge i < p)$ be a sequence in $\cantorspace$ with $b_s^i \in B_s^i$. Suppose for all $i < p$, $\delta_m^\omega(b_{0^k}^0, b_{0^k}^i) < r$. Suppose for each $i < p$ and for all $s \in {}^k2$, $\delta_m^\omega(b_{0^k}^i, b_s^i) < v$. Let $D$ be a dense open subset of $\bbP_z$. 

Then there is a sequence $(C_s^i : s \in {}^k2 \wedge i < p\}$ of $\lanalytic(z)$ sets so that for all $i < p$ and $s \in {}^k2$,  $\delta_m^\omega(C^0_{0^k}, C^i_{0^k}) < r$, $\delta_m^\omega(C_{0^k}^i, C_s^i) < v$, and $C_s^i \in D$. 
\end{lemma}

\begin{proof}
For each $i < p$, apply Lemma \ref{E2 splitting lemma 1} to $\{B_s^i : s \in {}^k2\}$ and $\{b_s^i : s \in {}^k2\}$ using $0^k$ as the distinguished index to obtain a sequence of $\lanalytic(z)$ sets $\{E_s^i : s \in {}^k2\}$ with the property that $E_s^i \subseteq B_s^i$, $b_s^i \in E_s^i$, and $\delta_m^\omega(E^i_{0^k}, E_s^i) < v$. 

Now apply Lemma \ref{E2 splitting lemma 1} to $\{E_{0^k}^i : i < p\}$ and $\{b_{0^k}^i : i < p\}$ with $0$ as the distinguished index to obtain $\lanalytic(z)$ sets $A_i \subseteq E_{0^k}^i$ so that $\delta_m^\omega(A_0,A_i) < r$. 

For each $i < p$, apply Lemma \ref{E2 splitting lemma 2} to $\{E_s^i : s \in {}^k2\}$ with $0^k$ as the distinguished index and $A_i \subseteq E_{0^k}^i$ to obtain $\lanalytic(z)$ sets $G_s^{i,-1}$ with the properties that $G_{0^k}^{i,-1} = A_i$ and $\delta_m^\omega(G_{0^k}^{i,-1}, G_s^{i,-1}) < v$. 

Note that since $G_{0^k}^{i,-1} = A_i$, the sequence $\{G_s^{i,-1} : i < p \wedge s \in {}^k2\}$ has the property that for all $i < p$, $\delta_m^\omega(G_{0^k}^{0, -1}, G_{0^k}^{i,-1}) < r$ and for each $i < p$ and $s \in {}^k2$, $\delta_m^\omega(G_{0^k}^{i,-1}, G_s^{i,-1}) < v$. 

One wants to create $G_s^{i,j}$ for $-1 \leq j < p$, $i < p$, and $s \in {}^k2$ so that 

\noindent (i) For each $i \in \omega$, $s \in {}^k2$, and $-1 \leq l \leq j < p$, $G^{i,j}_s \subseteq G^{i,l}_s$.  

\noindent (ii) For all $i < p$, $\delta_m^\omega(G_{0^k}^{0,j}, G_{0^k}^{i,j}) < r$.

\noindent (iii) For each $i < p$ and $s \in {}^k2$, $\delta_m^\omega(G_{0^k}^{i,j}, G_s^{i,j}) < v$. 

\noindent (iv) If $j \geq 1$ and $0 \leq l \leq j$, then $G_s^{l,j} \in D$. 

This already holds for $j = -1$. Suppose the construction worked up to stage $j < p - 1$ producing objects with the above properties. Apply Lemma \ref{E2 splitting lemma 3} to $\{G_s^{j + 1, j} : s \in {}^k2\}$ to get sets $\{G_s^{j + 1, j + 1} : s \in {}^k2\}$ each in $D$ and $\delta_m^\omega(G^{j + 1, j + 1}_{0^k}, G^{j + 1, j + 1}_s) < v$. 

Next apply Lemma \ref{E2 splitting lemma 2} to $\{G_{0^k}^{i,j} : i < p\}$ with $0$ as the distinguished index and $G^{j+1, j+1}_{0^k} \subseteq G^{j + 1,j}_{0^k}$ to obtain sets $G^{i,j + 1}_{0^k} \subseteq G^{i,j}_{0^k}$ with $\delta_m^\omega(G_{0^k}^{0,j + 1}, G^{i,j + 1}_{0^k}) < r$. (Note it is acceptable to use the notation $G^{j + 1, j + 1}_{0^k}$ since it is the same set as before by the statement of Lemma \ref{E2 splitting lemma 2}.) 

For each $i \neq j + 1$, apply Lemma \ref{E2 splitting lemma 2} on $\{G_s^{i,j} : s \in {}^k 2\}$ with $0^k$ as the distinguished index and $G_{0^k}^{i,j + 1} \subseteq G^{i,j}_{0^k}$ to obtain sets $\{G^{i,j + 1}_s : s \in {}^k2\}$ with the property that $\delta_m^\omega(G_{0^k}^{i,j+1}, G_{s}^{i,j + 1}) < v$. This completes the construction at stage $j + 1$. 

Let $C^i_{s} = G^{i,p-1}_s$. 
\end{proof}

\Begin{definition}{grainy sets}
(\cite{Borel-Equivalence-Relations} Definition 15.2.2) Let $q > 0$ be a rational number. Let $A \subseteq \cantorspace$. Let $a \in A$. The $q$-galaxy of $a$ in $A$, denoted $\Gal_A^q(a)$, is the set of all $b \in A$ so that there exists $a_0, ..., a_l \in A$ with $a = a_0$, $b = a_l$, and $\delta(a_i, a_{i + 1}) < q$ for all $0 \leq i < l - 1$. 

$A \subseteq \cantorspace$ is $q$-grainy if and only if for all $a \in A$ and $b \in \Gal_A^q(a)$, $\delta(a,b) < 1$. $A$ is grainy if and only if $A$ is $q$-grainy for some positive rational number $q$. 
\end{definition}

\Begin{fact}{covering analytic grainy by borel grainy}
Let $z \in \cantorspace$ and rational $q > 0$. If $A$ is a $\lanalytic(z)$ $q$-grainy set, then there is some $B \supseteq A$ which is $\lborel(z)$ $q$-grainy.
\end{fact}

\begin{proof}
(See \cite{Borel-Equivalence-Relations}, Claim 15.2.4 for a more constructive proof.)

Let $U \subseteq \omega \times \cantorspace$ be a universal $\lanalytic(z)$ set. The relation in variables $e$, $a$, and $b$ expressing $b \in \Gal_{U^e}^q(a)$ is $\lanalytic(z)$. 

Let $\mathcal{A}$ be the collection of all $\lanalytic(z)$ $q$-grainy subsets of $\cantorspace$. 
$$\{e : U^e \in \mathcal{A}\} = \{e : (\forall a)(\forall b)(b \in \Gal_{U^e}^q(a) \Rightarrow \delta(a,b) < 1)\}$$
This shows that $\mathcal{A}$ is a collection of $\lanalytic(z)$ sets which is $\lcoanalytic(z)$ in the code. By $\lanalytic(z)$ reflection, every $\lanalytic(z)$ $q$-grainy set is contained inside of a $\lborel(z)$ $q$-grainy set.
\end{proof}

\Begin{definition}{E2 dichotomy complex sets}
Let $z \in \cantorspace$. Let $S_z$ be the union of all $\lborel(z)$ grainy sets. Let $H_z = \cantorspace \setminus S_z$. 
\end{definition}

\Begin{fact}{E2 S and H properties}
Let $z \in \cantorspace$. $S_z$ is $\lcoanalytic(z)$. Hence $H_z$ is $\lanalytic(z)$.

Every nonempty $\lanalytic(z)$ subset of $H_z$ is not grainy.
\end{fact}

\begin{proof}
Similar to Fact \ref{property of E0 S and H sets}.
\end{proof}

\Begin{theorem}{structure through E2 sets satisfying effective condition}
Let $z \in \cantorspace$. Let $p \in \omega$. Suppose $(X_i : i < p)$ is a collection of $\lanalytic(z)$ subsets of $\cantorspace$ with the property that $\bigcap_{i < p} [X_i \cap H_z]_{E_2} \neq \emptyset$. Then there is a strictly increasing sequence $(m_k : k \in \omega)$ and functions $g^i : \finBinarySequence \rightarrow \finBinarySequence$ for each $i < p$ with the following properties:

1. If $|s| = k$, then for all $i < p$, $g^i(s) \in {}^{m_k}2$. 

2. If $s \subseteq t$, then for all $i < p$, $g^i(s) \subseteq g^i(t)$. 

3. If $|s| = |t| = k > 0$ and $s(k - 1) = t(k - 1)$, then $\delta_{m_{k - 1}}^{m_{k}}(g^i(s),g^i(t)) < 2^{-(k + 1)}$, 

4. If $|s| = |t| = k > 0$ and $s(k - 1) \neq t(k - 1)$, then $|\delta_{m_{k - 1}}^{m_k}(g^i(s), g^i(t)) - \frac{1}{k}| < 2^{-(k + 1)}$. 

5. For $i,j < p$ and $s \in {}^k2$ with $k > 0$, $\delta_{m_{k - 1}}^{m_k}(g^i(s),g^j(s)) < 2^{-(k + 1)}$. 

6. Define $\Phi^i(x) = \bigcup_{n \in \omega} g^i(x \upharpoonright n)$. $\Phi^i : \cantorspace \rightarrow X_i \cap H_z$ is a reduction witnessing $E_2 \leq_\borel E_2 \upharpoonright X_i \cap H_z$. 

Moreover, for $i, j < p$, $[\Phi^i[\cantorspace]]_{E_2} = [\Phi^j[\cantorspace]]_{E_2}$. 
\end{theorem}

\begin{proof}
During the construction, one will seek to create

\noindent (i) a strictly increasing sequence $(m_k : k \in \omega)$,

\noindent (ii) for each $i < p$ and $s \in \finBinarySequence$, $\lanalytic(z)$ sets $A^i_s$,

\noindent (iii) and for each $i < p$, $g^i(s) \in \finBinarySequence$. 

These objects will satisfy the following properties:

\noindent (I) If $|s| = k$, then $|g^i(s)| = m_k$. $s \subseteq t$ implies $g^i(s) \subseteq g^i(t)$. 

\noindent (II) $\emptyset \neq A_s^i \subseteq X^i \cap H_z \cap N_{g^i(s)}$. $s \subseteq t$ implies $A_t^i \subseteq A_s^i$. 

\noindent (III) If $k > 0$, $|s| = k$, then $\delta_{m_k}^\omega(A^i_{0^k},A^i_s) \leq 2^{-(k + 4)}$, where $0^k : k \rightarrow 2$ is the constant $0$ function. 

\noindent (IV) If $k > 0$, $|s| = |t| = k$, and $s(k - 1) = t(k - 1)$, then $\delta_{m_{k - 1}}^{m_k}(g^i(s), g^i(t)) \leq 2^{-(k + 1)}$. 

\noindent (V) If $k > 0$, $|s| = |t| = k$, $s(k - 1) \neq t(k - 1)$, then $|\delta_{m_{k - 1}}^{m_k}(g^i(s),g^i(t)) - \frac{1}{k}| < 2^{-(k + 1)}$.

\noindent (VI) If $|s| = k$ and $i < p$, $\delta_{m_{k}}^\omega(A_s^0, A_s^i) < 2^{-(k + 6)}$.

\noindent (VII) If $|s| = k$ and $k > 0$, then $\delta_{m_{k - 1}}^{m_k}(g^i(s),g^j(s)) < 2^{-(k + 1)}$. 

\noindent (VIII) Let $\mathcal{D} = (D_n : n \in \omega)$ be the countable collection of dense open subsets of $\bbP_z$ from Fact \ref{gandy-harrington dense open nonempty}. For all $s \in \finBinarySequence$ and $i < p$, $A_s^i \in D_{|s|}$. 

Suppose these objects having the above properties could be constructed. It only remains to verify 6: $\{\Phi^i(x)\} = \bigcap_{k \in \omega} A_{x \upharpoonright k}$ by (II) and (VIII). Hence $\Phi^i$ maps into $X_i \cap H_z$ by (II). Note that $\delta(\Phi^i(x), \Phi^j(y)) = \lim_{k \rightarrow \infty} \delta_0^{m_k}(g^i(x \upharpoonright k), g^j(y \upharpoonright k))$. Also using (IV) and (V), for any $k$
$$|\delta_0^{m_k}(g^i(x \upharpoonright k), g^i(y \upharpoonright k)) - \delta_0^k(x\upharpoonright k, y \upharpoonright k)| < \sum_{j < k} 2^{-(j + 1)} < 1$$ 
Hence $\Phi^i(x) \ E_2 \ \Phi^i(y) \Leftrightarrow \delta(\Phi^i(x),\Phi^i(y)) < \infty \Leftrightarrow \delta(x,y) < \infty \Leftrightarrow x \ E_2 \ y$. This shows each $\Phi^i$ witnesses $E_2 \leq_\borel E_2 \upharpoonright X_i \cap H_z$. Using (VII), for each $i,j < p$ and $x \in \cantorspace$,
$$\delta_{m_0}^{m_k}(g^i(x\upharpoonright k), g^j(x\upharpoonright j)) < \sum_{j < k} 2^{-(j + 1)} < 1$$
Hence $\Phi^i(x) \ E_2 \ \Phi^j(x)$. Hence $[\Phi^i[\cantorspace]]_{E_2} = [\Phi^j[\cantorspace]]_{E_2}$. 

Next the construction: Since $\bigcap_{i < p} [X_i \cap H_z]_{E_2} \neq \emptyset$, let $(b_\emptyset^i : i < p)$ be such that for all $i,j < p$, $b_\emptyset^i \ E_2 \ b_\emptyset^j$ and $b_\emptyset^i \in X_i \cap H_z$. Therefore, choose $m_0 \in \omega$ so that for all $i < p$, $\delta_{m_0}^\omega(b^0_\emptyset, b^i_\emptyset) < 2^{-6}$. For each $i < p$, let $g^i(\emptyset) = b^i_\emptyset \upharpoonright m_0$. 

Let $B_\emptyset^i = X_i \cap H_z \cap N_{g^i(\emptyset)}$. Apply Lemma \ref{E2 splitting lemma 4} to $\{B_\emptyset^i : i \in p\}$, $\{b_\emptyset^i : i < p\}$, and the dense open (in $\bbP_z$) set $D_0$ (where $r = 2^{-6}$ and $v = 1$) to obtain sets $A^i_\emptyset$ with the desired properties.

Suppose the objects from stage $k$ have been constructed with the desired properties. 

As $A^0_{0^k} \subseteq H_z$, Fact \ref{E2 S and H properties} implies that $A^0_{0^k}$ is not $2^{-(k + 5)}$-grainy. Hence there is a sequence $a_0, ..., a_M$ of points in $A_{0^k}^0$ so that for each $0 \leq j < M - 1$, $\delta(a_j, a_{j + 1}) < 2^{-(k + 5)}$ but $\delta(a_0, a_M) > 1$. Hence there is some $I$ so that $\delta(a_0, a_I) > \frac{1}{k + 1}$ and $\delta(a_0,a_I) - \frac{1}{k + 1} < 2^{-(k + 5)}$. 

Let $b^0_{0^k\hat{\ }0} = a_0$ and $b^0_{0^k\hat{\ }1} = a_I$. Since $\delta_{m_k}^\omega(A_s^0,A_s^i) < 2^{-(k + 6)}$ by (VI), find $b_{0^{k + 1}}^i, b_{0^{k}\hat{\ }1}^i \in A_{0^k}^i$ so that $\delta(b^0_{0^{k + 1}}, b^i_{0^{k + 1}}) < 2^{-(k + 6)}$ and $\delta(b_{0^k\hat{\ }1}^0, b_{0^k\hat{\ }1}^i) < 2^{-(k + 6)}$. 

The claim is that for all $i < p$, $|\delta_{m_k}^\omega (b_{0^{k + 1}}^i, b_{0^{k}\hat{\ }1}^i) - \frac{1}{k + 1}| < 2^{-(k + 4)}$: To see this, 
$$|\delta_{m_k}^\omega (b_{0^{k + 1}}^i, b_{0^{k}\hat{\ }1}^i) - \frac{1}{k + 1}|$$
$$\leq |\delta_{m_k}^\omega(b^i_{0^{k + 1}}, b^i_{0^k\hat{\ }1}) - \delta_{m_k}^\omega(b^0_{0^{k + 1}}, b^0_{0^{k}\hat{\ }1})| + |\delta_{m_k}^\omega(b^0_{0^{k + 1}}, b^0_{0^{k}\hat{\ }1}) - \frac{1}{k + 1}|$$
$$< |\delta_{m_k}^\omega(b^i_{0^{k + 1}}, b^i_{0^k\hat{\ }1}) - \delta_{m_k}^\omega(b^0_{0^{k + 1}}, b^i_{0^{k}\hat{\ }1})| + |\delta_{m_k}^\omega(b^0_{0^{k + 1}}, b^i_{0^{k}\hat{\ }1}) - \delta_{m_k}^\omega(b^0_{0^{k + 1}}, b^0_{0^{k}\hat{\ }1})| + 2^{-(k + 5)}$$
$$\leq \delta_{m_k}^\omega(b^i_{0^{k + 1}}, b^0_{0^{k + 1}}) + \delta_{m_k}^\omega(b^i_{0^{k} \hat{\ }1}, b^0_{0^{k}\hat{\ }1}) + 2^{-(k + 5)}$$ 
$$\leq 2^{-(k + 6)} + 2^{-(k + 6)} + 2^{-(k + 5)} = 2^{-(k + 4)}$$ 
This proves the claim.

Now fix a $i < p$. By (III), $\delta_{m_k}^\omega(A^i_{0^k}, A^i_t) < 2^{-(k + 4)}$ for each $t \in {}^k2$. For each $s \in {}^{k + 1}2$, let $b_s^i \in A_{s \upharpoonright k}^i$ be such that $\delta_{m_k}^\omega(b_{0^k\hat{\ }s(k)}^i, b_{s}^i) < 2^{-(k + 4)}$. 

Suppose $s \in {}^{k + 1}2$ and $s(k) = 1$. 
$$|\delta_{m_k}^\omega(b_{0^{k + 1}}^i, b^i_s) - \frac{1}{k + 1}|$$
$$\leq |\delta_{m_k}^\omega(b_{0^{k + 1}}^i, b_s^i) - \delta_{m_k}^\omega(b_{0^{k + 1}}^i, b_{0^k\hat{\ }1}^i)| + |\delta_{m_k}^\omega(b_{0^{k + 1}}^i, b_{0^{k}\hat{\ }1}^i) - \frac{1}{k + 1}|$$
$$< \delta(b_s^i, b^i_{0^k\hat{\ }1}) + 2^{-(k + 4)} \leq 2^{-(k + 4)} + 2^{-(k + 4)} = 2^{-(k + 3)}$$

By (I), (II), and the fact that $b^i_{0^{k + 1}}, b^i_{0^k\hat{\ }1} \in A_{0^k}^i$, there exists some $m_{k + 1} > m_k$ so that 

\noindent (i) $|\delta_{m_k}^{m_{k + 1}}(b^i_{0^{k + 1}}, b_s^i) - \frac{1}{k + 1}| < 2^{-(k + 3)}$ for all $s \in {}^{k + 1}2$ with $s(k) = 1$. 

\noindent (ii) $\delta_{m_{k + 1}}^\omega(b^0_{0^{k + 1}}, b^i_{0^{k + 1}}) < 2^{-(k + 7)}$ for all $i < p$. 

\noindent (iii) $\delta_{m_{k + 1}}^\omega(b^i_{0^{k + 1}}, b^i_s) < 2^{-(k + 5)}$ for all $s \in {}^{k + 1}2$. 

Let $g^i(s) = b^i_s \upharpoonright m_{k + 1}$. Suppose $s(k) = t(k)$. Without loss of generality, suppose $s(k) = t(k) = 1$. Then
$$\delta_{m_k}^{m_{k + 1}}(g^i(s),g^i(t)) \leq \delta_{m_k}^{m_{k + 1}}(b_s^i, b^i_{0^k \hat{\ } 1}) + \delta_{m_k}^{m_{k + 1}}(b^i_{0^k\hat{\ }1}, b^i_t) \leq 2^{-(k + 4)} + 2^{-(k + 4)} = 2^{-(k + 3)} < 2^{-(k + 2)}$$
This establishes (IV).

Suppose $s(k) \neq t(k)$. Without loss of generality, suppose $s(k) = 1$. Hence $t(k) = 0$. 
$$|\delta_{m_k}^{m_{k + 1}}(g^i(s),g^i(t)) - \frac{1}{k + 1}|$$
$$\leq|\delta_{m_k}^{m_{k + 1}}(b_s^i,b_t^i) - \delta_{m_k}^{m_{k + 1}}(b_s^i, b_{0^{k + 1}}^i)| + |\delta_{m_k}^{m_{k + 1}}(b_s^{i},b_{0^{k + 1}}^i) - \frac{1}{k + 1}|$$
$$< \delta_{m_k}^{m_{k + 1}}(b^i_t, b^i_{0^{k + 1}}) + 2^{-(k + 3)} < 2^{-(k + 4)} + 2^{-(k + 3)} < 2^{-(k + 2)}$$
This establishes (V).

Let $s \in {}^{k + 1} 2$. Without loss of generality suppose $s(k) = 0$. Suppose $i,j < p$. Observe:
$$\delta_{m_k}^{m_{k + 1}}(g^i(s),g^j(s)) \leq \delta_{m_k}^\omega(b^i_s,b^j_s) \leq \delta_{m_k}^\omega(b_s^i, b_{0^{k + 1}}^i) + \delta_{m_k}^\omega(b^i_{0^{k + 1}}, b^j_s)$$
$$\leq \delta_{m_k}^\omega(b_s^i, b_{0^{k + 1}}^i) + \delta_{m_k}^\omega(b_{0^{k + 1}}^i, b_{0^{k + 1}}^j) + \delta_{m_k}^\omega(b^j_{0^{k + 1}}, b^j_s)$$
$$\leq \delta_{m_k}^\omega(b_s^i, b_{0^{k + 1}}^i) + \delta_{m_k}^\omega(b_{0^{k + 1}}^i, b_{0^{k + 1}}^0) + \delta_{m_k}^\omega(b_{0^{k + 1}}^0, b^j_{0^{k + 1}}) + \delta_{m_k}^\omega(b^j_{0^{k + 1}}, b^j_s)$$
$$\leq 2^{-(k + 4)} + 2^{-(k + 6)} + 2^{-(k + 6)} + 2^{-(k + 4)} \leq 2^{-(k + 3)} + 2^{-(k + 5)} < 2^{-(k + 2)}$$ 
This establishes (VII).

For $s \in {}^{k + 1}2$, let $B_{s}^i = A_{s \upharpoonright k}^i \cap N_{g^i(s)}$. Apply Lemma \ref{E2 splitting lemma 4} on $\{B_s^i : i < p \wedge s \in {}^{k + 1}2\}$, $\{b_s^i : i < p \wedge s \in {}^{k + 1}2\}$, $r = 2^{-(k + 7)}$, $v = 2^{-(k + 5)}$, and $D_{k + 1}$ to obtain the desired objects $(A_s^i : i < p \wedge s \in {}^{k + 1}2)$ which satisfy the remaining conditions. 

This completes the proof.
\end{proof}

By relativizing to the appropriate parameter, one can obtain the following result.

\Begin{corollary}{E2 big sets with same saturation}
Let $p \in \omega$. Suppose $(X_i : i < p)$ is a collection of $\analytic$ subsets of $\cantorspace$ with the property that for all $i < p$, $E_2 \leq_\borel E_2 \upharpoonright X_i$ and for all $i, j < p$, $[X_i]_{E_2} = [X_j]_{E_2}$. Then there exists a sequence of strictly increasing integers $(m_k : k \in \omega)$ and maps $g^i: \finBinarySequence \rightarrow \finBinarySequence$ satisfying conditions 1 - 5 of Theorem \ref{structure through E2 sets satisfying effective condition}.
\end{corollary}

\Begin{fact}{power E2 characterization}
Let $B \subseteq \cantorspace$ be a $\analytic$ so that $E_2 \leq_\borel E_2 \upharpoonright B$. There there exists a strictly increasing sequence $(m_k : k \in \omega)$ with $m_0 = 0$ and function $g : \finBinarySequence \rightarrow \finBinarySequence$ with the following properties:

\noindent 1. If $|s| = k$, then $|g(s)| = m_k$.

\noindent 2. If $s \subseteq t$, then $g(s) \subseteq g(t)$. 

\noindent 3. If $|s| = |t| = k > 0$ and $s(k - 1) = t(k - 1)$, then $\delta_{m_{k - 1}}^{m_k}(g(s),g(t)) < 2^{-(k + 1)}$. 

\noindent 4. If $|s| = |t| = k > 0$ and $s(k - 1) \neq t(k - 1)$, then $|\delta_{m_{k - 1}}^{m_k}(g(s),g(t)) - \frac{1}{k}| < 2^{-(k + 1)}$.

\noindent 5. Let $\Phi : \cantorspace \rightarrow \cantorspace$ be defined by $\Phi(x) = \bigcup_{n \in \omega} g(x\upharpoonright n)$. Then $\Phi$ is a $\borel$ function such that $\Phi[\cantorspace] \subseteq B$ and $\Phi$ witnesses $E_2 \leq_\borel E_2 \upharpoonright B$. 
\end{fact}

\begin{proof}
This is implicit in \cite{Actions-by-the-Classical-Banach-Spaces}. Also see \cite{Borel-Equivalence-Relations} Theorem 15.4.1 and \cite{Canonical-Ramsey-Theory-on-Polish-Spaces} Theorem 7.43. The proof is quite similar to Theorem \ref{structure through E2 sets satisfying effective condition}.
\end{proof}

\section{$E_2$ Does Not Have the 2-Mycielski Property}\label{E2 Does Not Have the 2-Mycielski Property}

\Begin{theorem}{E2 does not have Mycielski property}
Let $D \subseteq {}^2(\cantorspace)$ be defined by 
$$D = \{(x,y) \in {}^2(\cantorspace) : (\exists i < j)(\delta_i^j(x,y) > 2 \wedge (\forall n)(i \leq n < j \Rightarrow x(n) \neq y(n)))\}$$
$D$ is dense open.

Let $(m_k : k \in \omega)$, $g$, and $\Phi$ be as in Fact \ref{power E2 characterization}. $(\Phi(\tilde 0), \Phi(\widetilde{01})) \notin D$. 

For any $\borel$ set $B$ so that $E_2 \upharpoonright B \equiv_\borel E_2$, $[B]^2_{E_2} \not\subseteq D$. 

$E_2$ does not have the 2-Mycielski property.
\end{theorem}

\begin{proof}
Let $(x,y) \in D$. There is some $i < j$ so that for all $n$ with $i \leq n < j$, $x(n) \neq y(n)$ and $\delta_i^j(x,y) > 2$. Let $\sigma = x \upharpoonright (j + 1)$ and let $\tau = y \upharpoonright (j + 1)$. Then $(x,y) \in N_{\sigma,\tau} \subseteq D$. $D$ is open. 

Let $\sigma,\tau \in \finBinarySequence$ with $|\sigma| = |\tau|$. Let $i = |\sigma|$. Find a $j > i$ so that $\sum_{i \leq n < j} \frac{1}{n + 1} > 2$. Let $\sigma',\tau' \in {}^{j + 1} 2$ be defined by
$$\sigma'(k) = \begin{cases}
\sigma(k) & \quad k < i \\
0 & \quad \text{otherwise}
\end{cases}
\ \ \ \text{ and } \ \ \ 
\tau'(k) = \begin{cases}
\tau(k) & \quad k < i \\
1 & \quad \text{otherwise}
\end{cases}.$$
$N_{\sigma',\tau'} \subseteq N_{\sigma,\tau}$ and $N_{\sigma',\tau'} \subseteq D$. $D$ is dense open.

Note that one must have
$$\sum_{m_{k - 1} \leq m < m_k} \frac{1}{m + 1} \geq 2^{-k - 1}$$
because $\frac{1}{k} - 2^{-k - 1} \geq 2^{-k} - 2^{-k - 1} = 2^{-k - 1}$ and so otherwise, condition 4 could not hold for any $s,t \in {}^k2$ with $s(k - 1) \neq t(k - 1)$. 

This implies that for any $s$ and $t$ so that $s(k - 1) = t(k - 1)$, there must be some $m$ with $m_{k - 1} \leq m < m_k$ so that $g(s)(m) = g(t)(m)$. 

Note that for all $s,t \in {}^{k + 1}2$.
$$\delta_{m_{k - 1}}^{m_{k + 1}}(g(s),g(t)) = \delta_{m_{k -1}}^{m_k}(g(s),g(t)) + \delta_{m_k}^{m_{k + 1}}(g(s),g(t)) \leq \frac{1}{k} + 2^{-k - 1} + \frac{1}{k + 1} + 2^{-k - 2} < 2.$$
Hence for any $s,t$, if there exists $i < j$ so that $\delta_i^j(g(s),g(t)) > 2$, then there is some $k \geq 1$ so that $i \leq m_{k - 1} < m_k < m_{k + 1} \leq j$. 

Now suppose that there is some $i < j$ so that $\delta_i^j(\Phi(\tilde 0),\Phi(\widetilde{01})) > 2$. There is some $k \geq 1$ so that $i \leq m_{k - 1} < m_k < m_{k + 1} \leq j$. Without loss of generality, suppose $k$ is even. $\tilde{0}(k) = 0 = \widetilde{01}(k)$. By the above, there is some $l$ with $m_k \leq l < m_{k + 1}$ so that $\Phi(\tilde{0})(l) = \Phi(\widetilde{01})(l)$. This shows $(\Phi(\tilde{0}), \Phi(\widetilde{01})) \notin D$. $\neg(\tilde{0} \ E_2 \ \widetilde{01})$ so $\neg(\Phi(\tilde{0}) \ E_2 \ \Phi(\widetilde{01}))$. Hence $[\Phi[\cantorspace]]^2_{E_2} \not\subseteq D$. 

It has been shown that for all $(m_k : k \in \omega)$, $g$, and associated $\Phi$, $[\Phi[\cantorspace]]^2_{E_2} \not\subseteq D$. If $B \subseteq \cantorspace$ is $\borel$ with the property that $E_2 \upharpoonright B \equiv_\borel E_2$, then Fact \ref{power E2 characterization} implies that there is some $(m_k : k \in \omega)$ and $g$ so that $\Phi[\cantorspace] \subseteq B$. This shows that for all such $B$, $[B]^2_{E_2} \not\subseteq D$. $E_2$ does not have the 2-Mycielski property. 
\end{proof}

\Begin{theorem}{E2 does not have weak 2-mycielski}
For each $n \in \omega$, let 
$$D_n = \{(x,y) \in {}^2(\cantorspace) : (\exists i < j)(n \leq i < j \wedge \delta_i^j(x,y) > 3 \wedge (\forall m)(i \leq m < j \Rightarrow x(m) \neq y(m)))\}.$$
Each $D_n$ is a dense open subset of ${}^2(\cantorspace)$. Hence $C = \bigcap_{n \in \omega} D_n$ is a comeager subset of ${}^2(\cantorspace)$. 

Suppose $(m_k : k \in \omega)$, $g^0$, $g^1$, $\Phi^0$, and $\Phi^1$ have properties 1 - 6 from Theorem \ref{structure through E2 sets satisfying effective condition}. Then $(\Phi^0(\tilde 0), \Phi^1(\widetilde{01})) \notin C$. 

For any $\borel$ sets $B_0$ and $B_1$ with $E_2 \leq_\borel E_2 \upharpoonright B_0$, $E_2 \leq_\borel E_2 \upharpoonright B_1$, and $[B_0]_{E_2} = [B_1]_{E_2}$, $B_0 \times_{E_2} B_1 \not\subseteq C$. 

$E_2$ does not have the weak $2$-Mycielski property.
\end{theorem}

\begin{proof}
Using Theorem \ref{structure through E2 sets satisfying effective condition}, if $|s| = |t| = k > 0$ and $s(k - 1) \neq t(k - 1)$, then $|\delta_{m_{k - 1}}^{m_k}(g^0(s),g^1(t)) - \frac{1}{k}| < 2^{-k}$. To see this: 
$$|\delta_{m_{k - 1}}^{m_k}(g^0(s),g^1(t)) - \frac{1}{k}|$$
$$\leq |\delta_{m_{k - 1}}^{m_k}(g^0(s), g^1(t)) - \delta_{m_{k - 1}}^{m_k}(g^1(s), g^1(t))| + |\delta_{m_{k - 1}}^{m_k}(g^1(s), g^1(t)) - \frac{1}{k}|$$
$$\leq \delta_{m_{k - 1}}^{m_k}(g^0(s),g^1(s)) + 2^{-(k + 1)} \leq 2^{-(k + 1)} + 2^{-(k + 1)} = 2^{-k}$$

Also if $|s| = |t| = k > 0$ and $s(k - 1) = t(k - 1)$, then $\delta_{m_{k - 1}}^{m_k}(g^0(s),g^1(t)) < 2^{-k}$. To see this: 
$$\delta_{m_{k - 1}}^{m_k}(g^0(s),g^1(t)) \leq \delta_{m_{k - 1}}^{m_k}(g^0(s),g^1(s)) + \delta_{m_{k - 1}}^{m_k}(g^1(s),g^1(t)) = 2^{-(k + 1)} + 2^{-(k + 1)} = 2^{-k}$$ 

$D_n$ is dense open by the same argument as in Theorem \ref{E2 does not have weak 2-mycielski}.

Note that if $k > 0$, then
$$\sum_{m_{k - 1} \leq m < m_k} \frac{1}{m + 1} \geq 2^{-k}$$
because $\frac{1}{k} - 2^{-k} \geq 2^{-(k - 1)} - 2^{-k} = 2^{-k}$ and so otherwise $|\delta_{m_{k - 1}}^{m_k}(g^0(s), g^1(s)) - \frac{1}{k}| < 2^{-k}$ could not hold. 

Therefore if $|s| = |t| > 0$, and $s(k - 1) = t(k - 1)$, then there must be some $m$ with $m_{k - 1} \leq m < m_k$ so that $g^0(s)(m) = g^1(t)(m)$. 

Note that for all $s,t \in {}^{k + 1}2$ with $k > 0$, 
$$\delta_{m_{k - 1}}^{m_{k + 1}}(g(s),g(t)) = \delta_{m_{k -1}}^{m_k}(g(s),g(t)) + \delta_{m_k}^{m_{k + 1}}(g(s),g(t)) \leq \frac{1}{k} + 2^{-k} + \frac{1}{k + 1} + 2^{-(k + 1)} < 3.$$
Hence for any $s,t$, if there exists $i < j$ so that $\delta_i^j(g^0(s),g^1(t)) > 3$, then there is some $k \geq 1$ so that $i \leq m_{k - 1} < m_k < m_{k + 1} \leq j$. 

Now by essentially the same argument as in Theorem \ref{E2 does not have Mycielski property}, $(\Phi^0(\tilde 0), \Phi^1(\widetilde{01})) \notin D_{m_0}$. Hence $(\Phi^0(\tilde 0), \Phi^{1}(\widetilde{01})) \notin C$. 

Now suppose that $B_0$ and $B_1$ are some $\borel$ sets so that $E_2 \leq_\borel E_2 \upharpoonright B_0$, $E_2 \leq_\borel E_2 \upharpoonright B_1$, and $[B_0]_{E_2} = [B_1]_{E_2}$. By Corollary \ref{E2 big sets with same saturation}, there is a sequence $(m_k : k \in \omega)$, $g^0$, $g^1$, $\Phi^0$ and $\Phi^1$ as above so that $\Phi^i[\cantorspace] \subseteq B_i$. By the earlier argument, $B_0 \times_{E_2} B_1 \not\subseteq C$. Hence $E_2$ does not have the weak $2$-Mycielski property.
\end{proof}

\section{Surjectivity and Continuity Aspects of $E_2$}\label{Surjectivity and Continuity Aspects of E2}

Fact \ref{power E2 characterization} states that every $\analytic$ set $B \subseteq \cantorspace$ so that $E_2 \leq_\borel E_2 \upharpoonright B$ has a closed set $C \subseteq B$ so that $E_2 \equiv_\borel E_2 \upharpoonright C$. Fact \ref{power E2 characterization} even asserts that $C$ is the body of a tree on $2$ with a specific structure:

\Begin{definition}{E2 tree}
A tree $p \subseteq \finBinarySequence$ is an $E_2$-tree if and only if there is some sequence $(m_k : k \in \omega)$ and map $g : \finBinarySequence \rightarrow \finBinarySequence$ satisfying the conditions of Fact \ref{power E2 characterization} so that $p$ is the downward closure of $g[\finBinarySequence]$. Note that if $\Phi$ is the map associated with $(m_k : k \in \omega)$ and $g$, then $[p] = \Phi[\cantorspace]$. 
\end{definition}

The following notation is used to avoid some very tedious superscripts and subscripts in the following results:

\Begin{definition}{E2 section specific distance notation}
If $x, y \in \cantorspace$ and $m, n \in \omega$ with $m \leq n$, then let $\varsigma(m,n,x,y) = \delta_{m}^n(x,y)$.
\end{definition}

\Begin{fact}{cont function surjective on 3-cantorspace on E2-tree}
There is a continuous function $P : [\cantorspace]^3_{E_2} \rightarrow {}^\omega 3$ so that for any $E_2$-tree $p$, $P[[[p]]^3_{E_2}] = {}^\omega 3$. 
\end{fact}

\begin{proof}
For  $(x,y) \in [\cantorspace]^2_{E_2}$ and any $n,m \in \omega$, define
$$S_{n,m}(x,y) = \min \{k \in \omega : \delta_n^k(x,y) > 3^{m + 2}\}$$
Each $S_{n,m}$ is continuous on $[\cantorspace]^3_{E_2}$. 

If $(x,y,z) \in [\cantorspace]^3_{E_2}$, then define a strictly increasing sequence of integers $(L_n : n \in \omega)$ by recursion as follows: Let $L_0 = 0$. Given $L_n$, let
$$L_{n + 1} = \min\{S_{L_n,n}(x,y), S_{L_n,n}(x,z), S_{L_n,n}(y,z)\}$$
(It is implicit that $L_n$ depends on the triple $(x,y,z)$.) By induction, it can be shown that each $L_n$ as a function of $(x,y,z)$ is continuous on $[\cantorspace]^3_{E_2}$. 

Define
$$P(x,y,z)(n) = \begin{cases}
0 & \quad S_{L_n,n}(x,y) \leq S_{L_n,n}(x,z) \text{ and } S_{L_n,n}(x,y) \leq S_{L_n,n}(y,z) \\
1 & \quad S_{L_n,n}(x,z) < S_{L_n,n}(x,y) \text{ and } S_{L_n,n}(x,z) \leq S_{L_n,n}(y,z) \\
2 & \quad S_{L_n,n}(y,z) < S_{L_n,n}(x,y) \text{ and } S_{L_n,n}(y,z) < S_{L_n,n}(x,z)
\end{cases}$$ 
$P$ is continuous on $[\cantorspace]^3_{E_2}$. 

Also define the sequence of integers $(N_n : n \in \omega)$ by recursion as follows: Let $N_0 = 0$ and if $N_n$ has been defined, then let
$$N_{n + 1} = \min\left\{k \in \omega : \sum_{N_i \leq i < k} \left(\frac{1}{i + 1} - 2^{-(i + 2)}\right) > 3^{n + 2}\right\}$$

Note that $N_{n + 1} > N_n + 2$ for each $n \in \omega$. By the definition of $N_{n + 1}$, one has that
$$\sum_{N_n \leq i < N_{n + 1} - 1} \left(\frac{1}{i + 1} - 2^{-i - 2}\right) \leq 3^{n + 2}.$$

These two facts imply
\begin{equation}\label{E2 eq 1}
\sum_{N_n \leq i < N_{n + 1}} \frac{1}{i + 1} \leq 3^{n + 2} + \frac{1}{N_{n + 1}} + \sum_{N_n \leq i < N_{n + 1} - 1} 2^{-i - 2} < 3^{n + 2} + 1
\end{equation} 

Let $k_n = N_{n + 1} - N_n$. Fix a $v \in {}^\omega 3$. Define $\sigma_n, \tau_n \in {}^{k_n}2$ by
$$\sigma_n = \begin{cases}
\tilde 1 \upharpoonright k_n & \quad v(n) = 0 \\
\widetilde{01} \upharpoonright k_n & \quad \text{otherwise}
\end{cases}
\ \ \ \ \ \ 
\tau_n = \begin{cases}
\tilde 1 \upharpoonright k_n & \quad v(n) = 1 \\
\widetilde{10} \upharpoonright k_n & \quad \text{otherwise}
\end{cases}$$

Let $x = \tilde 0$, $y = \sigma_0\hat{\ }\sigma_1\hat{\ } \sigma_2\hat{\ } ...$, and $z = \tau_0 \hat{\ } \tau_1 \hat{\ } \tau_2 ...$. Note $(x,y,z) \in [\cantorspace]^3_{E_2}$. 

Fix an $E_2$-tree $p$. Let $(m_k : k \in \omega)$, $g : \finBinarySequence \rightarrow \finBinarySequence$, and $\Phi : \cantorspace \rightarrow \cantorspace$ be the associated objects of $p$ coming from the definition of an $E_2$-tree.

Suppose $v(n) = 0$, then
\begin{equation}\label{E2 eq 2}\varsigma(m_{N_n}, m_{N_{n + 1}}, \Phi(x), \Phi(y)) = \sum_{N_n \leq i < N_{n + 1}} \varsigma(m_i, m_{i+ 1}, \Phi(x), \Phi(y)) > \sum_{N_n \leq i < n_{n + 1}} \left(\frac{1}{i + 1} - 2^{-i - 2}\right) > 3^{n + 2}\end{equation}
using the definition of $N_{n + 1}$. Also
\begin{equation}\label{E2 eq 3}\varsigma(m_{N_n}, m_{N_{n + 1}}, \Phi(x),\Phi(y)) < \sum_{N_n \leq i < N_{n + 1}} \left(\frac{1}{i + 1} + 2^{-n - 2}\right) < 3^{n + 2} + 1 + \sum_{N_n \leq i < N_{n + 1}} 2^{-i - 2} < 3^{n + 2} + \frac{3}{2}\end{equation}
using equation (\ref{E2 eq 1}) for the second inequality.

Note also
$$\varsigma(m_{N_n}, m_{N_{n + 1}},\Phi(x), \Phi(z)) = \varsigma(m_{N_n}, m_{N_n + 1},\Phi(x), \Phi(y)) + \sum_{N_n < i < N_{n + 1}} \varsigma(m_i,m_{i + 1},\Phi(x),\Phi(z))$$
$$< \frac{1}{N_n + 1} + 2^{-N_n - 2} + \frac{1}{2}\sum_{N_n < i < N_{n + 1}} \frac{1}{i + 1} + \sum_{N_n < i < N_{n + 1}} 2^{-i - 2}$$
$$= \frac{1}{2}\sum_{N_n \leq i < N_{n + 1}} \frac{1}{i + 1} + \frac{1}{2}\left(\frac{1}{N_n + 1}\right) + \sum_{N_n \leq i < N_{n + 1}} 2^{-i - 2} < \frac{1}{2}(3^{n + 2}) + \frac{3}{2}$$
using equation (\ref{E2 eq 1}). In summary,
\begin{equation}\label{E2 eq 4}
\varsigma(m_{N_n}, m_{N_{n + 1}}, \Phi(x), \Phi(z)) < \frac{1}{2}(3^{n + 2}) + \frac{3}{2}
\end{equation} 
Similarly, $\varsigma(m_{N_n}, m_{N_{n + 1}}, \Phi(y), \Phi(z)) < \frac{1}{2}(3^{n + 2}) + \frac{3}{2}$. The case for $v(n) = 1$ and $v(n) = 2$ are similar.

It remains to show that $P(\Phi(x),\Phi(y), \Phi(z)) = v$. In the following, let $(L_n : n \in \omega)$ be the sequence defined as above using $(\Phi(x), \Phi(y), \Phi(z))$. The following statements will be proved by induction on $n$:

\noindent(I) $m_{N_n} < L_{n + 1} \leq m_{N_{n + 1}}$.

\noindent(II) $P(\Phi(x), \Phi(y), \Phi(z))(n) = v(n)$.

\noindent (III) The following holds:
$$\max\{\varsigma(L_{n + 1}, m_{N_{n + 1}}, \Phi(x),\Phi(y)), \varsigma(L_{n + 1}, m_{N_{n + 1}}, \Phi(x),\Phi(z)), \varsigma(L_{n + 1}, m_{N_{n + 1}}, \Phi(y),\Phi(z))\} < \frac{1}{2}(3^{n + 2}) + \frac{3}{2}.$$

Suppose properties (I), (II), and (III) holds for all $k < n$. Suppose $v(n) = 0$. (The other cases are similar.) 

$L_n \leq m_{N_n}$ by definition if $n = 0$ and by the induction hypothesis otherwise. Therefore, 
$$\varsigma(L_n,m_{N_{n + 1}},\Phi(x),\Phi(y)) \geq \varsigma(m_{N_n},m_{N_{n + 1}},\Phi(x),\Phi(y)) > 3^{n + 2}$$
using equation (\ref{E2 eq 2}). This shows $L_{n + 1} \leq S_{L_n,n}(\Phi(x),\Phi(y)) \leq m_{N_{n + 1}}$. Using the induction hypothesis or the definition when $n = 0$, 
$$\max\left\{\varsigma(L_n,m_{N_n},\Phi(x),\Phi(y)), \varsigma(L_n,m_{N_n},\Phi(x),\Phi(z)), \varsigma(L_n,m_{N_n},\Phi(y),\Phi(z))\right\} < \frac{1}{2}(3^{n + 1}) + \frac{3}{2} < 3^{n + 2}.$$ 
Hence $L_{n + 1} \leq m_{N_n}$ is impossible. This proves (I).

Observe that
$$\varsigma(L_n,m_{N_{n + 1}},\Phi(x),\Phi(z)) = \varsigma(L_n,m_{N_n},\Phi(x),\Phi(z)) + \varsigma(m_{N_n},m_{N_{n + 1}},\Phi(x),\Phi(z))$$
$$ < \frac{1}{2}(3^{n + 1}) + \frac{3}{2} + \frac{1}{2}(3^{n + 2}) + \frac{3}{2} \leq 3^{n + 2}$$
using the induction hypothesis and equation (\ref{E2 eq 4}). This shows $S_{L_n,n}(\Phi(x),\Phi(z)) > m_{N_{n + 1}}$. Similarly, $S_{L_n,n}(\Phi(y),\Phi(z)) > m_{N_{n + 1}}$. $S_{L_n,n}(\Phi(x),\Phi(y)) \leq m_{N_{n + 1}}$ has already been shown above. Thus 
$$P(\Phi(x),\Phi(y),\Phi(z))(n) =  0 = v(n).$$
This shows (II).

Note that
$$\varsigma(L_{n + 1},m_{N_{n + 1}},\Phi(x),\Phi(z)) \leq \varsigma(m_{N_n},m_{N_{n + 1}},\Phi(x),\Phi(z)) < \frac{1}{2}(3^{n + 2}) + \frac{3}{2}$$
using equation (\ref{E2 eq 4}). Similarly, 
$$\varsigma(L_{n + 1},m_{N_{n + 1}},\Phi(y),\Phi(z)) < \frac{1}{2}(3^{n + 2}) + \frac{3}{2}.$$
Finally,
$$\varsigma(L_n,m_{N_{n + 1}},\Phi(x),\Phi(y)) = \varsigma(m_{N_n},m_{N_{n + 1}},\Phi(x),\Phi(y)) - \big(\varsigma(L_n,L_{n + 1 },\Phi(x),\Phi(y)) - \varsigma(L_n,m_{N_n},\Phi(x),\Phi(y))\big)$$
$$< 3^{n + 2} + \frac{3}{2} - 3^{n + 2} + \frac{1}{2}(3^{n + 1}) + \frac{3}{2} < \frac{1}{2}(3^{n + 2}) + \frac{3}{2}$$
using equation (\ref{E2 eq 3}), the definition of the sequence $(L_n : n \in \omega)$, and the induction hypothesis. This proves (III). 
\end{proof}

\Begin{theorem}{cont function surjective on E2-tree}
There is a continuous function $Q : [\cantorspace]^3_{E_2} \rightarrow \cantorspace$ so that for any $E_2$-tree $p$, $Q[[[p]]_{E_2}^3] = \cantorspace$. 

There is a $\lborel$ function $K : {}^3(\cantorspace) \rightarrow \cantorspace$ so that on any $\analytic$ set $A$ with $E_2 \leq_\borel E_2 \upharpoonright A$, $K[[A]^3_{E_2}] = \cantorspace$ (and in particular the image meets each $E_2$-equivalence class).
\end{theorem}

\begin{proof}
$Q$ can be obtained by composing the function from Fact \ref{cont function surjective on 3-cantorspace on E2-tree} with a homeomorphism from ${}^\omega 3 \rightarrow {}^\omega 2$. 

$K$ can be obtained by mapping elements of ${}^3(\cantorspace) \setminus [\cantorspace]^3_{E_2}$ to $\tilde 0$ and mapping elements in $[\cantorspace]^3_{E_2}$ according to $Q$. Note that by Fact \ref{power E2 characterization}, every such set $A$ contains an $E_2$-tree. 
\end{proof}

\Begin{lemma}{useful non-continuity lemma}
Let $X$ and $Y$ be topological spaces. Let $\mathscr{B} \subseteq \mathscr{P}(X)$ be a nonempty family of subsets of $X$. Let $A \subseteq Y$ be Borel. Suppose $f : X \rightarrow Y$ is a Borel function with the property that for all $B \in \mathscr{B}$ and open $U \subseteq X$ with $U \cap B \neq \emptyset$, there exists an $x \in B$ with $f(x) \notin A$ and an $x' \in U \cap B$ with $f(x') \in A$. Then there is a Borel function $g : X \rightarrow Y$ such that $g \upharpoonright B$ is not continuous for all $B \in \mathscr{B}$. 
\end{lemma}

\begin{proof}
The assumption above implies that $A$ is Borel but not equal to either $\emptyset$ or $Y$. The topology of $Y$ is not $\{\emptyset, Y\}$. There exists some $y_1,y_2 \in Y$ and open set $V \subseteq Y$ with $y_1 \in V$ and $y_2 \notin V$. Define
$$g(x) = \begin{cases}
y_1 & \quad f(x) \notin A \\
y_2 & \quad f(x) \in A
\end{cases}.$$

Suppose there was some $B \in \mathscr{B}$ so that $g \upharpoonright B$ is a continuous function. By the assumptions, there is a $x \in B$ so that $f(x) \notin A$. So $g(x) = y_1$. By continuity, $g^{-1}[V] \cap B$ is a nonempty open set containing $x \in B$. There is some $U \subseteq X$ open so that $g^{-1}[V] \cap B = U \cap B$. By the assumptions, there some $x' \in U \cap B$ so that $f(x') \in A$. Hence $g(x') = y_2 \notin V$. Contradiction.
\end{proof}

\Begin{fact}{non-continuity function into 3omega on E2-tree}
There is a $\borel$ function $P' : [\cantorspace]^3_{E_2} \rightarrow {}^\omega 3$ so that on any $E_2$-tree $p$, $P' \upharpoonright[[p]]^3_{E_2}$ is not continuous. 
\end{fact}

\begin{proof}
This result is proved by applying Lemma \ref{useful non-continuity lemma} to $X = [\cantorspace]^3_{E_2}$, $Y = {}^\omega 3$, $\mathscr{B} = \{[[p]]_{E_2}^3 : p \text{ is an $E_2$-tree}\}$, $f$ is the function $P$ from Fact \ref{cont function surjective on 3-cantorspace on E2-tree}, and $A = \{z \in {}^\omega 3 : (\exists k)(\forall n > k)(z(n) = 0)\}$. It remains to show that these objects satisfy the required properties of Lemma \ref{useful non-continuity lemma}.

Fix an $E_2$-tree $p$. Let $(m_k : k \in \omega)$, $g : \finBinarySequence \rightarrow \finBinarySequence$, and $\Phi$ be the objects associated with $p$ from the definition of an $E_2$-tree. Fact \ref{cont function surjective on 3-cantorspace on E2-tree} implies $P[[[p]]^3_{E_2}] = {}^\omega 3$. Hence there is some $(x,y,z) \in [[p]]^3_{E_2}$ so that $P(x,y,z) \notin A$. Let $U \subseteq [\cantorspace]^3_{E_2}$ be open so that $U \cap [[p]]^3_{E_2} \neq \emptyset$. There are some $s,t,u \in \finBinarySequence$ so that $\emptyset \neq N_{s,t,u} \cap [[p]]_{E_2}^3 \subseteq U \cap [[p]]^3_{E_2}$. Let $x' = s \hat{\ }\tilde 0$, $y' = t \hat{\ }\tilde 1$, and $z' = u\hat{\ }\widetilde{01}$. Using the computation from the proof of Fact \ref{cont function surjective on 3-cantorspace on E2-tree}, if $k$ is chosen so that $L_k \geq m_{|s|}$, then for all $n > k$, $P(\Phi(x'),\Phi(y'),\Phi(z'))(n) = 0$. Hence $(\Phi(x'), \Phi(y'), \Phi(z')) \in U \cap [[p]]_{E_2}^3$ and $P(\Phi(x'),\Phi(y'), \Phi(z')) \in A$. 
\end{proof}

\Begin{theorem}{borel function E2-big set not continuous}
There is a $\borel$ function $K : {}^3(\cantorspace) \rightarrow \cantorspace$ so that for any $\analytic$ set $A$ with $E_2 \leq_\borel E_2 \upharpoonright A$, $K \upharpoonright A$ is not continuous.
\end{theorem}

\begin{proof}
Use the usual arguments to adjust the domain and range of the function from Fact \ref{non-continuity function into 3omega on E2-tree}. Then apply Fact \ref{power E2 characterization}.
\end{proof}

\Begin{corollary}{a proof of failure of 3-Mycielski property for E2}
$E_2$ does not have the $3$-Mycielski proprety. 
\end{corollary}

\begin{proof}
Let $C \subseteq {}^3(\cantorspace)$ be any comeager set so that $K \upharpoonright C$ is continuous. Then $C$ witnesses the failure of the $3$-Mycielski property for $E_2$.
\end{proof}

\section{The Structure of $E_3$}\label{The Structure of E3}
This section will give the characterization of $E_3$-big $\analytic$ sets coming from its dichotomy result. See the references mentioned below for the details.

\Begin{definition}{E3 definition}
$E_3$ is the equivalence relation on $\pcantorspace$ defined by $x \ E_3 \ y$ if and only if $(\forall n)(x(n) \ E_0 \ y(n))$. 
\end{definition}

\Begin{definition}{grid notation}
Let $\langle \cdot, \cdot \rangle : {}^2 \omega \rightarrow \omega$ be some recursive pairing function. 

Let $\pi_1, \pi_2 : {}^2\omega \rightarrow \omega$ be projections onto the first and second coordinate, respectively.

Let $A \subseteq \omega$. Define $\dom(A) = \{(i,j) : \langle i, j \rangle \in A\}$. 

If $A \subseteq \omega$ is finite, let $L(A) = \sup \pi_1[\dom(A)]$. 

If $s \in {}^n 2$, let $\grid(s) : \mathrm{dom}(n) \rightarrow 2$ be defined by $\grid(s)(i,j) = s(\langle i, j \rangle)$.

\end{definition}

\Begin{definition}{E3 and group action}
$\Z_2$ is the group $(2, +^{\Z_2},0)$ where $+^{\Z_2}$ is modulo 2 addition and $0$ denotes the identity element. 

Let ${}^\omega \Z_2 = (\cantorspace, +^{{}^\omega \Z_2}, \tilde 0)$ where $+^{{}^\omega \Z_2}$ is the coordinate-wise addition of $+^{\Z_2}$ and $\tilde 0$ is the constant $0$ function. 

Let ${}^\omega({}^\omega \Z_2) = (\pcantorspace, +^{{}^\omega({}^\omega \Z_2)}, \bar{0})$ where $+^{{}^\omega({}^\omega \Z_2)}$ is the coordinate-wise addition of $+^{{}^\omega\Z_2}$ and $\bar{0} \in \pcantorspace$ is defined by $\bar{0}(k)(j) = 0$ for all $k,j \in \omega$. 

Let $Z \subseteq \cantorspace$ be defined by $Z = \{x \in \cantorspace : (\exists k)(\forall j > k)(x(j) = 0)\}$. 

$\group = (Z, +^{{}^\omega \Z_2}, \tilde{0})$ is the $\omega$-direct product of $\Z_2$. 

$\wgroup = ({}^\omega Z, +^{{}^\omega({}^\omega \Z_2)}, \bar{0})$ is the $\omega$-product of $\group$. 

If $g \in \wgroup$, define $\supp(g) = \{n \in \omega : g(n)\neq \tilde{0}\}$. 

$\wgroup$ acts on $\pcantorspace$ by left addition when $\pcantorspace$ is considered as ${}^\omega({}^\omega \Z_2)$. That is, $g \cdot x = g +^{{}^\omega({}^\omega\Z_2)} \ x$. 
\end{definition}

\Begin{fact}{E3 group action}
$x \ E_3 \ y$ if and only if there is a $g \in \wgroup$ so that $x = g \cdot y$. 
\end{fact}

\Begin{definition}{grid system}
A grid system is a sequence $(g_{s,t} : s,t \in \finBinarySequence \wedge |s| = |t|\})$ in $\wgroup$ with the following properties:

(I) If $s,t,u \in {}^n 2$ for some $n \in \omega$, then $g_{s,u} = g_{t,u} +^{\wgroup} g_{s,t}$.

(II) For all $m \leq n$, $s,t \in {}^n 2$, $u,v \in {}^m2$ and $l \in \pi_1[\dom(n)]$ with $u \subseteq s$, and $v \subseteq t$, if 
$$\grid(s) \upharpoonright \dom(n \setminus m) \cap ((l + 1)\times\omega) = \grid(t) \upharpoonright \dom(n \setminus m) \cap ((l + 1)\times\omega)$$
then for all $i \leq l$, $g_{s,t}(i) = g_{u,v}(i)$. 
\end{definition}

\Begin{fact}{E3 characterization}
Let $B \subseteq \pcantorspace$ be $\analytic$ so that $E_3 \upharpoonright B \equiv_\borel E_3$. Then there is a continuous injective map $\Phi : \pcantorspace \rightarrow \pcantorspace$, a grid system $(g_{s,t} : s,t \in \finBinarySequence \wedge |s| = |t|)$, and sequences $(k_i : i \in \omega)$ and $(p_{m,i} : m,i \in \omega)$ in $\omega$ with the following properties:

(i) $\Phi[\pcantorspace] \subseteq B$. 

(ii) If $s,t \in {}^n2$, then $\supp(g_{s,t}) \subseteq (k_{L(n)} + 1)$.

(iii) For each $m \in \omega$, $(k_i : i \in \omega)$ and $(p_{m,i} : i \in \omega)$ are strictly increasing sequences. 

(iv) For all $x,y \in \pcantorspace$ and $m,j \in \omega$, if $x(m)(j) = 0$ and $y(m)(j) = 1$, then $\Phi(x)(k_m)(p_{m,j}) = 0$ and $\Phi(y)(k_m)(p_{m,j}) = 1$. 

(v) Let $x,y \in \pcantorspace$ and $l \in \omega$. Suppose 
$$\Big(\forall (i,j) \in ((l + 1) \times \omega) \setminus dom(n)\Big) \Big(x(i)(j) = y(i)(j)\Big).$$
Let $s,t \in {}^n 2$ be such that for all $(i,j) \in \dom(n)$, $\grid(s)(i,j) = x(i)(j)$ and $\grid(t)(i,j) = y(i)(j)$. Then $(g_{s,t} \cdot \Phi(x))(l) = \Phi(y)(l)$. 
\end{fact}

\begin{proof}
This is implicit in \cite{Recent-Developments-in-the-Theory-of-Borel-Reducibilities}. See the presentation in \cite{Borel-Equivalence-Relations} Chapter 14, especially Section 14.5 and 14.6. 
\end{proof}

Note that $\Phi$ as above is an $E_3$ reduction.

\section{$E_3$ Does Not Have the $2$-Mycielski Property}\label{E3 Does Not Have the 2-Mycielski Property}

\Begin{definition}{grid basic open neighborhood}
For each $s \in \finBinarySequence$, let $N_{\grid(s)} = \{x \in \pcantorspace : (\forall (i,j) \in \dom(|s|))(x(i)(j) = s(\langle i, j \rangle))\}$. 

Each $N_{\grid(s)}$ is an open neighborhood of $\pcantorspace$ and also the collection $\{N_{\grid(s)} : s \in \finBinarySequence\}$ forms a basis for the topology of $\pcantorspace$. 
\end{definition}

When $\sigma : m \rightarrow \finBinarySequence$, then $N_\sigma$ will refer to the usual basic open neighborhood of $\pcantorspace$. Both types of open sets will be used in the proof of the following result.

\Begin{theorem}{E3 does not have mycielski property}
Let $D = \{(x,y) \in {}^2(\pcantorspace) : x(0) \neq y(0)\}$. $D$ is dense open. 

For all $\analytic$ sets $B \subseteq \pcantorspace$ with $E_3 \upharpoonright B \equiv_\borel E_3$, $[B]^2_{E_3} \not\subseteq D$. 

$E_3$ does not have the $2$-Mycielski property.
\end{theorem}

\begin{proof}
Suppose $(x,y) \in D$. There is some $n$ so that $x(0)(n) \neq y(0)(n)$. Let $\sigma,\tau : 1 \rightarrow \finBinarySequence$ be defined by $\sigma(0) = x(0) \upharpoonright (n + 1)$ and $\tau(0) = y(0) \upharpoonright (n + 1)$. Then $(x,y) \in N_{\sigma,\tau} \subseteq D$. $D$ is open. 

Suppose $\sigma,\tau : m \rightarrow \finBinarySequence$ have the property that for all $k < m$, $|\sigma(k)| = |\tau(k)|$. Define $\sigma', \tau' : m \rightarrow \finBinarySequence$ by
$$\sigma'(k) = \begin{cases}
\sigma(k) & \quad k \neq 0 \\
\sigma(k)\hat{\ }0 & \quad k = 0 
\end{cases}
\ \ \ \text{ and } \ \ \ 
\tau'(k) = \begin{cases}
\tau(k) & \quad k \neq 0 \\
\tau(k)\hat{\ }1 & \quad k = 0
\end{cases}$$
$N_{\sigma',\tau'} \subseteq N_{\sigma,\tau}$ and $N_{\sigma',\tau'} \subseteq D$. $D$ is dense open.

Fix $\Phi$ and the other objects specified by Fact \ref {E3 characterization}. Note that for any $s \in \finBinarySequence$, $g_{s,s} = \bar{0}$. In particular, $g_{\emptyset,\emptyset} = \bar{0}$. 

Let $\rho_n : 1 \rightarrow \finBinarySequence$ be defined by $\rho_n(0) = \Phi(\bar{0})(0)\upharpoonright n$. If $s \in \finBinarySequence$, then define $x_s \in \pcantorspace$ by
$$x_s(i)(j) = \begin{cases}
s(\langle i, j \rangle) & \quad (i,j) \in \dom(|s|) \\
0 & \quad \text{otherwise}
\end{cases}.$$

Let $s_0 = \emptyset$. 

Suppose $s_n \in \finBinarySequence$ has been defined so that $x_{s_n}(0) = \tilde 0$ and $\Phi(x_{s_n})(0) = \Phi(\bar{0})(0)$. By continuity, find some $u \in \finBinarySequence$ with $s_n \subseteq u$ and $x_{s_n} \in N_{\grid(u)}$ so that $N_{\grid(u)} \subseteq \Phi^{-1}[N_{\rho_{n + 1}}]$. Now find the least $k > |u|$ so that $k = \langle 1, q \rangle$ for some $q \in \omega$. Let $s_{n + 1} \supseteq u$ be of length $k + 1$ defined by
$$s_{n + 1}(j) = \begin{cases}
u(j) & \quad j < |u| \\
1 & \quad j = k \\
0 & \quad \text{otherwise}
\end{cases}.$$

Note that since $x_{s_{n + 1}}(0) = \tilde 0 = \bar{0}(0)$, $(g_{\emptyset,\emptyset} \cdot \Phi(x_{s_{n + 1}}))(0) = \Phi(\bar{0})(0)$ by condition (v) of Definition \ref{E3 characterization}. This implies that $\Phi(x_{s_{n + 1}})(0) = \Phi(\bar{0})(0)$. 

Now define $x \in \pcantorspace$ by
$$x(i)(j) = \left(\bigcup_{n \in \omega} \grid(s_n)\right)(i,j).$$

Since $N_{\grid(s_n)} \subseteq \Phi^{-1}[N_{\rho_n}]$ for all $n \in \omega$, $\Phi(x)(0) = \Phi(\bar{0})(0)$. Hence $(\Phi(x), \Phi(\bar{0})) \notin D$. However, there are infinitely many $q \in \omega$ so that $x(1)(q) = 1$. Since $\Phi$ is an $E_3$ reduction, $\neg(\Phi(x)\ E_3 \  \Phi(\bar{0}))$. 

It has been shown that for any map $\Phi$ as in Fact \ref{E3 characterization}, $[\Phi[\pcantorspace]]^2_{E_3} \not\subseteq D$. Since any $\analytic$ set $B \subseteq \pcantorspace$ with the property that $E_3 \upharpoonright B \equiv_{\borel} E_3$ has some such map $\Phi$ so that $\Phi[\pcantorspace] \subseteq B$, no such $B$ can have the property that $[B]^2_{E_3} \subseteq D$. $E_3$ does not have the 2-Mycielski property.
\end{proof}

\section{Completeness of Ultrafilters on Quotients}\label{Completeness of Ultrafilters on Quotients}
Without the axiom of choice, the notion of completeness of ultrafilters needs to be defined with care.

\Begin{definition}{completeness of ultrafilter}
Let $X$ be a set. Let $U$ be an ultrafilter on $X$. Let $I$ be a set. $U$ is $I$-complete if and only if for any set $J$ which inject into $I$ but is not in bijection with $I$, and any injective function $f : J \rightarrow U$, $\bigcap_{j \in J} f(j) \in U$. 

$U$ is $I^+$-complete if and only if for all $J$ which inject into $I$ and all injective functions $f : J \rightarrow U$, $\bigcap_{j \in J} f(j) \in U$. 
\end{definition}

$\aleph_1$-complete is often called countably complete. A well-known result is that there are no countably complete ultrafilters on $\cantorspace$. There are countably complete ultrafilters on quotients of Polish spaces by equivalence relations.

\Begin{fact}{comeager ultrafilter}
Let $\mathcal{C} \subseteq \mathcal{P}(\cantorspace \slash E_0)$ be defined by $A \in \mathcal{C}$ if and only if $\bigcup A$ belongs to the comeager filter on $\cantorspace$. $\mathcal{C}$ is a countably complete ultrafilter on $\cantorspace \slash E_0$. 
\end{fact}

\begin{proof}
$\mathcal{C}$ is an ultrafilter follows from the generic ergodicity of $E_0$. Countable completeness is clear; in fact under $\mathsf{AD}$, every ultrafilter is countably complete.
\end{proof}

A natural question is whether this ultrafilter or any ultrafilter on $\cantorspace \slash E_0$ could be more than just countably complete. $\bbR$ injects into $\cantorspace \slash E_0$. Is $\mathcal{C}$ $\bbR^+$-complete? Note the function $f$ in Definition \ref{completeness of ultrafilter} is required to be injective.  Otherwise this notion becomes clearly trivial  using the function $f : \bbR \rightarrow \mathcal{C}$ defined by $x \mapsto (\cantorspace \slash E_0) \setminus \{[x]_{E_0}\}$. The next fact will show using a modification of the above function that there are no nonprincipal $\bbR^+$-complete ultrafilters on quotients of Polish spaces. 

\Begin{fact}{no R+ complete ultrafilter}
$(\mathsf{ZF + AD})$ Suppose $E$ is an equivalence relation on a Polish space $X$ so that $= \leq E$ (where $\leq$ denotes the existence of a reduction). Then no nonprincipal ultrafilter on $X \slash E$ is $\bbR^+$-complete.
\end{fact}

\begin{proof}
Let $U$ be a nonprincipal $(\bbR)^+$-complete ultrafilter on $X \slash E$. Let $\Psi : \cantorspace \rightarrow X$ be a reduction witnessing $= \leq E$. Let $\Phi : \cantorspace \rightarrow X \slash E$ be defined by $\Phi(x) = [\Psi(x)]_E$. $\Phi$ is an injective function.

Let $\tilde L = (X \slash E) \setminus \Phi[\cantorspace] = \bigcap_{x \in \cantorspace} (X \slash E) \setminus \{\Phi(x)\}$. $\tilde L \in U$ since $U$ is both nonprincipal and $\bbR^+$-complete.

Let $L = \bigcup \tilde L$. $L$ must be uncountable. Hence $L$ is in bijection with $\cantorspace$. Define $f : L \rightarrow (X \slash E)$ by 
$$f(x) = (X \slash E) \setminus \{[x]_E, \Phi(x)\}$$
To show $f$ is injective, it suffices to show that the map on $L$ defined by $x \mapsto \{[x]_E, \Phi(x)\}$ is injective: Suppose $x \neq y$ and $\{[x]_E, [\Psi(x)]_E\} = \{[y]_E, [\Psi(y)]_E\}$. Since $\Psi$ is a reduction, $\neg(\Psi(x) \ E \ \Psi(y))$. Therefore, one must have that $x \ E \ \Psi(y)$. This is impossible since $[x]_E \in (X \slash E) \setminus \Phi[\cantorspace]$. This shows $f$ is injective. 

Since for all $x \in L$, $[x]_E \notin f(x)$, $\tilde L \cap \bigcap_{x \in L} f(x) = \emptyset$. Since $\tilde L \in U$, $\bigcap_{x \in L} f(x) \notin U$. $U$ is not $\bbR^+$-complete. Contradiction.
\end{proof}

\Begin{fact}{ADR no R+ complete non-principal on quotient non-wellorderable}
$(\mathsf{ZF + AD_\bbR} \text{ or } \mathsf{ZF + AD^+ + V = L(\mathscr{P}(\reals)})$ Let $X$ be a Polish space and $E$ be an equivalence relation on $\cantorspace$. If $\cantorspace \slash E$ is not well-ordered, then there is no $\bbR^+$-complete nonprincipal ultrafilter on $X \slash E$. 
\end{fact}

\begin{proof}
Under $\mathsf{ZF + AD_\bbR}$, results of Woodin and Martin show that every set of reals is $\kappa$-Suslin for some $\kappa < \Theta$. So the complement of $E$ is $\kappa$-Suslin for some $\kappa < \Theta$. In $\mathsf{ZF + AD}$, \cite{Equivalence-Relations-Projective-and-Beyond} showed that if the complement of $E$ is $\kappa$-Suslin, then either the identity reduces into $E$ or $\cantorspace \slash E$ is in bijection with a cardinal less than or equal to $\kappa$ (and hence can be well-ordered). 

Under $\mathsf{ZF + AD^+ + V = L(\mathscr{P}(\reals)}$, \cite{A-Trichotomy-Theorem-in-Natural} Theorem 1.4 (along with \cite{A-Trichotomy-Theorem-in-Natural} Corollary 3.2) states that for any set $X$, either $X$ is wellordered or $\reals$ injects into $X$. 

In either case, the result now follows from \ref{no R+ complete ultrafilter}. 
\end{proof}

\section{Conclusion}\label{Conclusion}
This section includes some questions. 

\Begin{question}{ACR jonsson map for R}
Under $\mathsf{ZF} + \neg\mathrm{AC}_\omega^\bbR$, can there be $\omega$-J\'onsson functions for $\cantorspace$?

In particular, is there an $\omega$-J\'onsson function for $\cantorspace$ in the Cohen-Halpern-L\'evy model $H$ (see Question \ref{ACR and omega jonsson map})?
\end{question}

\Begin{question}{E2 jonsson surjectivity question}
It was shown that $E_2$ does not have the $2$-Mycielski property. An interesting question would be: what is the relation between the $n$-Mycielski property, $n$-J\'onsson property, and the surjectivity properties in dimension $n$ for $E_2$?

In particular, does the 2-dimensional version of the results in Section \ref{Surjectivity and Continuity Aspects of E2} hold?

Does $\cantorspace \slash E_2$ have the 2-J\'onsson property, 3-J\'onsson property, or full J\'onsson property?
\end{question}

\Begin{question}{other equivalence relation Jonsson}
For $E_1$ and $E_3$, this paper only considers the Mycielski property. One can ask about some of the other properties of $E_1$ or $E_3$ which had been studied for $E_0$ and $E_2$. For example: 

Does $\pcantorspace \slash E_1$ or $\pcantorspace \slash E_3$ have the J\'onsson property? 
\end{question}

\Begin{question}{no other set jonsson question}
Assuming determinacy, if $R \slash E_0$ injects into a set $X$, can $X$ have the J\'onsson property? More specifically, if $E$ is a $\borel$ equivalence relation on $\reals$ so that $E_0 \leq_\borel E$, can $\reals \slash E$ have the J\'onsson property?
\end{question}

\bibliographystyle{amsplain}
\bibliography{references}
\end{document}